\documentclass[11pt]{article}
\pdfoutput=1
\usepackage{enumitem}
\usepackage{setspace}
\usepackage{xcolor}
\usepackage{colortbl, multirow}
\usepackage[
bookmarksopen,
bookmarksdepth=2,
breaklinks=true
colorlinks=true,
urlcolor=blue]{hyperref}

\usepackage{fancyhdr}
\fancypagestyle{firstpage}{
\fancyhf{}
}

\usepackage{amsmath}
\usepackage{amssymb}
\usepackage{amsfonts}
\usepackage{xfrac}
\usepackage{graphicx}
\usepackage[linesnumbered,ruled,vlined]{algorithm2e}
\let\oldnl\nl% Store \nl in \oldnl
\newcommand{\nonl}{\renewcommand{\nl}{\let\nl\oldnl}}% Remove line number for one line

\usepackage{pgf}
\usepackage{tikz}
\usetikzlibrary{arrows,automata}
\usepackage[latin1]{inputenc}

\newtheorem{theorem}{Theorem}
\newtheorem{proposition}[theorem]{Proposition}
\newtheorem{lemma}[theorem]{Lemma}
\newtheorem{corollary}[theorem]{Corollary}
\newtheorem{definition}{Definition}

\newtheorem{conjecture}{Conjecture}

\begin{document}
%\title{Extended Bottleneck Division and Assignment Problems and An Optimal Solution for the Muffin Problem}
\title{An Optimal Solution for the Muffin Problem}
\author{Richard E. Chatwin
\\
%Zeta Global\\
%rchatwin1@gmail.com
%\date{October 23, 2016
}
\maketitle
\thispagestyle{firstpage}

\rule{417pt}{.7pt}

\begin{abstract}
The muffin problem asks us to divide $m$ muffins into pieces and assign each of those pieces to one of $s$ students so that the sizes of the pieces assigned to each student total $m/s$, with the objective being to maximize the size of the smallest piece in the solution. Muffin problems are a special type of variant of extended bottleneck transportation problem in which the transportation time is simply the quantity transported between any source and sink and the objective is to maximize the minimum transportation time. Of particular interest are Three Matrix Division and Assignment Problems (3M-DAP), for which all sources have the same supply and sinks are divided into two subsets having the same demand within each subset. Muffin problems are 3M-DAP in which all sinks have the same demand. We present a recursive algorithm for solving any 3M-DAP, and hence any muffin problem, and demonstrate that it always produces an optimal solution. The nature of the recursive algorithm allows us to identify interesting relationships between families of such problems.
%Problems with $m < s$ can easily be mapped to problems with $m > s$ so focus on the latter. Then an optimal solution has smallest piece $\ge 1/3$. Investigate a class of problems in which each muffin is divided into $t$ pieces and each student receives either $u$ pieces that total $x_u$ or $v$ pieces that total $x_v$. We present a recursive algorithm for solving any problem in this class, in which the problem is either solved directly or is reduced to a smaller problem with the same optimal value. A muffin problem lies
\end{abstract}

%\newpage
%  
%\tableofcontents

%\newpage

\section{Introduction}

We investigate in detail the muffin assignment problem in which we are asked to divide a given number $m$ of muffins into pieces and assign each piece to one of $s$ students so that each student receives the same total assignment $x = m/s$. Our objective is to maximize the size of the smallest piece of muffin in the solution. This is a variant of an extended bottleneck transportation problem (EBTP).

The bottleneck transportation problem (BTP) (Garfinkel and Rao, 1971) is a transportation problem for which the objective is to minimize the maximum time for transport of all goods in the network. Typically, the transport time from a source to a sink is considered fixed, independent of the quantity transported (provided the quantity is positive). The extended bottleneck transportation (EBTP) (Dai et al., 1996) removes this restriction, allowing the transportation time to be an affine function of the quantity transported. 

We focus on a particular set of variants (EBTP$^{\prime}$) of the EBTP in which the transportation time is simply the quantity transported and the objective is to maximize the minimum transportation time. Such problems may be formulated as:
\begin{alignat*}{3}
  \text{(EBTP$^{\prime}$):} \quad &&\max \min_{\{(i, j) \, | \, x_{ij} \, > \, 0\}} x_{ij} &&\\
  \text{subject to} \quad &&\sum_{j \in D} x_{ij} = s_i &\quad \quad i \in S &\\
  &&\sum_{i \in S} x_{ij} = d_j &\quad \quad j \in D &\\
  &&x_{ij} \ge 0 &\quad \quad i \in S, j \in D, &
\end{alignat*} 
where $S = \{1, 2, \ldots, m\}$ is the set of sources, $s_i$ is the supply at source $i$, $D = \{1, 2, \ldots, n\}$ is the set of sinks, $d_j$ is the demand at sink $j$, $\sum_{i \in S} s_i = \sum_{j \in D} d_j$, and $x_{ij}$ is the quantity transported from source $i$ to sink $j$.

A muffin problem is an EBTP$^{\prime}$ in which the supplies at all sources are identical and the demands at all sinks are identical. This regularity in the problem formulation is one critical element that admits the derivation of a solution algorithm and proof of its optimality.

For given numbers $m$ of muffins and $s$ of students, let $f(m,s)$ represent the optimal solution to the problem, i.e., the largest size of the smallest piece of muffin in any division of the muffins into pieces and assignment of the muffin pieces to the students.

\subsection{History of the Muffin Problem}

The muffin problem was created by recreational mathematician Alan Frank in 2009. He shared it with a private math email group, whose members established the following results, which we state without proof.
\begin{theorem} \label{thm:elem_results} Let $k, m, s \in \mathbb{N}$.
\begin{enumerate}
\item $f(m, s)$ exists, is rational, and is computable;
\item $f(m,1) = 1$ for any $m \ge 1$;
\item $f(ks, s) = 1$ for any $k \ge 1, s \ge 1$;
\item $f((2k+1)s/2,s) = \frac{1}{2}$ for any $k \ge 1, s \ge 2$ and even; and
\item $f(m, s) = (m/s)f(s,m)$ for any $m < s$.
\end{enumerate}
\end{theorem}
Veit Elser was the first to prove part 1 by showing that $f(m, s)$ is the solution to a mixed integer linear program with rational coefficients (versions of this result can be found as Theorem 11 in Cui et al. 2018 and as Theorem C.10 in Gasarch et al. 2020). An alternative proof that uses elementary topology was developed by Caleb Stanford (this approach shows that $f(m, s)$ exists and is rational but does not establish that is computable; a write-up can be found as Theorem C.15 in Gasarch et al. 2020). Part 5 was first proved by Erich Friedman (see Theorem 3 in Cui et al. 2018 and Theorem 2.9 in Gasarch et al. 2020).

Because Theorem \ref{thm:elem_results} addresses the problem when $s = 1$ and $2$, when $2x$ is integral, and when $m = s$, and also shows that the solution for any problem with $m < s$ can be derived from a corresponding problem with $m > s$, in the sequel we focus on $m > s > 2$, with $2x$ non-integral.

The group also made the following two conjectures:
\begin{conjecture} \label{conj:f(km,ks)=f(m,s)} For all $k, m, s \in \mathbb{N}$, $f(km, ks) = f(m, s)$.
\end{conjecture}

\begin{conjecture} \label{f>=1/3_conj} For all $m > s \ge 2$, $f(m, s) \ge 1/3$.
\end{conjecture}
Both of these conjectures are indeed true. We provide the first proof of Conjecture \ref{conj:f(km,ks)=f(m,s)} in \S \ref{conj1_proof}. The second conjecture has been proved by Cui et al. (2018) (see Theorem 20 therein). We will provide a simpler proof that integrates into a broader solution for the muffin problem.

The muffin problem was described by Jeremy Copeland in The New York Times Numberplay Online Blog (Antonick 2013) and appeared in a booklet of the Julia Robinson Mathematics Festival containing a sample of mathematical puzzles (Blachman 2016). This latter piqued the interest of William Gasarch, a professor at the University of Maryland, who in concert with several of his students has subsequently developed a suite of algorithms for solving the muffin problem. Some of their preliminary work is described in Ciu et al. (2018); a complete description is contained in their book (Gasarch, et al. 2020). Gasarch and his colleagues conjecture that their suite of algorithms may optimally solve any muffin problem but they have not yet been able to prove this. 

Gasarch el al. (2020) also includes a description of an algorithm for solving the muffin problem developed in 2010 by Scott Huddleston, a member of the private math email group. Though described very differently, Huddleston's algorithm is essentially identical to the algorithm we present in this paper. Section \ref{sec:scott} describes Huddleston's algorithm in the terminology of the present paper and highlights the situations in which the two algorithms employ different strategies. Two major advantages of the exposition herein, and in particular, of placing the muffin problem in the context of the broader class of extended bottleneck transportation problems, are that: (i) we are able to motivate an intuitive derivation of the algorithm; and (ii) we are able to leverage that intuition to prove the optimality of the algorithm.

\subsection{Preliminaries}

A solution with $f(m, s) > 1/3$ must have each muffin divided into at most two pieces. In fact, we can assume that all muffins are divided into two pieces because: (a) $x$ is non-integral so some muffin must be divided into two pieces, so $f(m, s) \le 1/2$; and (b) if we have a solution in which a muffin is not divided, then that whole muffin must be assigned to one student, so we can achieve an equivalent solution by dividing that muffin into two halves and assigning the two halves to the same student. 

\begin{definition} A \emph{supply-constrained muffin problem} is a muffin problem in which each muffin must be divided into exactly two pieces. Let $f(m, s, 2)$ represent the optimal solution to such a problem.
\end{definition}
Write
\[
  n = \left \lfloor \frac{2m}{s} \right \rfloor = \left \lfloor 2x \right \rfloor.
\]
If each muffin is divided into two pieces, then there are $2m$ pieces to be assigned. Then there is some student who receives at most $n$ pieces and some student who receives at least $n+1$ pieces. It is natural to make the following conjecture:

\begin{conjecture} \label{n_or_n+1_conj} For a supply-constrained muffin problem there exists an optimal solution in which all students receive either $n$ or $n + 1$ pieces.
\end{conjecture}

\begin{definition} A \emph{fully-constrained muffin problem} is a muffin problem in which each muffin is to be divided into exactly two pieces and each student must receive either $n$ or $n + 1$ pieces. Let $f(m, s, 2, n)$ represent the optimal solution to such a problem.
\end{definition}
For now we will assume Conjectures \ref{f>=1/3_conj} and \ref{n_or_n+1_conj} and restrict attention to fully-constrained muffin problems. These conjectures imply that:
\begin{enumerate}
\item $f(m, s) = \max\{1/3, f(m, s, 2)\}$
\item $f(m, s, 2) = f(m, s, 2, n)$.
\end{enumerate}

\section{Three-Matrix Division and Assignment Problems and Families Thereof}

A fully-constrained muffin problem requires dividing each supply (muffin) into two pieces and assigning the pieces so that each demand (student) receives either $n$ or $n+1$ pieces. Consider applying similar constraints to a more general EBTP$^{\prime}$. 

\begin{definition} A \emph{division and assignment problem (DAP)} is an EBTP\,$^{\prime}$ problem in which each source $i \in S$ must divided into $t_i \le n$ pieces and each sink $j \in D$ must be assigned $u_j \le m$ pieces, where $\sum_{i \in S} t_i = \sum_{j \in D} u_j$. 
\end{definition}

We can represent a DAP as $(\{\boldsymbol{t}_i : i \in S\}; \{\boldsymbol{u}_j : j \in D\})$, where the vector $\boldsymbol{t}_i$ represents the $i^{th}$ source and has length $t_i$ and the vector $\boldsymbol{u}_j$ represents the $j^{th}$ sink and has length $u_j$. Adapting notation, the supplies $\{x_{t_i} : i \in S\}$ and the demands $\{x_{u_j} : j \in D\}$ are implied and must satisfy $\sum_{i \in S} x_{t_i} = \sum_{j \in D} x_{u_j}$. The problem is to fill all the vectors so that
\begin{itemize}
\item the sum of the elements in $\boldsymbol{t}_i$ is $x_{t_i}$ for each $i \in S$;
\item the sum of the elements in $\boldsymbol{u}_j$ is $x_{u_j}$ for each $j \in D$;
\item the multiset of the elements in the $\{\boldsymbol{t}_i : i \in S\}$ is the same as the multiset of elements in the $\{\boldsymbol{u}_j : j \in D\}$;
\item the size of the smallest element is maximized.
\end{itemize}

When two or more sources have the same supply and must be divided into the same number of pieces, we can combine their vectors together to form a matrix and likewise for sinks, and alternatively represent a DAP as $(\{T_k\}; \{U_l\})$, where matrix $T_k$ has dimensions $s_{t_k} \times t_k$, each row of $T_k$ must sum to $x_{t_k}$, matrix $U_l$ has dimensions $s_{u_l} \times u_l$, each row of $U_l$ must sum to $x_{u_l}$, $\sum_k s_{t_k}t_k = \sum_l s_{u_l}u_l$, and $\sum_k s_{t_k}x_{t_k} = \sum_l s_{u_l}x_{u_l}$. The problem is to fill all the matrices so that
\begin{itemize}
\item the sum of the elements in each row of $T_k$ is $x_{t_k}$ for each $k$;
\item the sum of the elements in each row of $U_l$ is $x_{u_l}$ for each $l$;
\item the multiset of the elements in the $\{T_k\}$ is the same as the multiset of elements in the $\{U_l\}$;
\item the size of the smallest element is maximized.
\end{itemize}

Finally, inspired by the fully-constrained muffin problems, we focus on those DAPs that can be represented by three matrices, one for the supplies and two for the demands.

\begin{definition} A \emph{three-matrix division and assignment problem (3M-DAP)} $P = (T; U, V)$ takes as given:
\begin{itemize}
\item The matrix dimensions: $(s_t, t)$, $(s_u, u)$, and $(s_v, v)$. Let $n_t = s_tt$, $n_u = s_uu$, and $n_v = s_vv$;
\item The required rowsums for each matrix: $x_t$, $x_u$, and $x_v$, all rational;
\item The following requirements:
\[
\begin{split}
  t \ge 2, u \ge 2, \text{ and } v &\ge 1,\\
  \text{if } v = 1, \text{ then } n_v &\le (t - 2)s_t,\\
  s_t > 0, s_u > 0, \text{ and } s_v &\ge 0,\\
  n_u + n_v &= n_t,\\
  s_ux_u + s_vx_v &= s_tx_t, \quad \text{and}\\
  \frac{x_u}{u} &< \frac{x_v}{v}.
\end{split}
\]
\end{itemize}
Then
\[
  \frac{s_ux_u + s_vx_v}{s_uu + s_vv} = \frac{s_tx_t}{s_tt},
\]
so we must have
\begin{alignat*}{2}
  \frac{x_u}{u} &< \frac{x_t}{t} < \frac{x_v}{v}, &\quad \text{if } s_v &> 0\\
  \frac{x_u}{u} &= \frac{x_t}{t} < \frac{x_v}{v}, & \text{if } s_v &= 0.
\end{alignat*}
A \emph{feasible solution to a 3M-DAP} selects the values of the elements of each matrix so that the required rowsums are met and so that the multiset of elements of $T$ is the same is the multiset of elements of $U$ and $V$. The \emph{objective for a 3M-DAP} is to find a feasible solution with maximum size of the smallest element in $T$. Write $f(P) = f(T; U, V)$ for the optimal value.
\end{definition}

For any 3M-DAP specify the parameters $\lambda = x_u - x_v$ and $\gamma = x_t/t$. Note that we must have $\lambda < \gamma(u - v)$. 

\begin{definition} A \emph{3M-DAP parameter set} is a set of parameters $(t, u, v, \lambda, \gamma)$ satisfying $t, u, v$ integer, $\lambda, \gamma$ rational, $t \ge 2$, $u \ge 2$, $v \ge 1$, and $\lambda < \gamma(u - v)$.
\end{definition}

Focus on problems of specific interest is enabled by specifying the following classes of problem parameters.

\begin{definition} Classes of 3M-DAP parameter sets are defined by specifying additional constraints on the problem parameters:
\begin{enumerate}
\item \emph{$u$-weighted}: $2 \le v \le u$
\item \emph{weakly $u$-weighted}: $2 \le v \le u \le (v - 1)t$
\item \emph{strongly $u$-weighted}: $2 \le v < u$
\item \emph{$v$-weighted}: $2 \le u < v$
\item \emph{$v1$}: $v = 1$.
\end{enumerate}
A \emph{($u$-weighted, weakly $u$-weighted, strongly $u$-weighted, $v$-weighted, $v1$) 3M-DAP} is a 3M-DAP for which the parameter set is a ($u$-weighted, weakly $u$-weighted, strongly $u$-weighted, $v$-weighted, $v1$) 3M-DAP parameter set.
\end{definition}

\begin{definition} Given a 3M-DAP parameter set $(t, u, v, \lambda, \gamma)$ and a positive integer $s$, the \emph{3M-DAP family} $F_s(t, u, v, \lambda, \gamma)$ is the set of 3M-DAPs $(T; U, V)$ with the given number of columns ($t$ for $T$, $u$ for $U$, and $v$ for $V$), with $x_u - x_v = \lambda$ and $x_t/t = \gamma$, and with $s_u + s_v = s$. 
\end{definition}

Any 3M-DAP family can be parametrized by $x = x_u$; varying $x$ will lead to corresponding variance in $s_u$ and $s_v$. 

\begin{lemma} \label{lem:x_vals} The 3M-DAP family $F_s(t, u, v, \lambda, \gamma)$ contains at most $s$ problems. For each such problem $\lambda + \gamma v < x \le \gamma u$.
\end{lemma}
\textbf{Proof.} Because we require that $s_u > 0$, $s_u$ can take on one of the $s$ values $1, 2, \dots, s$. Each such corresponds to a different problem in the family only when $s_t = (s_uu + s_vv)/t$ is integer. Now,
\[
  s_ux_u + s_vx_v = s_tx_t \quad \Leftrightarrow \quad sx - \lambda s_v = \gamma n_t \quad \Leftrightarrow \quad x = \frac{1}{s}(\lambda s_v + \gamma n_t).
\]
Further, $n_t = n_u + n_v = us_u + vs_v = us - (u - v)s_v$, so
\begin{equation} \label{eq:x}
  x = \gamma u - [\gamma (u - v) - \lambda] \frac{s_v}{s}.
\end{equation}
Because $\lambda < \gamma(u - v)$ and $0 \le s_v < s$, we have $\lambda + \gamma v < x \le \gamma u$. 
%Not all possible rational values of $x$ in the interval are contained in the family as only integral values of $s_v$ for which also $s_t = (n_u + n_v)/t$ is integral are possible.
\hspace{\stretch{1}}$\blacksquare$

\begin{definition} Given a 3M-DAP parameter set $(t, u, v, \lambda, \gamma)$, the \emph{extended 3M-DAP family}\\ $F(t, u, v, \lambda, \gamma)$ is the union over all positive integer values of $s$ of the 3M-DAP families\\ $F_s(t, u, v, \lambda, \gamma)$.
\end{definition}

Now we can write $s_t$, $s_u$, and $s_v$ in terms of $t$, $u$, $v$, $\lambda$, $\gamma$, $s$, and $x$, using that $n_u + n_v = n_t$ and $s_u + s_v = s$:
\begin{align}
  s_v &= \frac{\gamma u - x}{\gamma (u - v) - \lambda} s \label{s_v}\\
  s_u &= \frac{x - \lambda - \gamma v}{\gamma (u - v) - \lambda} s \label{s_u}\\
  s_t &= \frac{n_t}{t} = \frac{n_u + n_v}{t} = \frac{(u - v)x - \lambda u}{\gamma (u - v) - \lambda} \frac{s}{t}. \label{s_t}
\end{align}

\begin{proposition} \label{prop:rational} The values of $x$ for the problems in the extended 3M-DAP family $F(t, u, v, \lambda, \gamma)$ are all possible rational numbers in the interval $(\lambda + \gamma v, \gamma u]$.
\end{proposition}
\textbf{Proof.} For any problem in the extended family, $x$ is given by equation \eqref{eq:x} so is rational. By Lemma \ref{lem:x_vals}, $x$ lies in the given interval.

For the converse, $\lambda$ and $\gamma$ are rational, so there are integers $c, d, e,$ and $f$ with $d$ and $f$ positive such that $\gamma = c/d$, and $\lambda = e/f$. Suppose $a$ and $b$ are integer with $b$ positive and $\lambda + \gamma v < a/b \le \gamma u$. We want to show that there exists a member of the extended family with $x = a/b$. Set
\[
\begin{split}
  s &= [\gamma(u - v) - \lambda]bdft\\
  s_v &= \left(\gamma u - \frac{a}{b}\right)bdft\\
  s_u &= \left(\frac{a}{b} - \lambda - \gamma v\right)bdft\\
  s_t &= \left[(u - v)\frac{a}{b} - \lambda u\right]bdf.
\end{split}
\]
Then each is an integer with $s_v$ nonnegative and the other three positive. From \eqref{eq:x} we obtain $x = a/b$, as desired.
\hspace{\stretch{1}}$\blacksquare$

\begin{definition} The set $\{t, \{u_i\}_{i \ge k}, \{v_j\}_{j \ge k}, \lambda, \gamma\}$ defines a \emph{complete 3M-DAP space} if for all $i \ge k$:
\begin{itemize}
\item $(t, u_i, v_i, \lambda, \gamma)$ is a 3M-DAP parameter set;
\item $v_{i + 1} = u_i - \lambda/\gamma$;
\item $u_{i + 1} > u_i$.
\end{itemize}
\end{definition}

\begin{definition} Given a complete 3M-DAP space $\{t, \{u_i\}_{i \ge k}, \{v_j\}_{j \ge k}, \lambda, \gamma\}$, the \emph{complete 3M-DAP family} $F(t, \{u_i\}_{i \ge k}, \{v_j\}_{j \ge k}, \lambda, \gamma)$ is the union over all $i \ge k$ of the extended 3M-DAP families $(t, u_i, v_i, \lambda, \gamma)$.
\end{definition}

\subsection{Fully-constrained muffin problems as three-matrix division and assignment problems}

We now describe how fully-constrained muffin problems fit within the class of 3M-DAPs. A fully-constrained muffin problem is a 3M-DAP in which the matrix $T$ represents the muffins: $t = 2$, $s_t = m$, and $x_t = 1$; and the matrices $U$ and $V$ represent the students: $u = n + 1$, $v = n$, $s_u + s_v = s$, and $x_u = x_v = x$. To ensure that $n_u + n_v = n_t$, we must have 
\begin{align*}
  s_u &= 2m - ns,\\
  s_v &= (n + 1)s - 2m.
\end{align*}
A fully-constrained muffin problem has $\lambda = 0$ and $\gamma = 1/2$.

\begin{definition} An \emph{extended fully-constrained-muffin-problem family of order $n$} is precisely the extended 3M-DAP family $F(2, n + 1, n, 0, 1/2)$. This extended family includes all possible rational values of $x$ in the interval $(n/2, (n + 1)/2]$.
\end{definition}

\begin{definition} The \emph{complete fully-constrained-muffin-problem family} is precisely the complete 3M-DAP family $F(2, \{n + 1\}_{n \ge 2}, \{n\}_{n \ge 2}, 0, 1/2)$. This family contains all possible fully-constrained muffin problems.
\end{definition}

\subsection{Equivalence and Standard Form}

\begin{definition} Two 3M-DAPs are \emph{equivalent} if the dimensions of the matrices are the same in the two problems and there exists an increasing, linear \emph{equivalence function} that maps the elements in a solution to the first problem to the elements in a solution to the second problem.
\end{definition} 

\begin{definition} A strongly $u$-weighted 3M-DAP is in \emph{standard form} if $\lambda = 0$ and $\gamma = 1/2$. Then $x_u = x_v = x$ and $x_t = t/2$. An extended 3M-DAP family in standard form has $x \in (v/2, u/2]$. Such problems are of interest because fully-constrained muffin problems have this form.
\end{definition}

\begin{theorem} \label{standard_form_equiv} Any strongly $u$-weighted 3M-DAP $P = (T; U, V)$ is equivalent to a 3M-DAP $\hat{P} = (\hat{T}; \hat{U}, \hat{V})$ in standard form. The equivalence function $h(\cdot)$ is:
\begin{equation} \label{eqfunc}
  h(y) = \frac{1}{2} \frac{(u - v)y - (x_u - x_v)}{(u - v)x_t/t - (x_u - x_v)} = \frac{1}{2} \frac{(u - v)y - \lambda}{\gamma(u - v) - \lambda}.
\end{equation}
\end{theorem}
\textbf{Proof.} The equivalence function $h(\cdot)$ maps the matrices $T$, $U$, and $V$ to the matrices $\hat{T}$, $\hat{U}$, and $\hat{V}$, respectively. First, if $(y_1, \ldots, y_u)$ is a row of $U$ and $(z_1, \ldots, z_v)$ is a row of $V$, then $y_1 + \cdots + y_u = x_u$, $z_1 + \cdots + z_v = x_v$, $(h(y_1), \ldots, h(y_u))$ is a row of $\hat{U}$, and $(h(z_1), \ldots, h(z_v))$ is a row of $\hat{V}$. Then
\[
\begin{split}
  h(y_1) + \cdots + h(y_u) &= \frac{1}{2} \frac{(u - v)x_u - u(x_u - x_v)}{(u - v)x_t/t - (x_u - x_v)}\\
  &= \frac{1}{2} \frac{ux_v - vx_u}{(u - v)x_t/t - (x_u - x_v)}\\
  &= \frac{1}{2} \frac{(u - v)x_v - v(x_u - x_v)}{(u - v)x_t/t - (x_u - x_v)}\\
  &= h(z_1) + \cdots + h(z_v).
\end{split}
\]
Second, if $(w_1, \ldots, w_t)$ is a row of $T$, then $w_1 + \cdots + w_t = x_t$, and $(h(w_1), \ldots, h(w_t))$ is a row of $\hat{T}$. Then
\[
  h(w_1) + \cdots + h(w_t) = \frac{1}{2} \frac{(u - v)x_t - t(x_u - x_v)}{(u - v)x_t/t - (x_u - x_v)} = \frac{t}{2}.
\]
Then the mapped problem is a 3M-DAP and has $\lambda = 0$ and $\gamma = 1/2$ and so is in standard form, as claimed.
%Because we must map $x_t/t$ to $1/2$, look for an equivalence function of the form
%\[
%  h(y) = \frac{1}{2} + \left(y - \frac{x_t}{t}\right)\theta.
%\]
%The other condition is that if $y_1 + \cdots + y_u = x_u$ and $z_1 + \cdots + z_v = x_v$, then $h(y_1) + \cdots + h(y_u) = h(z_1) + \cdots + h(z_v)$. To ensure this, we need
%\[
%\begin{split}
%  h(y_1) + \cdots + h(y_u) &= h(z_1) + \cdots + h(z_v)\\
%  \Leftrightarrow \frac{u}{2} - \left(u\frac{x_t}{t} - x_u\right)\theta &= \frac{v}{2} - \left(v\frac{x_t}{t} - x_v\right)\theta\\
%  \Leftrightarrow \frac{1}{2}t(u - v) &= \left[(u - v)x_t - t(x_u - x_v)\right]\theta\\
%  \Leftrightarrow \theta &= \frac{1}{2} \frac{t(u - v)}{(u - v)x_t - t(x_u - x_v)} = \frac{1}{2} \frac{u - v}{(u - v)x_t/t - (x_u - x_v)}.
%\end{split}
%\]
%Note that because $x_u/u < x_t/t < x_v/v$, we have $x_t/t > x_u/u > (x_u - x_v)/(u - v)$, so $(u - v)x_t > t(x_u - x_v)$, so $\theta > 0$.
%
%Now
%\[
%  2\frac{x_t}{t}\theta - 1 = \frac{x_u - x_v}{(u - v)x_t/t - (x_u - x_v)},
%\]
%so we have
%\[
%  h(y) = \frac{1}{2} \frac{(u - v)y - (x_u - x_v)}{(u - v)x_t/t - (x_u - x_v)}.
%\]
\hspace{\stretch{1}}$\blacksquare$

A 3M-DAP with $u = v$ must have $\lambda < 0$ and so cannot be in standard form. This is reflected in the form of the equivalence function in \eqref{eqfunc}: the requirement in Theorem \ref{standard_form_equiv} that $u > v$ is necessary; if $u = v$, it is not possible to transform the problem to standard form. When $u < v$, the problem can be transformed to standard form but the equivalence function will be decreasing so the smallest piece in a solution to the original problem will be mapped to the largest piece in a solution to the problem in standard form. Accordingly, in the sequel, only problems with $u > v$ will be mapped to standard form.

For later reference, note that the inverse transformation is 
\begin{equation} \label{eqfunc_inv}
  h^{-1}(z) = 2\left[\frac{x_t}{t} - \frac{x_u - x_v}{u - v}\right]z + \frac{x_u - x_v}{u - v} = 2\frac{x_t}{t}z + \frac{x_u - x_v}{u - v}(1 - 2z) = 2\gamma z + \frac{\lambda}{u - v}(1 - 2z).
\end{equation}

\section{Motivating an Algorithm for Solving 3M-DAPs}

In the following sections we will develop a recursive algorithm for optimally solving any 3M-DAP. This algorithm optimally solves any fully-constrained muffin problem because such a problem is a 3M-DAP.

An immediate upper bound for any 3M-DAP is given by the average size of an element in the U matrix: $x_u/u$. Call a 3M-DAP for which this upper bound can be achieved a \emph{0-problem}. We shall see that the set of 0-problems can be naturally divided into two, which we call \emph{types 1 and 2}. Any 0-problem can be solved directly. 

When the upper bound of $x_u/u$ is not achievable, the 3M-DAP is not a 0-problem, we reduce it to a smaller problem that is then solved recursively: if the reduced problem is a 0-problem, solve it directly, else further reduce to a yet smaller problem; continue until a 0-problem is reached. The value of the original problem will be the same as the value of the smaller problem. In this section we motivate this solution approach.

\subsection{$(T, V)$-pairs}

Consider a solution to a DAP $P = (T; \{U_1, \ldots, U_p\}, V)$.

\begin{definition}
A \emph{$(T, V)$-pair} is a pair $(A, B)$ of sets of rows of $T$ and $V$, respectively, such that each row in $A$ contains at least one element of the rows of $V$ in $B$ and collectively, the rows in $A$ contain all the elements of the rows of $V$ in $B$. Write $a = |A|$ and $b = |B|$. 
\end{definition}
Note that it is possible to have $b = 0$, so that $A$ consists entirely of elements of $U$.

\begin{definition}
A $(T, V)$-pair is called \emph{inseparable} if there do not exist $(T, V)$-pairs $(A_1, B_1)$ and $(A_2, B_2)$ such that $A_1 \cap A_2 = \emptyset$, $A = A_1 \cup A_2$, $B_1 \cap B_2 = \emptyset$, and $B = B_1 \cup B_2$. 
\end{definition}

\begin{definition}
A \emph{maximal} $(T, V)$-pair is a $(T, V)$-pair $(A, B)$ for which all the remaining elements in the rows of $A$ belong to $U$. Then $(T, V)$ is a maximal $(T, V)$-pair.
\end{definition}

\begin{lemma} \label{lem:maximal_inseparable_pairs}
Given a solution to a DAP $P = (T; \{U_1, \ldots, U_p\}, V)$, the pair $(T, V)$ can be divided into a set $S =  \{(A_1, B_1), (A_2, B_2), \ldots, (A_M, B_M)\}$ of mutually exclusive and collectively exhaustive, maximal, inseparable $(T, V)$-pairs.
\end{lemma}
\textbf{Proof.} If $(T, V)$ is inseparable, we can simply set $S = \{(T, V)\}$. Otherwise, $(T, V)$ can be separated into $(T, V)$-pairs $(A_1, B_1)$ and $(A_2, B_2)$ such that $A_1 \cap A_2 = \emptyset$, $A = A_1 \cup A_2$, $B_1 \cap B_2 = \emptyset$, and $B = B_1 \cup B_2$. Then each pair $(A_1, B_1)$ and $(A_2, B_2)$ is maximal. If one of these is separable, then we can separate it into further maximal pairs. We can repeat the process until all resulting pairs are inseparable.
\hspace{\stretch{1}}$\blacksquare$

\begin{lemma} \label{lem:inseparable_pair} 
If $(A, B)$ is an inseparable $(T, V)$-pair, then $a \le (v - 1)b + 1$.
\end{lemma}
\textbf{Proof.} Consider inserting the rows of $B$ into $A$ sequentially. Let $C_i \subset A$ be the rows of $A$ containing at least one element of the rows of $B$ after having inserted $i$ of the rows. Insert the first row; then $|C_1| \le v$. After inserting $i$ rows, because $(A, B)$ is inseparable, there must be a row in $B$ that has yet to be inserted and that has an element in a row in $C_i$. Insert such a row next into $A$. Then $|C_{i + 1}| \le |C_i| + v - 1$. Conclude that
\[
  a \le v + (b - 1)(v - 1) = (v - 1)b + 1.
\]
\hspace{\stretch{1}}$\blacksquare$

\begin{lemma} \label{lem:pairs_in_opt_soln}
Consider a DAP $P = (T; \{U_1, \ldots, U_p\}, V)$. If either $t > 2$ or $v > 2$, suppose that 
\[
  f(P) > \frac{(v - 1)x_t - x_v}{(v - 1)t - v}.
\]
Then in any optimal solution, any $(T, V)$-pair must have $a \ge (v - 1)b + 1$ and any inseparable $(T, V)$-pair must have $a = (v - 1)b + 1$.
\end{lemma}
\textbf{Proof.} Consider such a problem containing a $(T, V)$-pair $(A, B)$. If $v = 1$, the result is immediate. If $t = v = 2$, then because the number of elements in $A$ must exceed the number of elements in $B$, we must have $a > b = (v - 1)b$. Otherwise, $(v - 1)t - v > 0$. Suppose that $a \le (v - 1)b$. Then an upper bound on the minimum piece size of the solution is given by the average size of the elements in $A$ that are not in $B$, which is
\[
  \frac{ax_t - bx_v}{at - bv} = \frac{(a/b)x_t - x_v}{(a/b)t - v} \le \frac{(v - 1)x_t - x_v}{(v - 1)t - v} < f(P),
\]
where the last inequality is by assumption. Conclude that such a solution is not optimal. Because the solution is assumed optimal, we must have $a \ge (v - 1)b + 1$. The conclusion follows upon referencing the previous lemma.
\hspace{\stretch{1}}$\blacksquare$

\subsection{Conjectures}

Now focus on the solution to a 3M-DAP $P = (T; U, V)$. A maximal, inseparable $(T, V)$-pair $(A, B)$ provides an upper bound on the value of the problem $P$, namely the average value of the $U$-elements in the pair:
\[
  \frac{ax_t - bx_v}{at - bv} \le \frac{[(v - 1)b + 1]x_t - bx_v}{[(v - 1)b + 1]t - bv} = \frac{x_t + [(v - 1)x_t - x_v]b}{t + [(v - 1)t - v]b}.
\]
The left-hand term is decreasing in the ratio $b/a$; the right-hand term is decreasing in $b$. Because any solution must contain a $(T, V)$-pair with $a/b \le s_t/s_v$, and when $s_t \le (v - 1)s_v$ the constraint of Lemma \ref{lem:inseparable_pair} does not apply, we make the following conjecture.

\begin{conjecture} \label{conj:0prob_type1} \emph{0-problems of type 1}. For a 3M-DAP $P = (T; U, V)$ with $s_t \le (v - 1)s_v$, then the optimal solution consists of $(T, V)$-pairs with $a/b = s_t/s_v$ ($(T, V)$ is one such pair), and the optimal value is given by
\[
  \frac{s_tx_t - s_vx_v}{s_tt - s_vv} = \frac{x_u}{u},
\]
i.e., all values in $U$ must be $x_u/u$.
\end{conjecture}

Now suppose that $s_t > (v - 1)s_v$. Assuming $(v - 1)t - v > 0$,
\[
  \frac{x_u}{u} = \frac{s_tx_t - s_vx_v}{s_tt - s_vv} > \frac{(v - 1)x_t - x_v}{(v - 1)t - v}.
\]
This explains the assumption in the statement of Lemma \ref{lem:pairs_in_opt_soln}, which we formalize here.
\begin{conjecture} \label{NprobLB_conj} For a 3M-DAP $P = (T; U, V)$ with $(v - 1)t - v > 0$ and $s_t > (v - 1)s_v$,
\[
  f(P) > \frac{(v - 1)x_t - x_v}{(v - 1)t - v}.
\]
\end{conjecture}
In light of Conjecture \ref{NprobLB_conj} and Lemma \ref{lem:pairs_in_opt_soln}, we seek an optimal solution in which each pair $(A_i, B_i) \in S$, being inseparable, has $a_i = (v - 1)b_i + 1$.

\begin{definition}
A \emph{$b$-pair} is a maximal inseparable $(T, V)$-pair $(A, B)$ with $|B| = b$ and $|A| = (v - 1)b + 1$. It contains $b$ rows from $V$, $(v - 1)b + 1$ rows from $T$, and those rows from $T$ contain $[(v - 1)b + 1]t - bv$ elements from $U$.
\end{definition}

Define
\[
  r_v = n_u - (t - 2)s_t - (v - 2)s_v = 2[s_t - (v - 1)s_v],
\]
so that when $s_t > (v - 1)s_v$, $r_v$ is positive and even. Further define
\begin{align}
  b^* &= \frac{2s_v}{r_v} = \frac{s_v}{s_t - (v - 1)s_v} \label{eq:b*}\\
  \hat{b} &= \left \lceil \frac{2s_v}{r_v} \right \rceil. \label{eq:b'}
\end{align}
If $b^*$ is integer, a $b^*$-pair has $a^* = (v - 1)b^* + 1$ so that $a^*/b^* = s_t/s_v$. Then the average value of a $U$-element in the $b^*$-pair is $x_u/u$. Then it is natural to seek a solution with all pairs $(A_i, B_i) \in S$ having $b_i = b^*$ and we might anticipate that it is possible to achieve the $x_u/u$ upper bound. 

If $b^*$ is not integer, the set $S$ must contain at least one pair $(A_i, B_i)$ with $b_i \ge \hat{b}$ and at least one pair $(A_j, B_j)$ with $b_j \le \hat{b} - 1$. When $b > b^*$, the average value of a $U$-element in the $b$-pair is less than $x_u/u$, so this upper bound cannot be achieved. It is now natural to seek a solution with all pairs $(A_i, B_i) \in S$ having either $b_i = \hat{b}$ or $b_i = \hat{b} - 1$.

\begin{conjecture} \label{conj:b*} \emph{0-problems of type 2}. For a 3M-DAP $P = (T; U, V)$ with $s_t > (v - 1)s_v$, if $b^* = 2s_v/r_v$ is integer, then there exists an optimal solution in which all maximal inseparable pairs are $b^*$-pairs.
\end{conjecture}

\begin{conjecture} \label{conj:b'_or_b'-1} \emph{Reduced problems}. For a 3M-DAP $P = (T; U, V)$ with $s_t > (v - 1)s_v$, if $b^* = 2s_v/r_v$ is not integer, then there exists an optimal solution in which all maximal inseparable pairs are either $\hat{b}$- or $(\hat{b} - 1)$-pairs.
\end{conjecture}

\begin{lemma} \label{lem:v1}
Consider a $v1$ 3M-DAP, i.e., with $v = 1$. Then the problem is not a 0-problem of type 1 and $b^* \le t - 2$.
\end{lemma}
\textbf{Proof.} The first claim follows on observing that clearly $s_t > (v - 1)s_v$. For the second claim, $b^* = s_v/s_t$ and by assumption for such a problem $s_v = n_v \le (t - 2)s_t$.
\hspace{\stretch{1}}$\blacksquare$

\subsection{Sketch of Solution Strategy for 0-problems of type 1}

In the next few sections we focus on problems with $s_t \le (v - 1)s_v$ and demonstrate that the $x_u/u$ upper bound can be achieved. As an example of the approach, suppose $n_u < (t-2)s_t$ and $v \ge 2$. Set all elements of $U$ to $x_u/u$ and insert these into $T$ so that the rows of $T$ have similar numbers of elements filled. Let $q_t \ge 0$ and $0 \le r_t < s_t$ be such that
\[
  n_u = q_ts_t + r_t.
\]
Then $r_t$ rows of $T$ will have $q_t + 1$ elements set to $x_u/u$ and $s_t - r_t$ rows of $T$ will have $q_t$ elements set to $x_u/u$. Because $q_t < t - 2$, the remaining problem looks familiar: it is a weakly $u$-weighted 3M-DAP $(T^{\prime}; U^{\prime}, V^{\prime})$ where $T^{\prime}$ is the matrix $V$ with required rowsums $x_v$, $U^{\prime}$ is the empty $(s_t - r_t) \times (t - q_t)$ submatrix of $T$ with required rowsums $x_t - q_tx_u/u$, and $V^{\prime}$ is the empty $r_t \times (t - q_t - 1)$ submatrix of $T$ with required rowsums $x_t - (q_t + 1)x_u/u$. One key element to our solution strategy is to demonstrate lower bounds on 3M-DAPs; in this case demonstrating a lower bound on this latter problem that is larger than $x_u/u$ establishes that the upper bound of $x_u/u$ is achievable for the original problem.

We begin the process of demonstrating lower bounds in \S \ref{sec:greedy_algo}, where we describe a greedy algorithm for solving 3M-DAPs with $t=u=v=2$. This algorithm establishes a lower bound on the optimal solution to such problems. The algorithm and corresponding lower bound is used in \S \ref{sec:1/3} in an algorithm that produces a solution with value of 1/3 for any muffin problem, thus establishing Conjecture \ref{f>=1/3_conj}.

We will later be able to establish Conjecture \ref{n_or_n+1_conj}. Then any muffin problem can be solved by first applying the recursive algorithm to solve the fully-constrained version of the problem. If the solution has value greater than or equal to 1/3, the solution is optimal for the unconstrained muffin problem. Otherwise, apply the algorithm of \S \ref{sec:1/3} to obtain a solution with value 1/3. In fact, we will be able to identify a priori which muffin problems have optimal solutions with value equal to 1/3, and hence can apply the appropriate algorithm.

\section{A Greedy Algorithm for 3M-DAPs with $t=u=v=2$} \label{sec:greedy_algo}

For a 3M-DAP with $t=u=v=2$, we have $s_t = s_u + s_v$ so the number of rows in $T$ is equal to the sum of numbers of rows in $U$ and $V$. We must have $x_u < x_v \Leftrightarrow \lambda < 0$ and $x_t < x_v \Leftrightarrow x_u > x_t + \lambda$.

Algorithm \ref{alg:greedy} is a greedy algorithm for this type of problem. Construct two vectors $\boldsymbol{y}$ and $\boldsymbol{z}$, each of length $s_t$. The first row of $T$ is $[y_1, z_{s_t}] = [x_t/2, x_t/2]$. For $2 \le i \le s_t$, the $i^{th}$ row of $T$ is $[y_i, z_{i-1}]$. Each row of $U$ and $V$ is $[y_j, z_j]$ for some $j$.

\fontsize{8}{8}
\begin{algorithm}
\SetKwInOut{Input}{input}
\SetKwInOut{Output}{output}
\Input{$T$, $U$, $V$, empty matrices, each with two columns and with numbers of rows $s_t$,\\
\quad \quad $s_u$, and $s_v = s_t - s_u$, respectively\\
$x_t$, $x_u$, $x_v$, the required row sums for each row of $T$, $U$, and $V$, respectively\\
%$T$, an empty matrix of dimensions $s_t \times 2$\\
%$x_t$, the required row sum for each row of $T$\\
%$U$, an empty matrix of dimensions $s_u \times 2$\\
%$x_u$, the required row sum for each row of $U$\\
%$V$, an empty matrix of dimensions $s_v \times 2$, where $s_v = s_t - s_u$\\
%$x_v$, the required row sum for each row of $V$\\
}
\Output{the filled matrices $T$, $U$, $V$\\
the multiset of entries in $T$ is the union of the multiset of entries in $U$ and $V$\\
}
\BlankLine

\Begin(\textbf{part 1}){
divide the first row of $T$, $[x_t/2, x_t/2]$, so that $y_1 = z_{s_t} = x_t/2$\\

\For{$1 \le i \le s_t$}{
	compute $z_i^u = x_u - y_i$ and $z_i^v = x_v - y_i$\\
	\uIf{$z_i^u > (x_u + x_t - x_v)/2 = (x_t + \lambda)/2$}{
		set $z_i = z_i^u$ and assign $[y_i, z_i]$ to the next empty row of $U$\\
	}
	\uElse{
		set $z_i = z_i^v$ and assign $[y_i, z_i]$ to the next empty row of $V$\\
	}
	\uIf{$i < s_t$}{
		compute $y_{i+1} = x_t - z_i$\\
		assign $[y_{i+1}, z_i]$ to the next empty row of $T$\\
	}
}
}

\Begin(\textbf{part 2}){
compute $y_{\min} = \min_i\{y_i\}$ and $z_{\min} = \min_i\{z_i\}$\\
compute $\varepsilon = (z_{\min} - y_{\min})/2$\\
$y_i \leftarrow y_i + \varepsilon, \forall i$\\
$z_i \leftarrow z_i - \varepsilon, \forall i$\\
}

\caption{Greedy algorithm for 3M-DAP with $t=u=v=2$}
\label{alg:greedy}
\end{algorithm}

\normalsize

At the $i^{th}$ step of the algorithm, the piece $y_i$ from $T$ must be assigned to either $U$ or $V$. If $y_i$ is assigned to $U$, then the other piece in that row of $U$ has size $z_i = z_i^u = x_u - y_i$. 
Alternatively, if $y_i$ is assigned to $V$, then the other piece in that row of $V$ has size $z_i = z_i^v = x_v - y_i$. The greedy algorithm assigns $y_i$ so as to make $z_i$ as close to $x_t/2$ as possible. Now
\[
  z_i^u - z_i^v = x_u - x_v = \lambda < 0,
\]
so if $z_i^u > (x_u + x_t - x_v)/2 = (x_t + \lambda)/2$, then $y_i$ is assigned to $U$ and $z_i = z_i^u$, else $y_i$ is assigned to $V$ and $z_i = z_i^v$. Note that if only rows of $U$ remain, then it must be the case that $z_i^u > (x_t + \lambda)/2$, so that $y_i$ is assigned to $U$, and likewise if only rows of $V$ remain, then $z_i^u < (x_t + \lambda)/2$, and $y_i$ is assigned to $V$, as required.

The second part of the algorithm adjusts the sizes of all the pieces so that feasibility is maintained and the size of the smallest piece is increased (or remains unchanged).

The following lemma provides bounds on the sizes of the pieces in the solution given by the greedy algorithm. The lower bound in particular will be useful.

\begin{lemma} \label{greedy_t=u=v=2} A lower bound for a 3M-DAP with $t = u = v = 2$ is given by $(x_t + \lambda)/2$: Algorithm \ref{alg:greedy} produces $\boldsymbol{y}$ and $\boldsymbol{z}$ with $(x_t + \lambda)/2 < y_i, z_i < (x_t - \lambda)/2$, for all $i$. 
\end{lemma}
\textbf{Proof.} We begin by showing that part 1 produces $\boldsymbol{y}$ and $\boldsymbol{z}$ with $(x_t + \lambda)/2 \le y_i < (x_t - \lambda)/2$ and $(x_t + \lambda)/2 < z_i \le (x_t - \lambda)/2$, for all $i$. 

The proof is by induction. We certainly have $(x_t + \lambda)/2 \le y_1 < (x_t - \lambda)/2$. Suppose it is true of $y_i$. Then we have
\[
\begin{split}
  x_t + \lambda &< x_u < x_t\\
  \Leftrightarrow (x_t + 3\lambda)/2 = x_t + \lambda - (x_t - \lambda)/2 < x_t + \lambda - y_i &< x_u - y_i < x_t - y_i \le (x_t - \lambda)/2\\
  \Leftrightarrow (x_t + 3\lambda)/2 &< z_i^u < (x_t - \lambda)/2.
\end{split}
\]
We assign $y_i$ to $U$ when $z_i^u > (x_t + \lambda)/2$, in which case $z_i = z_i^u$ and $(x_t + \lambda)/2 < z_i < (x_t - \lambda)/2$. Alternatively, we assign $y_i$ to $V$ when $z_i^u \le (x_t + \lambda)/2$, in which case $z_i = z_i^v$ and $(x_t + \lambda)/2 < z_i \le (x_t - \lambda)/2$. Then $y_{i+1} = x_t - z_i$ and we have $(x_t + \lambda)/2 \le y_{i+1} < (x_t - \lambda)/2$.

Finally, in part 2, the size of the smallest piece is only unchanged if $z_{\min} = y_{\min}$, but then both are larger than $(x_t + \lambda)/2$. Otherwise, the size of the smallest piece strictly increases and so must be larger than $(x_t + \lambda)/2$.
\hspace{\stretch{1}}$\blacksquare$

It is in fact possible to show, as we do in Appendix \ref{App_Greedy_Opt}, that Algorithm \ref{alg:greedy} produces an optimal solution for any 3M-DAP with $t = u = v = 2$ and to determine the optimal value precisely. However, the strict lower bound of Lemma \ref{greedy_t=u=v=2} will suffice in the sequel.

\begin{lemma} \label{lem:alg_greedy_complexity} Algorithm \ref{alg:greedy} has complexity $\Theta(n_t)$.
\end{lemma}
\textbf{Proof.} Algorithm \ref{alg:greedy} is initialized by computing $x_t/2$ and $(x_t + \lambda)/2 = (x_u + x_t - x_v)/2$, and by setting $y_{\min} = +\infty$ and $z_{\min} = +\infty$, for a total of 8 operations.

In part 1 of the algorithm, for each row $i$ of $T$, Algorithm \ref{alg:greedy}:
\begin{itemize}
\item computes $y_i$ (either as $x_t/2$ when $i = 1$ or as $x_t - z_{i-1}$ when $i > 1$);
\item computes $z_i^u = x_u - y_i$ and $z_i^v = x_v - y_i$;
\item compares $z_i^u$ to $(x_t + \lambda)/2$;
\item sets $z_i$ (to either $z_i^u$ or $z_i^v$);
\item assigns $[y_i, z_i]$ to the next empty row of either $U$ or $V$;
\item assigns $[y_i, z_{i-1}]$ to the next empty row of $T$;
\item compares $y_i$ to $y_{\min}$ and updates $y_{\min}$ if necessary;
\item compares $z_i$ to $z_{\min}$ and updates $z_{\min}$ if necessary.
\end{itemize}
The number of operations is at most $8n_t$.

Part 2 of the algorithm computes $\varepsilon = (z_{\min} - y_{\min})/2$, then either adds or subtracts $\varepsilon$ to or from every element in $T$, $U$, and $V$. The number of operations is $2n_t + 3$.

Then Algorithm \ref{alg:greedy} takes at most $10n_t + 11$ operations. The conclusion follows.
\hspace{\stretch{1}}$\blacksquare$

\section{$f(m,s) \ge \frac{1}{3}$} \label{sec:1/3}

Using the lower bound established in the previous section allows us to prove Conjecture 1.

\fontsize{8}{8}
\begin{algorithm}
\SetKwInOut{Input}{input}
\SetKwInOut{Output}{output}
\Input{$m$, the number of muffins\\
$s$, the number of students\\
}
\Output{the division of muffins and the assignment of pieces to students}
\BlankLine
cut each of $m - s$ muffins into three equal pieces of size $1/3$\\
\Begin(\textbf{Case 1:} $3x = k > 3$, $k$ an integer){
	cut the remaining $s$ muffins into two equal pieces of size $1/2$\\
	assign each student $(k - 3)$ pieces of size $1/3$ and two pieces of size $1/2$
}

\Begin(\textbf{Case 2:} $k < 3x < k + 1$, $k$ an integer){
	assign the $3(m - s)$ pieces of size 1/3 so that each student receives either $k - 3$ or $k - 2$ such pieces\\
	use Algorithm \ref{alg:greedy} to complete the assignment by cutting each of the remaining $s$ muffins into two pieces and assigning each student two of those pieces
}

\caption{Algorithm to achieve lower bound of 1/3 for any muffin problem}
\label{alg:1/3}
\end{algorithm}

\normalsize

\begin{lemma} \label{1/3_v1} For any muffin problem with $m > s$, $f(m,s) \ge 1/3$.
\end{lemma}
\textbf{Proof.} We demonstrate that Algorithm \ref{alg:1/3} results in a feasible solution with smallest piece at least 1/3 for all muffin problems. Feasibility is straightforward. That the smallest piece is at least 1/3 is immediate in Case 1. For Case 2, recognize that the final step is a 3M-DAP with $t=u=v=2$, $x_t = 1$, $x_u = x - k/3$, and $x_v = x - (k - 1)/3$. Apply the greedy algorithm of the previous section: $\lambda = -1/3$ and by Lemma \ref{greedy_t=u=v=2}, the smallest piece in the remaining assignment will be of size greater than $(1+\lambda)/2 = 1/3$.
\hspace{\stretch{1}}$\blacksquare$

\begin{lemma} \label{lem:1/3_v1_complexity} Algorithm \ref{alg:1/3} has complexity $\Theta(m)$.
\end{lemma}
\textbf{Proof.} Compute $x = m/s$ and $k = \lfloor 3x \rfloor$. The division of muffins is represented by two matrices, $T_1$ having dimensions $(m - s) \times 3$ and $T_2$ having dimensions $s \times 2$.

Compare $k$ to $3x$: if equal, then the problem solution can be represented by a DAP $(T_1, T_2; U)$, with $U$ having dimensions $s \times (k - 1)$. Assign value $\sfrac{1}{3}$ to all elements of $T_1$ and to all elements of the first $k - 3$ columns of $U$. Assign value $\sfrac{1}{2}$ to all elements of $T_2$ and to all elements of the last two columns of $U$. Initialization and assignment of values to the elements of the three matrices requires $4s(k - 1) = 4(3m - s)$ operations.

If $k$ is not equal to $3x$, then compute $s_u = 3m - ks$ and $s_v = (k + 1)s - 3m$, so that $s_u + s_v = s$ and $(k - 2)s_u + (k - 3)s_v = 3(m - s)$. The problem solution can be represented by a DAP $(T_1, T_2; U, V)$, with $U$ having dimensions $s_u \times k$ and $V$ having dimensions $s_v \times (k - 1)$. Initialization of the four matrices requires $2(3m - s)$ operations. Assign value $\sfrac{1}{3}$ to all elements of $T_1$, to all elements of the first $k - 2$ columns of $U$, and to all elements of the first $k - 3$ columns of $V$. These assignments require $6(m - s)$ operations. Letting $U^{\prime}$ represent the last two columns of $U$ and $V^{\prime}$ represent the last two columns of $V$, Algorithm \ref{alg:greedy} is employed to solve the 3M-DAP $(T_2; U^{\prime}, V^{\prime})$. By Lemma \ref{lem:alg_greedy_complexity}, this last step has complexity $\Theta(s)$.

Recalling that $s < m$, the conclusion follows.
\hspace{\stretch{1}}$\blacksquare$

\section{Additional Lower Bounds}

We can leverage the lower bound on the optimal solution for any 3M-DAP with $t = u = v = 2$ to establish useful lower bounds on more general problems. 
% Note that Algorithm \ref{alg:t_even} and Lemma \ref{lem:LB_on_f;t_even} apply also to problems with $u \le v$.

\fontsize{8}{8}
\begin{algorithm}
\SetKwInOut{Input}{input}
\SetKwInOut{Output}{output}
\Input{$t$ even\\
$T$, $U$, $V$, empty matrices of dimensions $s_t \times t$, $s_u \times u$, and $s_v \times v$, respectively\\
$x_t$, $x_u$, $x_v$, the required row sums for each row of $T$, $U$, and $V$, respectively\\
}
\Output{the filled matrices $T$, $U$, $V$\\
the multiset of entries in $T$ is the union of the multiset of entries in $U$ and $V$\\
}
\BlankLine

divide each row of $T$ into $t/2$ pairs of elements\\
each pair is a row in matrix $T^{\prime}$: dimensions $(ts_t/2) \times 2$, required row sums $x_t^{\prime} = 2x_t/t$\\
divide $[(u - 2)s_u + (v - 2)s_v]/2 = ts_t/2 - s_u - s_v$ rows of $T^{\prime}$ $[x_t/t, x_t/t]$\\
set all elements in the first $u - 2$ columns of $U$ to value $x_t/t$\\
set all elements in the first $v - 2$ columns of $V$ to value $x_t/t$\\
form a $(T^{\prime \prime}; U^{\prime}, V^{\prime})$ 3M-DAP where:\\
\quad $T^{\prime \prime}$ is the matrix consisting of the remaining rows of $T^{\prime}$; $x_t^{\prime \prime} = 2x_t/t$\\
\quad $U^{\prime}$ is the matrix consisting of the remaining two columns of $U$; $x_u^{\prime} = x_u - (u - 2)x_t/t$\\
\quad $V^{\prime}$ is the matrix consisting of the remaining two columns of $V$; $x_v^{\prime} = x_v - (v - 2)x_t/t$\\
use Algorithm \ref{alg:greedy} to solve the $(T^{\prime \prime}; U^{\prime}, V^{\prime})$ 3M-DAP\\
\caption{Algorithm to achieve lower bound for 3M-DAP with $t$ even}
\label{alg:t_even}
\end{algorithm}

\normalsize

\begin{lemma} \label{lem:LB_on_f;t_even} Given a 3M-DAP $P = (T; U, V)$ with $t$ even, $f(P) > \lambda/2 + x_t/t - (u - v)x_t/(2t)$. When $t = 2$, in an optimal solution the largest element has size less than $(x_t - \lambda)/2 + (u - v)x_t/4$.
\end{lemma}
\textbf{Proof.} Apply Algorithm \ref{alg:t_even} to solve $P$. Now 
\[
\begin{split}
  x_u/u &< x_t/t = x_t^{\prime}/2 \Leftrightarrow x_u < ux_t^{\prime}/2 \Leftrightarrow x_u^{\prime} = x_u - (u - 2)x_t^{\prime}/2 < x_t^{\prime}, \quad \text{and}\\
  x_v/v &> x_t/t = x_t^{\prime}/2 \Leftrightarrow x_v > vx_t^{\prime}/2 \Leftrightarrow x_v^{\prime} = x_v - (v - 2)x_t^{\prime}/2 > x_t^{\prime}.
\end{split}
\]
So assigning the remaining rows of $T^{\prime}$ to the last two columns of $U$ and $V$ equates to a 3M-DAP $(T^{\prime \prime}; U^{\prime}, V^{\prime})$ where $t^{\prime \prime} = u^{\prime} = v^{\prime} = 2$ and $x_u^{\prime} = x_u - (u - 2)x_t^{\prime}/2 < x_t^{\prime} = x_t^{\prime \prime} < x_v - (v - 2)x_t^{\prime}/2 = x_v^{\prime}$. Then $\lambda^{\prime} = \lambda - (u - v)x_t^{\prime}/2 = \lambda - (u - v)x_t/t$.

By Lemma \ref{greedy_t=u=v=2}, the smallest piece in an optimal solution to this reduced problem is greater than
\[
  \frac{1}{2}(x_t^{\prime \prime} + \lambda^{\prime}) = \frac{1}{2}\lambda + \frac{x_t}{t} - \frac{1}{2}(u - v)\frac{x_t}{t}.
\] 
\hspace{\stretch{1}}$\blacksquare$

\begin{lemma} \label{lem:alg_t_even_complexity} Algorithm \ref{alg:t_even} has complexity $\Theta(n_t)$.
\end{lemma}
\textbf{Proof.} Algorithm \ref{alg:t_even} begins by constructing an empty matrix $T^{\prime}$ of dimensions $(n_t/2) \times 2$, then assigning the value $x_t/t$ to every element in the first $n_t/2 - s_u - s_v$ rows of $T^{\prime}$, to every element in the first $u - 2$ columns of $U$, and to every element in the first $v - 2$ columns of $V$. This requires $2(n_t - s_u - s_v) + 1$ operations.

Algorithm \ref{alg:t_even} then constructs a $(T^{\prime \prime}; U^{\prime}, V^{\prime})$ 3M-DAP where $n_t^{\prime \prime} = 2(s_u + s_v)$, requiring a total of $4(s_u + s_v) + 7$ operations, then applies Algorithm \ref{alg:greedy} to solve this subproblem, which has complexity $\Theta(s_u + s_v)$ by Lemma \ref{lem:alg_greedy_complexity}. The conclusion follows.
\hspace{\stretch{1}}$\blacksquare$

We can use the lower bound of Lemma \ref{lem:LB_on_f;t_even} to establish a lower bound on 3M-DAPs for which $u = v + 1$. Lemma \ref{lem:lambdaLB_ge} shows, using Algorithm \ref{alg:u=v+1_ge}, that the optimal value is at least $\lambda$ for such problems. Lemma \ref{lem:lambdaLB} strengthens the lower bound to be strict; the algorithm and proof are more complicated and can be found in Appendix \ref{lem:lambdaLB_proof} (this more complicated algorithm does not have linear running time).

\fontsize{8}{8}
\begin{algorithm}
\SetKwInOut{Input}{input}
\SetKwInOut{Output}{output}
\Input{$T$, $U$, $V$, empty matrices of dimensions $s_t \times t$, $s_u \times u$, and $s_v \times v$, respectively\\
$x_t$, $x_u$, $x_v$, the required row sums for each row of $T$, $U$, and $V$, respectively\\
$u = v + 1$\\
}
\Output{the filled matrices $T$, $U$, $V$\\
the multiset of entries in $T$ is the union of the multiset of entries in $U$ and $V$\\
}
\BlankLine

\uIf{$t$ even}{
use Algorithm \ref{alg:t_even} to solve the $(T; U, V)$ problem\\
}
\uElseIf{$s_t \ge (u - 2)s_u + (v - 2)s_v$}{
// $u = 3, v = 2, s_t \ge s_u$\\
\uIf{$s_t = s_u$}{
fill the first column of $U$ with value $\lambda$ and insert these into the first column of $T$\\
fill all remaining elements of $U$ and all elements of $V$ with value $x_v/v$ and insert these into the remaining columns of $T$\\
}
\uElse{
// $s_t > s_u$\\
fill the first column of $U$ with value $\lambda$ and insert these into the first column of $T$\\
form a $(T^{\prime}; U^{\prime}, V^{\prime})$ 3M-DAP where:\\
$\quad T^{\prime}$ is a $(s_u + s_v) \times 2$ matrix where each row is either the unfilled elements of a row of $U$ or a row of $V$; $x_t^{\prime} = x_v$\\ 
$\quad U^{\prime}$ is the $(s_t - s_u) \times t$ unfilled submatrix of $T$; $x_u^{\prime} = x_t$\\
$\quad V^{\prime}$ is the $s_u \times (t - 1)$ unfilled submatrix of $T$; $x_t^{\prime} = x_t - \lambda$\\
$t^{\prime} = 2$ and $u^{\prime} = v^{\prime} + 1$, so use Algorithm \ref{alg:t_even} to solve the $(T^{\prime}; U^{\prime}, V^{\prime})$ 3M-DAP\\
}
}
\uElse{
// $s_t < (u - 2)s_u + (v - 2)s_v$\\
fill the first column of $T$ with value $\lambda$ and insert these values into $U$ and $V$ so that the difference in the numbers of unfilled elements between any two rows is at most 1\\
\uIf{all rows of $U$ and $V$ have the same number of unfilled elements}{
complete the assignment by filling all remaining elements with value $(x_t - \lambda)/(t - 1) > \lambda$\\
}
\uElse{
let $v^{\prime}$ be such that all rows of $U$ and $V$ have either $v^{\prime}$ unfilled elements or $u^{\prime} \equiv v^{\prime} + 1$ unfilled elements\\
form a $(T^{\prime}; U^{\prime}, V^{\prime})$ 3M-DAP where:\\
$\quad T^{\prime}$ is the submatrix of $T$ consisting of all columns except the first; $x_t^{\prime} = x_t - \lambda$\\
$\quad U^{\prime}$ is the matrix consisting of the unfilled elements of those rows of $U$ and $V$ with $u^{\prime}$ unfilled elements; $x_u^{\prime} = x_u - (u - u^{\prime})\lambda = x_v - (v - u^{\prime})\lambda$\\
$\quad V^{\prime}$ is the matrix consisting of the unfilled elements of those rows of $U$ and $V$ with $v^{\prime}$ unfilled elements; $x_v^{\prime} = x_v - (v - v^{\prime})\lambda = x_u - (u - v^{\prime})\lambda = x_u - \lambda - (u - u^{\prime})\lambda = x_u^{\prime} - \lambda$\\
$t^{\prime}$ is even and $u^{\prime} = v^{\prime} + 1$, so use Algorithm \ref{alg:t_even} to solve the $(T^{\prime}; U^{\prime}, V^{\prime})$ 3M-DAP\\
}
}
\caption{Algorithm to achieve lower bound for 3M-DAP with $u=v+1$}
\label{alg:u=v+1_ge}
\end{algorithm}

\normalsize

\begin{lemma} \label{lem:lambdaLB_ge} Given a 3M-DAP $P = (T; U, V)$ with $u = v + 1$, $f(P) \ge \lambda = x_u - x_v$.
\end{lemma}
\textbf{Proof.} Apply Algorithm \ref{alg:u=v+1_ge} to solve $P$. When $t$ is even, apply the previous lemma to conclude that
\[
  f(P) > \frac{1}{2}\lambda + \frac{1}{2}\frac{x_t}{t} > \lambda,
\]
because $x_t/t > x_u/u > (x_u - x_v)/(u - v) = \lambda$.

Now we investigate when $t$ is odd. First, suppose that $s_t \ge (u - 2)s_u + (v - 2)s_v$. If $v \ge 3$, then $u \ge 4$, so $(v - 2)s_v \ge n_v/3$ and $(u - 2)s_u \ge n_u/2$, so $(u - 2)s_u + (v - 2)s_v > (n_u + n_v)/3 = n_t/3$. But $t \ge 3$, so $s_t \le n_t/3$. But this contradicts our supposition so we must have $v = 2$ and $u = 3$. Then the supposition is that $s_t \ge s_u$.

If $s_t = s_u$, Algorithm \ref{alg:u=v+1_ge} produces a solution with only two distinct values: $\lambda$ and $x_v/v > \lambda$.

If $s_t > s_u$, we obtain $f(P) \ge \min\{\lambda, f(P^{\prime})\}$, where the subproblem $P^{\prime} = (T^{\prime}; U^{\prime}, V^{\prime})$ has $t^{\prime} = 2$, $u^{\prime} = v^{\prime} + 1$, and $\lambda^{\prime} = \lambda$. By the first part of the proof, $f(P^{\prime}) > \lambda$, so the conclusion follows.

Second, suppose that $s_t < (u - 2)s_u + (v - 2)s_v$. After filling the first column of $T$ with value $\lambda$ and inserting these values into $U$ and $V$, we have two cases. When all rows of $U$ and $V$ have the same number of unfilled elements, Algorithm \ref{alg:u=v+1_ge} produces a solution with only two distinct values: $\lambda$ and $(x_t - \lambda)/(t - 1) > \lambda$.

Otherwise, we obtain $f(P) \ge \min\{\lambda, f(P^{\prime})\}$, where the subproblem $P^{\prime} = (T^{\prime}; U^{\prime}, V^{\prime})$ has $t^{\prime}$ even, $u^{\prime} = v^{\prime} + 1$, and $\lambda^{\prime} = \lambda$. By the first part of the proof, $f(P^{\prime}) > \lambda$, so the conclusion follows.\hspace{\stretch{1}}$\blacksquare$

\begin{lemma} \label{lem:alg:u=v+1_ge_complexity} Algorithm \ref{alg:u=v+1_ge} has complexity $\Theta(n_t)$.
\end{lemma}
\textbf{Proof.} Algorithm \ref{alg:u=v+1_ge} operates in one of two ways, depending on the problem characteristics. In the first way, when $s_t = k(s_u + s_v) + s_u$ for some nonnegative integer $k$, the problem is solved directly by assigning the value $\lambda$ to all elements of the first column of $T$, of the first $k + 1$ columns of $U$, and of the first $k$ columns of $V$. All the remaining elements of $T$, $U$, and $V$ are assigned the value $(x_t - \lambda)/(t - 1)$ (which equals $x_v/v$ when $s_t = s_u$. Clearly, this direct solution has complexity $\Theta(n_t)$.

The second method involves assigning the value $\lambda$ to all elements of the first column of either $U$ or $T$, and similarly assigning this value to elements in either $T$ or in $U$ and $V$, so that the remaining elements form a reduced 3M-DAP that can be solved in linear time using Algorithm \ref{alg:t_even}.
\hspace{\stretch{1}}$\blacksquare$

\begin{lemma} \label{lem:lambdaLB} Given a 3M-DAP $P = (T; U, V)$ with $u = v + 1$, $f(P) > \lambda = x_u - x_v$.
\end{lemma}

\section{0-Problems of Type 1}

In this section we investigate 0-problems of type 1, namely any 3M-DAP with $s_t \le (v - 1)s_v$ (note that this requires that $s_v > 0$ and $v > 1$). Observe that $n_u < (t - 1)s_t$ and $n_u \le (t - 2)s_t + (v - 2)s_v$. Insert the elements of $U$ into $T$, as evenly as possible across the rows of $T$ so that the difference in the number of elements inserted from $U$ into any two rows of $T$ is at most 1. Write
\[
\begin{split}
  q_t &= \left \lfloor \frac{n_u}{s_t} \right \rfloor\\
  r_t &= n_u - q_t s_t.
\end{split}
\]
Here $q_t$ represents the number of columns of $T$ that are completely filled, so that each row receives either $q_t$ or $q_t + 1$ elements, and $r_t$ represents the number of rows that receive $q_t + 1$ elements.

\fontsize{8}{8}
\begin{algorithm}
\SetKwInOut{Input}{input}
\SetKwInOut{Output}{output}
\Input{$T$, $U$, $V$, empty matrices of dimensions $s_t \times t$, $s_u \times u$, and $s_v \times v$, respectively\\
$x_t$, $x_u$, $x_v$, the required row sums for each row of $T$, $U$, and $V$, respectively\\
$s_t \le (v - 1)s_v$\\
}
\Output{the filled matrices $T$, $U$, $V$\\
the multiset of entries in $T$ is the union of the multiset of entries in $U$ and $V$\\
}
\BlankLine

set all elements of $U$ to value $x_u/u$\\
$q_t \leftarrow \left \lfloor n_u/s_t \right \rfloor$\\
$r_t \leftarrow n_u - q_t s_t$\\

\uIf{$r_t = 0$}{
	set all elements of $V$ to value $x_v/v$\\
	fill each row of $T$ with $q_t$ elements of value $x_u/u$ and $t - q_t$ elements of value $x_v/v$\\
}
\uElseIf{$q_t < t - 2$}{
	for $r_t$ rows of $T$, set $q_t + 1$ elements to value $x_u/u$; call the unfilled submatrix $V^{\prime}$\\
	set $x_v^{\prime} = x_t - (q_t + 1)x_u/u$\\
	for the remaining $s_t - r_t$ rows of $T$, set $q_t$ elements to value $x_u/u$; call the unfilled submatrix $U^{\prime}$\\
	set $x_u^{\prime} = x_t - q_tx_u/u$\\
	set $T^{\prime} = V$\\
	use Algorithm \ref{alg:u=v+1} to solve the $(T^{\prime}; U^{\prime}, V^{\prime})$ 3M-DAP\\
}
\uElse{
%	\tcc{$q_t = t - 2$}
	\tcp*[h]{$q_t = t - 2$}\\
	for $r_t$ rows of $T$, set $t - 1$ elements to value $x_u/u$; set the last element to $\rho = x_t - (t - 1)x_u/u > x_u/u$\\
	for the remaining $s_t - r_t$ rows of $T$, set the first $t - 2$ elements to value $x_u/u$; call the unfilled submatrix $T^{\prime}$: then $t^{\prime} = 2$\\
	$x_t^{\prime} \leftarrow x_t - (t - 2)x_u/u = \rho + x_u/u$\\

	$q_v \leftarrow \left \lfloor r_t/s_v \right \rfloor$\\
	$r_v \leftarrow r_t - q_v s_v$\\

	\uIf{$r_v = 0$}{
		set all elements of $T^{\prime}$ to value $\sigma = [x_t - (t - 2)x_u/u]/2 > x_u/u$\\
		fill each row of $V$ with $q_v$ elements of value $\rho$ and $v - q_v$ elements of value $\sigma$\\
	}
	\uElse{
		\tcp*[h]{$q_v < v - 2$}\\
		for $r_v$ rows of $V$, set $q_v + 1$ elements to value $\rho$; call the unfilled submatrix $U^{\prime}$\\
		set $x_u^{\prime} = x_v - (q_v + 1)\rho$\\
		for the remaining $s_v - r_v$ rows of $V$, set $q_v$ elements to value $\rho$; call the unfilled submatrix $V^{\prime}$\\
		set $x_v^{\prime} = x_v - q_v\rho$\\
		use Algorithm \ref{alg:t_even} to solve the $(T^{\prime}; U^{\prime}, V^{\prime})$ 3M-DAP\\
	}
}

\caption{Algorithm for 3M-DAP that are 0-problems of type 1: achieves $x_u/u$ upper bound}
\label{alg:0prob_type1}
\end{algorithm}

\normalsize

\begin{lemma} \label{qt<t-2} Given a 3M-DAP $P = (T; U, V)$ with $s_t \le (v - 1)s_v$ and either $r_t = 0$, or $q_t < t - 2$ and $r_t > 0$, $f(P) = x_u/u$. Algorithm \ref{alg:0prob_type1} produces an optimal solution in which all elements of $V$ have size strictly greater than $x_u/u$.
\end{lemma}
\textbf{Proof.} Apply Algorithm \ref{alg:0prob_type1} to solve $P$. When $r_t = 0$, it suffices to confirm that the rowsums of $T$ are correct:
\[
  q_t\frac{x_u}{u} + (t - q_t)\frac{x_v}{v} = \frac{1}{s_t}(s_ux_u + s_vx_v) = x_t.
\]
Then the $x_u/u$ bound is achieved and the solution is optimal.

When $q_t < t - 2$ and $r_t > 0$, the algorithm constructs a subproblem $P^{\prime} = (T^{\prime}; U^{\prime}, V^{\prime})$ that is a 3M-DAP. Check:
\[
\begin{split}
  q_ts_t &< n_u < (q_t + 1)s_t\\
  \Leftrightarrow \frac{q_t}{t} &< \frac{n_u}{n_t} < \frac{q_t + 1}{t}\\
  \Leftrightarrow \frac{q_t}{t}\frac{x_u}{u} + \left(1 - \frac{q_t}{t}\right)\frac{x_v}{v} &> \frac{x_t}{t} > \frac{q_t + 1}{t}\frac{x_u}{u} + \left(1 - \frac{q_t + 1}{t}\right)\frac{x_v}{v}\\
  \Leftrightarrow \frac{x_u^{\prime}}{u^{\prime}} = \frac{x_t - q_tx_u/u}{t - q_t} &< \frac{x_v}{v} = \frac{x_t^{\prime}}{t^{\prime}} < \frac{x_t - (q_t + 1)x_u/u}{t - q_t - 1} = \frac{x_v^{\prime}}{v^{\prime}}.
\end{split}
\]
Because $u^{\prime} = v^{\prime} + 1$, we can apply Lemma \ref{lem:lambdaLB} to conclude that $f(P^{\prime}) > x_u^{\prime} - x_v^{\prime} = x_u/u$. Then the $x_u/u$ bound is achieved and the solution is optimal.
\hspace{\stretch{1}}$\blacksquare$

Otherwise $q_t = t - 2$ and $r_t > 0$. Then $r_t$ rows of $T$ have $t - 1$ elements filled; the last element must take value $\rho = x_t - (t - 1)x_u/u > x_u/u$. Insert these elements into $V$, as evenly as possible across the rows of $V$ so that the difference in the number of elements of value $\rho$ inserted into any two rows of $V$ is at most 1. Write
\[
\begin{split}
  q_v &= \left \lfloor \frac{r_t}{s_v} \right \rfloor\\
  r_v &= r_t - q_v s_v.
\end{split}
\]
Here $q_v$ represents the number of columns of $V$ that are completely filled, so that each row receives either $q_v$ or $q_v + 1$ elements of value $\rho$, and $r_v$ represents the number of rows that receive $q_v + 1$ such elements.

\begin{lemma} \label{qt=t-2} Given a 3M-DAP $P = (T; U, V)$ with $v > 1$ and $n_u = (t - 2)s_t + r_t$, where $0 < r_t \le (v - 2)s_v$ or $r_t = (v - 1)s_v$, $f(P) = x_u/u$. Algorithm \ref{alg:0prob_type1} produces an optimal solution in which all elements of $V$ have size strictly greater than $x_u/u$.
\end{lemma}
\textbf{Proof.} Apply Algorithm \ref{alg:0prob_type1} to solve $P$. When $r_v = 0$, it suffices to check that the rowsums of $V$ are correct. First,
\[
  n_u = (t - 2)s_t + r_t \Leftrightarrow n_v = n_t - n_u = 2s_t - r_t = 2s_t - q_vs_v \Leftrightarrow (v + q_v)s_v = 2s_t.
\]
Then
\[
\begin{split}
  q_v\rho + (v - q_v)\sigma &= q_v\left[x_t - (t - 1)\frac{x_u}{u}\right] + (v - q_v)\frac{1}{2}\left[x_t - (t - 2)\frac{x_u}{u}\right]\\
  &= \frac{1}{2}(v + q_v)\left(x_t - t\frac{x_u}{u}\right) + v\frac{x_u}{u}\\
  &= \frac{s_t}{s_v}\left(x_t - t\frac{x_u}{u}\right) + v\frac{x_u}{u}\\
  &= \frac{1}{s_v}\left[s_ux_u + s_vx_v - (us_u + vs_v)\frac{x_u}{u}\right] + v\frac{x_u}{u}\\
  &= x_v - v\frac{x_u}{u} + v\frac{x_u}{u} = x_v.
\end{split}
\]
Then the $x_u/u$ bound is achieved and the solution is optimal.

When $r_v > 0$, the algorithm constructs a subproblem $P^{\prime} = (T^{\prime}; U^{\prime}, V^{\prime})$ that is a 3M-DAP with $t^{\prime} = 2$ and then employs Algorithm \ref{alg:t_even} to solve this subproblem.

The elements of $V$ will be $r_t$ having value $\rho$ and $2(s_t - r_t)$ with average value $\sigma = (\rho + x_u/u)/2 < \rho$, so the average value of an element in $V$ is less than $\rho$, i.e., $x_v/v < \rho$.

The unfilled elements of $V$ are divided into two matrices $U^{\prime}$ and $V^{\prime}$, where $U^{\prime}$ has $r_v$ rows, $u^{\prime} = v - q_v - 1$ columns, and required row sum $x_u^{\prime} = x_v - (q_v + 1)\rho$, while $V^{\prime}$ has $s_v - r_v$ rows, $v^{\prime} = v - q_v$ columns, and required row sum $x_v^{\prime} = x_v - q_v\rho$. The elements of $U^{\prime}$ and $V^{\prime}$ are inserted into the matrix $T^{\prime}$, which represents the last two columns of the remaining rows of $T$, i.e., has $s_t - r_t$ rows, $t^{\prime} = 2$ columns, and required row sum $x_t^{\prime} = x_t - (t - 2)x_u/u = \rho + x_u/u$. Because $x_v/v < \rho$, we have
\[
  \frac{x_v - (q_v + 1)\rho}{v - q_v - 1} < \frac{x_v - q_v\rho}{v - q_v} < \frac{x_v}{v}.
\]
Then it must be the case that
\[
  \frac{x_u^{\prime}}{u^{\prime}} = \frac{x_v - (q_v + 1)\rho}{v - q_v - 1} < \frac{x_t - (t - 2)x_u/u}{2} = \frac{x_t^{\prime}}{t^{\prime}} < \frac{x_v - q_v\rho}{v - q_v} = \frac{x_v^{\prime}}{v^{\prime}}.
\]
Then we can apply Lemma \ref{lem:LB_on_f;t_even} to conclude that
\[
\begin{split}
  f(P^{\prime}) &> \frac{1}{2}(x_t^{\prime} + \lambda^{\prime}) - \frac{1}{4}(u^{\prime} - v^{\prime})x_t^{\prime}\\
  &= \frac{3}{4}x_t^{\prime} + \frac{1}{2}\lambda^{\prime}\\
  &= \frac{3}{4}\left(\rho + \frac{x_u}{u}\right) - \frac{1}{2}\rho\\
  &= \frac{1}{4}\rho + \frac{3}{4}\frac{x_u}{u} > \frac{x_u}{u}.
\end{split}
\]
Then the $x_u/u$ bound is achieved and this solution is optimal.
\hspace{\stretch{1}}$\blacksquare$

\begin{theorem} \label{thm:0type1} Any 3M-DAP $P = (T; U, V)$ with $s_t \le (v - 1)s_v$ is a 0-problem: $f(P) = x_u/u$. Algorithm \ref{alg:0prob_type1} produces an optimal solution in which all elements of $V$ have size strictly greater than $x_u/u$.
\end{theorem}

\begin{lemma} \label{lem:alg_0prob_type1_complexity} Algorithm \ref{alg:0prob_type1} has complexity $\Theta(n_t)$.
\end{lemma}
\textbf{Proof.} Algorithm \ref{alg:0prob_type1} operates in one of two ways, depending on the problem characteristics. In the first way, when either $r_t = 0$ or $q_t = t - 2$ and $r_v = 0$, the problem is solved directly. The value $x_u/u$ is assigned to all elements of $U$ and to the same number of elements of $T$. The remaining elements of $T$ and all elements of $V$ are assigned one of the values $x_v/v$, $\rho$, and $\sigma$. Clearly, this direct solution has complexity $\Theta(n_t)$.

The second method involves assigning the value $x_u/u$ to all elements of $U$ and to the same number of elements of $T$, and further, if $q_t = t - 2$, assigning the value $\rho$ to some elements of $T$ and to the same number of elements of $V$. The remaining elements form a reduced 3M-DAP that can be solved in linear time using either Algorithm \ref{alg:u=v+1} (by Lemma \ref{lem:alg_u=v+1_complexity}) or Algorithm \ref{alg:t_even} (by Lemma \ref{lem:alg_t_even_complexity}).
\hspace{\stretch{1}}$\blacksquare$

\subsection{0-Problems of Type 1 in the Context of Extended Families}

Here we put 0-problems of type 1 into the context of extended families by investigating for what values of $x$ a 3M-DAP $P = (T; U, V)$ is a 0-problem of type 1. Define
\[
  x_{\infty} \equiv \frac{\lambda + \gamma (v - 1)t}{(v - 1)t + u - v}u = \frac{x_u + (v - 1)x_t - x_v}{u + (v - 1)t - v}u.
\]

\begin{theorem} \label{thm:0type1_x<=x_infty} A 3M-DAP $P = (T; U, V)$ is a 0-problem of type 1 if and only if $x \le x_{\infty}$.
\end{theorem}
\textbf{Proof.} If $v = 1$, there are no 0-problems. If $t = v = 2$, then $s_t \le (v - 1)s_v \Leftrightarrow n_t = n_v$, so $n_u = 0$ and $x_t = x_v$, so the problem is trivial and the claim of the theorem is a tautology. Otherwise, either $t > 2$ or $v > 2$, so $(v - 1)t - v > 0$. Then $s_t \le (v - 1)s_v$ implies that
\[
  \frac{x}{u} = \frac{s_tx_t - s_vx_v}{s_tt - s_vv} \le \frac{(v - 1)x_t - x_v}{(v - 1)t - v},
\]
from which it follows that
\[
  \frac{x}{u} \le \frac{x_u + (v - 1)x_t - x_v}{u + (v - 1)t - v} = x_{\infty}.
\]
\hspace{\stretch{1}}$\blacksquare$

\section{Sketch of Solution Strategy for Remaining Problems}

The previous section showed how to find an optimal solution when $s_t \le (v - 1)s_v$. The following sections address the more complicated case when $s_t > (v - 1)s_v$. Recall that given any solution, the matrices $T$ and $V$ can be divided into maximal, inseparable pairs, and given the value of the problem exceeds some lower bound, that those pairs must be $b$-pairs for some values of $b$.

The next section describes an algorithm (Algorithm \ref{alg:(T,V)-pair}) for completing a $b$-pair: given the values of the $U$ elements in the pair, the algorithm identifies the values of the $V$-elements in the pair and inserts them into the rows of $T$ and $V$ so that the required rowsums are achieved. Provided the smallest value of any of the $U$-elements exceeds some lower bound, the values of all the $V$-elements in pair will exceed $x_u/u$.

When $b^*$ is integer, we show in the Section \ref{0-probs_type2} that the 3M-DAP is a 0-problem: the $x_u/u$ bound is achieved by setting all values in $U$ to $x_u/u$, then using Algorithm \ref{alg:(T,V)-pair} to complete a $b^*$-pair, and inserting multiple copies of this pair into $T$ and $V$.

When $b^*$ is not integer, we seek a solution in which all maximal inseparable pairs are either $\hat{b}$- or $(\hat{b} - 1)$-pairs. The number of the $U$-elements in a $b$-pair must be $[(v - 1)b + 1]t - bv$ and these must sum to $[(v - 1)b + 1]x_t - bx_v$. This number and sum are the same for all $\hat{b}$-pairs in the solution and likewise are the same for all $(\hat{b} - 1)$-pairs in the solution.

To determine the values in the matrix $U$, we can create a reduced problem that is a new 3M-DAP $P^{\prime} = (T^{\prime}; U^{\prime}, V^{\prime})$, where $T^{\prime} = U$, there is one row in $U^{\prime}$ for each $\hat{b}$-pair, and one row in $V^{\prime}$ for each $(\hat{b} - 1)$-pair. Given a solution to this reduced problem, we apply Algorithm \ref{alg:(T,V)-pair} to complete each $\hat{b}$-pair and each $(\hat{b} - 1)$-pair. If the smallest value in the reduced problem, i.e., in $U$, exceeds a lower bound specified in the statement of Proposition \ref{prop:Vsub>x_u/u}, then all elements in $V$ will be greater than $x_u/u$, and the smallest value of the reduced problem will also be the smallest value of the original problem.

The reduced problem can be solved in a similar manner. Either the reduced problem is a 0-problem, in which case we solve it directly, or we can reduce it further. The reductions continue until we reach a 0-problem that is solved directly. Solutions to the larger problems are constructed from those to the smaller problems using Algorithm \ref{alg:(T,V)-pair} to complete each $b$-pair.

\section{Completing a $b$-pair}

Let $(A, B)$ be a $b$-pair from $(T, V)$. Suppose we are given the values of the $[(v - 1)b + 1]t - bv$ elements of $U$ that are contained in $A$. Here we describe an algorithm for inserting these elements into $A$, constructing the matrix $B$, and inserting its elements into $A$ so that the required rowsums for all rows in the two submatrices are obtained.

Arrange the given elements of $U$ in any order into a vector $\boldsymbol{u}$ of length $[(v - 1)b + 1]t - bv$ (for our analysis only the minimum value of these elements is important, not the order). 

When $b = 0$, the matrix $B$ has no rows and the matrix $A$ has one row and contains only the elements of $U$ contained in $\boldsymbol{u}$: we can simply set the one row of $A$ to be $\boldsymbol{u}$.

Otherwise, $b > 0$. When $v = 1$, all elements of $V$ take value $x_v$, so
\begin{equation} \label{V_sub_v=1}
B = \begin{bmatrix} 
    x_v\\
    \vdots\\
    x_v
    \end{bmatrix}.
\end{equation}
The matrix $A$ has only one row: the given elements of $U$ and the elements of $B$ are inserted into $A$ as follows:
\begin{equation} \label{T_sub_v=1}
A = \begin{bmatrix} 
    u_1 & \cdots & u_{t-b} & x_v & \cdots & x_v
    \end{bmatrix}.
\end{equation}

When $v > 1$ and $b = 1$ the matrix $B$ has only one row:
\begin{equation} \label{V_sub_b=1}
B = \begin{bmatrix} 
    w_1 & w_2 & \cdots  & w_{v - 1} & w_v
    \end{bmatrix}.
\end{equation}
The given elements of $U$ and the elements of $B$ are inserted into $A$ as follows:
\begin{equation} \label{T_sub_b=1}
A = \begin{bmatrix} 
    u_{11} & \cdots  & u_{1, t - 2} & u_{1, t - 1} & w_1 \\
    \vdots & \ddots & \vdots & \vdots & \vdots \\
    u_{v1} & \cdots  & u_{v, t - 2} & u_{v, t - 1} & w_v
    \end{bmatrix},
\end{equation}
where for ease of notation we relabel the elements of $\boldsymbol{u}$ by their position in the $A$ matrix.

Otherwise, $v > 1$ and $b > 1$. In this case, divide the elements of $B$ into three vectors: $\boldsymbol{w}$ of length $b(v - 2) + 2$, and $\boldsymbol{y}$ and $\boldsymbol{z}$, each of length $b - 1$. Then $B$ contains these elements as follows:
\begin{equation} \label{V_sub}
B = \begin{bmatrix} 
    w_1 & w_2 & \cdots  & w_{v - 1} & y_1 \\
    z_1 & w_v & \cdots  & w_{2v - 3} & y_2 \\
    \vdots & \vdots & \ddots & \vdots & \vdots \\
    z_{b - 2} & w_{(b - 2)(v - 2) + 2} & \cdots  & w_{(b - 1)(v - 2) + 1} & y_{b - 1} \\
    z_{b - 1} & w_{(b - 1)(v - 2) + 2} & \cdots  & w_{b(v - 2) + 1} & w_{b(v - 2) + 2}
    \end{bmatrix}.
\end{equation}
The given elements of $U$ and the elements of $B$ are inserted in to $A$ as follows:
\begin{equation} \label{T_sub}
A = \begin{bmatrix} 
    u_{11} & \cdots  & u_{1, t - 2} & u_{1, t - 1} & w_1 \\
    \vdots & \ddots & \vdots & \vdots & \vdots \\
    u_{b(v - 2) + 2, 1} & \cdots  & u_{b(v - 2) + 2, t - 2} & u_{b(v - 2) + 2, t - 1} & w_{b(v - 2) + 2} \\
    u_{b(v - 2) + 3, 1} & \cdots  & u_{b(v - 2) + 3, t - 2} & y_1 & z_1 \\
    \vdots & \ddots & \vdots & \vdots & \vdots \\
    u_{b(v - 1) + 1, 1} & \cdots  & u_{b(v - 1) + 1, t - 2} & y_{b - 1} & z_{b - 1}
    \end{bmatrix},
\end{equation}
where again we relabel the elements of $\boldsymbol{u}$ by their position in the $A$ matrix.

\fontsize{8}{8}
\begin{algorithm}
\SetKwInOut{Input}{input}
\SetKwInOut{Output}{output}
\Input{$A$, an empty matrix of dimensions $[(v - 1)b + 1] \times t$\\
$x_t$, the required row sum for each row of $A$\\
$B$, an empty matrix of dimensions $b \times v$\\
$x_v$, the required row sum for each row of $B$\\
$\boldsymbol{u}$, a vector of length $[(v - 1)b + 1]t - bv$, containing the given elements\\
\quad \quad of $U$ in any order\\
}
\Output{the filled matrices $A$ and $B$\\
insertion of elements into $(A, B)$ is as in (\eqref{T_sub_v=1}, \eqref{V_sub_v=1}), (\eqref{T_sub_b=1}, \eqref{V_sub_b=1}), or (\eqref{T_sub}, \eqref{V_sub})\\
the multiset of entries in $A$ is the union of the multiset of entries in $\boldsymbol{u}$ and $B$\\
}
\BlankLine
\uIf{$b = 0$}{
	insert $\boldsymbol{u}$ into the one row of $A$\\
}
\uElse{
	insert the elements of $\boldsymbol{u}$ into $A$ as in \eqref{T_sub_v=1}, \eqref{T_sub_b=1} or \eqref{T_sub}, relabeling as necessary\\
	\uIf{$v = 1$}{
	  set remaining elements of $A$ and all elements of $B$ to $x_v$\\
	}
	\uElse{
		\For{$1 \le i \le b(v - 2) + 2$}{compute $w_i = x_t - \sum_{j = 1}^{t - 1} u_{ij}$ and insert into $A$ and $B$\\
		}
		\uIf{$b > 1$}{
	  		compute $y_1 = x_v - \sum_{j = 1}^{v - 1} w_j$ and insert into $A$ and $B$\\
	  		\For{$1 \le i \le b - 2$}{
				compute $z_i = x_t - y_i - \sum_{j = 1}^{t - 2} u_{b(v - 2) + 2 + i, j}$ and insert into $A$ and $B$\\
	  			compute $y_{i + 1} = x_v - z_i - \sum_{j = 2}^{v - 1} w_{i(v - 2) + j}$ and insert into $A$ and $B$\\
	  		}
			compute $z_{b - 1} = x_t - y_{b - 1} - \sum_{j = 1}^{t - 2} u_{b(v - 1) + i, j}$ and insert into $A$ and $B$\\
		}
	}
}
\caption{Algorithm for completing a maximal, inseparable $(T, V)$-pair}
\label{alg:(T,V)-pair}
\end{algorithm}

\normalsize

\begin{proposition} \label{prop:Vsub>x_u/u} Consider $(A, B)$, a $b$-pair for a 3M-DAP $P = (T; U, V)$ problem with $x > x_{\infty}$. Suppose that the elements of $U$ in the $(T, V)$-pair are given and Algorithm \ref{alg:(T,V)-pair} is used to fill the elements of $A$ and $B$. Let $Y$ be the smallest element among the given elements of $U$ and, if $v > 1$, suppose that
\[
  Y > \frac{x_u/u + (v - 1)x_t - x_v}{(v - 1)(t - 1)}.
\]
Then all elements of $B$ have value greater than $x_u/u$.
\end{proposition}
\textbf{Proof.} If $v = 1$, all elements of $B$ have value $x_v > x_u/u$. Otherwise, we first show that $x_t - (t - 1)Y$ is an upper bound on all elements of $B$ by using the fact that all the given elements of $U$ have value at least $Y$. First, when $b = 1$, this follows directly by examining the rows of $A$:
\[
  w_i \le x_t - (t - 1)Y \quad \quad \text{for } 1 \le i \le v.
\]
Otherwise, $b > 1$. Then from the rows of $A$ we have
\begin{alignat*}{2}
  w_i &\le x_t - (t - 1)Y &\quad \quad& \text{for } 1 \le i \le b(v - 2) + 2\\
  z_j &\le x_t - (t - 2)Y - y_j && \text{for } 1 \le j \le b - 1.
\end{alignat*}
Using these, from the rows of $B$ we obtain
\[
\begin{split}
  y_j &\ge x_v - (v - 2)[x_t - (t - 1)Y] - z_{j - 1}\\
  &\ge x_v - (v - 2)[x_t - (t - 1)Y] - x_t + (t - 2)Y + y_{j - 1}\\
  &= x_v - (v - 1)x_t + [(v - 1)t - v]Y + y_{j - 1},\\
  z_{j - 1} &\ge x_v - (v - 2)[x_t - (t - 1)Y] - y_j\\
  &\ge x_v - (v - 2)[x_t - (t - 1)Y] - x_t + (t - 2)Y + z_j\\
  &= x_v - (v - 1)x_t + [(v - 1)t - v]Y + z_j.
\end{split}
\]
Now $x > x_{\infty}$, so $x_u/u > [(v - 1)x_t - x_v]/[(v - 1)t - v]$, so given the assumed lower bound on $Y$, it follows that
\[
  Y > \frac{(v - 1)x_t - x_v}{(v - 1)t - v} \Leftrightarrow x_v - (v - 1)x_t + [(v - 1)t - v]Y > 0.
\]
Then $y_j > y_{j - 1}$ and $z_j < z_{j - 1}$. The smallest element of $\boldsymbol{y}$ is $y_1$ and
\[
  y_1 = x_v - \sum_{i=1}^{v-1} w_i \ge x_v - (v - 1)x_t + (v - 1)(t - 1)Y > \frac{x_u}{u} > Y,
\]
where the second-to-last inequality follows from the assumed lower bound on $Y$. Similarly, the smallest element of $\boldsymbol{z}$ is $z_{b - 1}$ and $z_{b - 1} > x_u/u > Y$.

The largest element of $\boldsymbol{y}$ is $y_{b - 1}$, and
\[
  y_{b - 1} \le x_t - (t - 2)Y - z_{b - 1} < x_t - (t - 1)Y.
\]
Similarly, the largest element of $\boldsymbol{z}$ is $z_1$, and $z_1 < x_t - (t - 1)Y$. Then $x_t - (t - 1)Y$ is an upper bound on all elements of $B$.

Then a lower bound on any element of $B$ is given by supposing that all other elements in the same row of $B$ have size $x_t - (t - 1)Y$; this lower bound is
\[
  x_v - (v - 1)[x_t - (t - 1)Y] = x_v - (v - 1)x_t + (v - 1)(t - 1)Y > \frac{x_u}{u},
\]
where the inequality follows from the assumed lower bound on $Y$.
\hspace{\stretch{1}}$\blacksquare$

\begin{lemma} \label{lem:alg_(T,V)-pair_complexity} Suppose that whenever Algorithm \ref{alg:(T,V)-pair} is applied to complete a $(T, V)$-pair with $v = 1$, then $b \le t - 2$. Then Algorithm \ref{alg:(T,V)-pair} has complexity $\Theta(n_a) = \Theta([(v - 1)b + 1]t)$. More specifically, the number of operations taken by the algorithm is bounded below by $n_a$ and above by $3n_a$.
\end{lemma}
\textbf{Proof.} The vector $\boldsymbol{u}$ is inserted into $A$, requiring one operation per element for a total of $n_a - bv$ operations. When $b = 0$, the algorithm is complete and requires $n_a$ operations.

When $b > 0$ and $v = 1$, each element of $B$ must be set to $x_v$ and inserted into $A$, requiring two operations per element for a total of $2b$ operations. Therefore, the total number of operations required is $n_a + b$, which lies between $n_a$ and $2n_a$ (by the supposition in the statement of the lemma and on observing that $n_a = t$ in this instance).

\begin{table}
\centering
\begin{tabular}{ccc}
  \hline
  Element type& No. elements& No. operations/element\\
  \hline
  $w_i$& $(v - 2)b + 2$& $t - 1$\\
  $y_i$& $b - 1$& $v - 1$\\
  $z_i$& $b - 1$& $t - 1$\\
  \hline
\end{tabular}
\caption{number of operations to compute elements of $B$ using Algorithm \ref{alg:(T,V)-pair}}
\label{tab:B_elts_ops}
\end{table}

When $v > 1$, the elements of $B$ must be computed and inserted into $A$ and $B$. The insertions require two operations per element for a total of $2bv$ insertions. 

When $b = 1$, each element of $B$ is in a separate row of $A$ (is a $w_i$-element) so is computed using $t - 1$ operations for a total of $(t - 1)v$ operations. Then the total number of operations is $n_a + tv = 2n_a$.

When $b > 1$, there are three types of elements computed by Algorithm \ref{alg:(T,V)-pair}, as shown in Table \ref{tab:B_elts_ops}. Then the total number of operations to compute the elements of $B$ is
\[
  [(v - 1)b + 1](t - 1) + (b - 1)(v - 1) = [(v - 1)b + 1]t - v = n_a - v.
\]
Then the total number of operations is $2n_a + (b - 1)v < 3n_a$, so this lies between $2n_a$ and $3n_a$. ($(b - 1)v < n_a$ because $t \ge 2$ and $v \ge 2$, so $n_a - (b - 1)v = [bv - (b - 1)]t - (b - 1)v \ge 4$).

Then under any circumstances, the number of operations taken by the algorithm is bounded below by $n_a$ and above by $3n_a$, so has complexity $\Theta(n_a)$.
\hspace{\stretch{1}}$\blacksquare$

\section{0-Problems of Type 2} \label{0-probs_type2}

In this section we discuss 0-problems of type 2, namely any 3M-DAP with $r_v = n_u - (t - 2)s_t - (v - 2)s_v > 0$ and $b^* = 2s_v/r_v$ integer. We prove that any such problem is indeed a 0-problem, thereby establishing Conjecture \ref{conj:b*}.

\fontsize{8}{8}
\begin{algorithm}
\SetKwInOut{Input}{input}
\SetKwInOut{Output}{output}
\Input{$T$, $U$, $V$, empty matrices of dimensions $s_t \times t$, $s_u \times u$, and $s_v \times v$, respectively\\
$x_t$, $x_u$, $x_v$, the required row sums for each row of $T$, $U$, and $V$, respectively\\
$s_uu + s_vv = s_tt$\\
$s_ux_u + s_vx_v = s_tx_t$\\
$r_v = n_u - (t - 2)s_t + (v - 2)s_v > 0$\\
$b^* = 2s_v/r_v$ integer\\
}
\Output{the filled matrices $T$, $U$, $V$\\
the multiset of entries in $T$ is the union of the multiset of entries in $U$ and $V$\\
}
\BlankLine
set all elements of $U$ to $x_u/u$\\
initialize a $b^*$-pair $(A, B)$\\
set $\boldsymbol{u}$ to be a vector of length $[(v - 1)b^* + 1]t - b^*v$ with all elements equal to $x_u/u$\\
apply Algorithm \ref{alg:(T,V)-pair} with inputs $A$, $x_t$, $B$, $x_v$, and $\boldsymbol{u}$ to complete the pair $(A, B)$\\
insert $A$ into $T$ $r_v/2$ times and $B$ into $V$ $r_v/2$ times\\

\caption{Algorithm for 3M-DAP that are 0-problems of type 2: achieves $x_u/u$ upper bound}
\label{alg:0prob_type2}
\end{algorithm}

\normalsize

\begin{theorem} \label{thm:0type2} Any 3M-DAP $P = (T; U, V)$ with $r_v = n_u - (t - 2)s_t - (v - 2)s_v > 0$ and $b^* = 2s_v/r_v$ integer is a 0-problem: $f(P) = x_u/u$, and there exists an optimal solution in which all elements of $V$ have size strictly greater than $x_u/u$.
\end{theorem}
\textbf{Proof.} Apply Algorithm \ref{alg:0prob_type2} to solve $P$. Because $r_v > 0$, it follows that $x > x_{\infty}$. Because $b^*$ is integer, the matrices $T$ and $V$ can be divided into $r_v/2$ identical $(T, V)$-pairs, $(A, B)$. $B$ contains $b^*$ rows from $V$ and $A$ contains $(v - 1)b^* + 1 = 2(v - 1)s_v/r_v + 1 = 2[(v - 1)s_v + r_v/2]/r_v = 2s_t/r_v$ rows from $T$, and those rows from $T$ contain $2n_t/r_v - 2n_v/r_v = 2n_u/r_v$ elements from $U$.

All elements of $U$ are set to $x_u/u$. All elements of $V$ are elements of $B$. Because $x > x_{\infty}$, we have
\[
  \frac{x_u}{u} > \frac{x_u/u + (v - 1)x_t - x_v}{(v - 1)(t - 1)}.
\]
Then all the conditions of Proposition \ref{prop:Vsub>x_u/u} are satisfied. The conclusions follow.
\hspace{\stretch{1}}$\blacksquare$

\begin{lemma} \label{lem:alg_0prob_type2_complexity} Algorithm \ref{alg:0prob_type2} has complexity $\Theta(n_t)$.
\end{lemma}
\textbf{Proof.} Assigning value $x_u/u$ to all elements of $U$ and to the vector $\boldsymbol{u}$ requires $n_u + 2n_u/r_v$ operations. Completing the pair $(A, B)$ using Algorithm \ref{alg:(T,V)-pair} requires between $2n_t/r_v$ and $6n_t/r_v$ operations by Lemma \ref{lem:alg_(T,V)-pair_complexity}. Inserting $A$ and $B$ into $T$ and $V$ requires $n_t + n_v$ operations. The total number of operations is bounded below by $2n_t$ and bounded above by $6n_t$.
\hspace{\stretch{1}}$\blacksquare$

\subsection{0-Problems of Type 2 in the Context of Extended Families}

Here we put 0-problems of type 2 into the context of extended families by investigating for what values of $x$ a 3M-DAP $P = (T; U, V)$ is a 0-problem of type 2. Define for any integer $b \ge 0$:
\[
  x_b \equiv \frac{[\lambda + \gamma (v - 1)t] ub + \gamma tu}{[(v - 1)t + u - v]b + t}.
\]
Observe that $(v - 1)t + u - v > 0$ and $\lambda < \gamma(u - v)$, so using mediants
\[
\begin{split}
  \lambda + \gamma(v - 1)t &< [(v - 1)t + u - v]\gamma\\
  \Leftrightarrow \frac{\lambda + \gamma (v - 1)t}{(v - 1)t + u - v} &< \gamma\\
  \Leftrightarrow \frac{[\lambda + \gamma (v - 1)t] ub + \gamma tu}{[(v - 1)t + u - v]b + t} &< \frac{[\lambda + \gamma (v - 1)t] u(b - 1) + \gamma tu}{[(v - 1)t + u - v](b - 1) + t}\\
  \Leftrightarrow x_b &< x_{b-1},
\end{split}
\]
\[
  x_0 = \frac{\gamma tu}{t} = \gamma u,
\]
and
\[
  \lim_{b \rightarrow \infty} x_b = \frac{[\lambda + \gamma (v - 1)t] u}{(v - 1)t + u - v} = x_{\infty}.
\]

\begin{lemma} \label{lem:b'<=>x=x_b'} Let $P = (T; U, V)$ be a 3M-DAP with $r_v > 0$ and $b - 1 \le 2s_v/r_v < b$ for some integer $b \ge 1$. Then $x = x_u \in (x_b, x_{b - 1}]$. If $2s_v/r_v = b$ for some integer $b \ge 1$, then $x = x_u = x_b$.
\end{lemma}
\textbf{Proof.} Observe that, if
\[
  b \ge \frac{s_v}{s_t - (v - 1)s_v},
\]
so that $s_tb \ge [(v - 1)b + 1]s_v$, then
\begin{align*}
  \frac{x}{u} &= \frac{s_tx_t - x_vx_v}{s_tt - s_vv}\\
  &\ge \frac{[(v - 1)b + 1]x_t - bx_v}{[(v - 1)b + 1]t - bv}\\
  &= \frac{[(v - 1)b + 1]t\gamma + b\lambda - bx}{[(v - 1)b + 1]t - bv}.
\end{align*}
Then
\[
  x \ge \frac{[(v - 1)b + 1]t\gamma + b\lambda}{[(v - 1)b + 1]t + (u - v)b}u = \frac{[(v - 1)t\gamma + \lambda]b + t\gamma}{[(v - 1)t + u - v]b + t}u = x_b.
\]
The conclusions follow.
\hspace{\stretch{1}}$\blacksquare$

\begin{theorem} \label{thm:0type2_x=x_b} A 3M-DAP $P = (T; U, V)$ is a 0-problem of type 2 if and only if $x = x_b$ for some $b \ge 1$.
\end{theorem}

To summarize what we have established to date: the 3M-DAP extended family $(t, u, v, \lambda, \gamma)$ includes all rational values of $x$ in the interval $(\lambda + \gamma v, \gamma u = x_0]$. The problem is a 0-problem whenever:
\begin{itemize}
\item 0-problem of type 1: $x \in (\lambda + \gamma v, x_{\infty}]$;
\item 0-problem of type 2: $x = x_b$ for $b = 0, 1, 2, \ldots$.
\end{itemize}
It remains to understand the problems for which $x_b < x < x_{b-1}$ for any $b \ge 1$.

\section{A Recursive Algorithm for 3M-DAPs}

Suppose that $x_b < x < x_{b - 1}$ for some $b > 0$. Then $n_u > (t - 2)s_t + (v - 2)s_v$ and recall that
\[
  r_v = n_u - (t - 2)s_t - (v - 2)s_v = ts_t - vs_v - (t - 2)s_t - (v - 2)s_v = 2s_t - 2(v - 1)s_v,
\]
so $r_v$ is positive and even.

\subsection{Constructing a reduced problem} \label{sec:reduced_problem}

We seek a solution with $(T, V)$ divided into $s_u^{\prime}$ maximal inseparable $b$-pairs and $s_v^{\prime}$ maximal inseparable $(b-1)$-pairs. By matching the number of $T$ rows and the number of $V$ rows, we have
\begin{align}
  [(v - 1)b + 1]s_u^{\prime} + [(v - 1)(b - 1) + 1]s_v^{\prime} &= s_t \label{eq:st}\\
  bs_u^{\prime} + (b - 1)s_v^{\prime} &= s_v \label{eq:sv}.
\end{align}
Then
\[
\begin{split}
  s_u^{\prime} + s_v^{\prime} &= s_t - (v - 1)s_v = \frac{1}{2}r_v\\
  s_v^{\prime} &= \frac{1}{2}br_v - s_v\\
  s_u^{\prime} &= s_v - \frac{1}{2}(b - 1)r_v.
\end{split}
\]

Now we can form the elements in $U$ from the $s_u^{\prime}$ maximal inseparable $b$-pairs into a matrix $U^{\prime}$, with each row corresponding to a maximal inseparable pair. Likewise, we can form the elements in $U$ from the $s_v^{\prime}$ maximal inseparable $(b-1)$-pairs into a matrix $V^{\prime}$, with each row corresponding to a maximal inseparable pair. The elements from $U^{\prime}$ and $V^{\prime}$ can be reorganized into the matrix $T^{\prime} = U$. To summarize:
\begin{alignat*}{2}
  s_u^{\prime} &= s_v - \frac{1}{2}(b - 1)r_v \quad \quad &&\text{number of $b$-pairs}\\
  u^{\prime} &= [(v - 1)b + 1] t - bv \quad \quad &&\text{number of $U$-elements in a $b$-pair}\\
  x_u^{\prime} &= [(v - 1)b + 1] x_t - bx_v \quad \quad &&\text{sum of $U$-elements in a $b$-pair}\\
  s_v^{\prime} &= \frac{1}{2}br_v - s_v \quad \quad &&\text{number of $(b - 1)$-pairs}\\
  v^{\prime} &= [(v - 1)(b - 1) + 1] t - (b - 1)v \quad \quad &&\text{number of $U$-elements in a $(b - 1)$-pair}\\
  x_v^{\prime} &= [(v - 1)(b - 1) + 1] x_t - (b - 1)x_v \quad \quad &&\text{sum of $U$-elements in a $(b - 1)$-pair}.
\end{alignat*}
We call the $P^{\prime} = (T^{\prime}, U^{\prime}, V^{\prime})$ problem the \textit{reduced problem} for problem $P = (T; U, V)$. It is smaller in the sense that the number of elements in $T^{\prime} = U$ is $n_u < n_t$, the number of elements in $T$.

\begin{proposition} \label{prop:reduced_problem_properties}
Let $P$ be a 3M-DAP with $x_{b} \le x < x_{b-1}$ for some $b \ge 1$. Let $P^{\prime}$ be the reduced problem for problem $P$. Then $P^{\prime}$ is a 3M-DAP. Further, 
\begin{enumerate}
\item[(i)] if $P$ is a $v1$ 3M-DAP, then $P^{\prime}$ is $v$-weighted with $v^{\prime} = u^{\prime} + 1$, otherwise $P^{\prime}$ is $u$-weighted;
\item[(ii)] if $t > 2$ or $v > 2$, then $P^{\prime}$ is strongly $u$-weighted;
\item[(iii)] $P^{\prime}$ is weakly $u$-weighted if and only if $P$ has $b > 1$ or $P$ is $u$-weighted and has $b = 1$.
\end{enumerate}
\end{proposition}
\textbf{Proof.} We show that the reduced problem satisfies the properties of a 3M-DAP. First, $t^{\prime} = u \ge 2$. If $v = 1$, then by assumption $n_v = s_v \le (t - 2)s_t$. Also, $r_v = 2s_t$, so $b^* = 2s_v/r_v \le t - 2$, so $\hat{b} = b \le t - 2$. Then $u^{\prime} = t - b \ge 2$, and $v^{\prime} = u^{\prime} + 1$. When $v \ge 2$, both $u^{\prime}$ and $v^{\prime}$ are nondecreasing functions of $b$, so are smallest when $b = 1$: then $u^{\prime} = v(t - 1) \ge 2$ and $v^{\prime} = t \ge 2$. 

Second, $x_t^{\prime}, x_u^{\prime}$, and $x_v^{\prime}$ are all clearly rational, so $\lambda^{\prime}$ and $\gamma^{\prime}$ are rational. Third, $s_t^{\prime} = s_u > 0$, and because $x_{b} \le x < x_{b-1}$, we must have $s_u^{\prime} > 0$ and $s_v^{\prime} \ge 0$. Fourth, we can confirm that:
\[
\begin{split}
  n_u^{\prime} + n_v^{\prime} &= [s_v - (b - 1)r_v/2]\{[(v - 1)b + 1]t - bv\}\\
  &\quad \quad \quad \quad \quad \quad + (br_v/2 - s_v)\{[(v - 1)(b - 1) + 1]t - (b - 1)v\}\\
  &= [(v - 1)t - v]s_v + tr_v/2\\
  &= (v - 1)ts_v - n_v + n_t - t(v - 1)s_v\\
  &= n_u = n_t^{\prime},
\end{split}
\]
and also that:
\[
\begin{split}
  s_u^{\prime}x_u^{\prime} + s_v^{\prime}x_v^{\prime} &= [s_v - (b - 1)r_v/2]\{[(v - 1)b + 1]x_t - bx_v\}\\
  &\quad \quad \quad \quad \quad \quad + (br_v/2 - s_v)\{[(v - 1)(b - 1) + 1]x_t - (b - 1)x_v\}\\
  &= [(v - 1)x_t - x_v]s_v + x_tr_v/2\\
  &= (v - 1)s_vx_t - s_vx_v + s_tx_t - (v - 1)s_vx_t\\
  &= s_tx_t - s_vx_v\\
  &= s_ux_u = s_t^{\prime}x_t^{\prime}.
\end{split}
\]
Finally, we have
\[
  \frac{x_u^{\prime}}{u^{\prime}} = \frac{[(v - 1)b + 1] x_t - bx_v}{[(v - 1)b + 1] t - bv} = \frac{(v - 1 + 1/b)x_t - x_v}{(v - 1 + 1/b)t - v}
\]
and
\[
  \frac{x_v^{\prime}}{v^{\prime}} = \frac{[(v - 1)(b - 1) + 1] x_t - (b - 1)x_v}{[(v - 1)(b - 1) + 1] t - (b - 1)v} = \frac{[v - 1 + 1/(b - 1)]x_t - x_v}{[v - 1 + 1/(b - 1)]t - v},
\]
where the second equality applies only when $b > 1$. If this is the case, then we can argue as follows. The mediant of $[x_v^{\prime}/(b - 1)]/[v^{\prime}/(b - 1)]$ and $x_v/v$ is $x_t/t$, so $x_v^{\prime}/v^{\prime} < x_t/t < x_v/v$. Likewise, the mediant of $[1/(b - 1) - 1/b]x_t/[1/(b - 1) - 1/b]t = x_t/t$ and $x_u^{\prime}/u^{\prime}$ is $x_v^{\prime}/v^{\prime}$, so we must have $x_u^{\prime}/u^{\prime} < x_v^{\prime}/v^{\prime}$, as desired. If $b = 1$, then $x_u^{\prime}/u^{\prime} = (vx_t - x_v)/(vt - v) < x_t/t = x_v^{\prime}/v^{\prime}$.

Now we investigate the properties of the parameter set for the reduced problem. First, consider a $v1$-problem. Then $v = 1$ and recall that $b^* \le t - 2$, so
\[
\begin{split}
  t^{\prime} &= u \ge 2;\\
  u^{\prime} &= [(v - 1)b + 1] t - bv  = t - b \ge 2, \quad \text{with equality if and only if } b = t - 2;\\
  v^{\prime} &= [(v - 1)(b - 1) + 1] t - (b - 1)v = t - b + 1 > 2.
\end{split}
\]
Then $v^{\prime} = u^{\prime} + 1$, so the problem is $v$-weighted.

Otherwise $b \ge 1$, $t \ge 2$, $u \ge 2$, and $v \ge 2$, so:
\[
\begin{split}
  t^{\prime} &= u \ge 2;\\
  u^{\prime} &= [(v - 1)b + 1] t - bv \ge (v - 2)b + 2 \ge 2, \quad \text{with equality if and only if } t = v = 2;\\
  v^{\prime} &= [(v - 1)(b - 1) + 1] t - (b - 1)v \ge (v - 2)(b - 1) + 2 \ge 2,\\
  &\quad \quad \quad \quad \quad \quad \quad \quad \quad \quad \quad \text{with equality if and only if } t = 2 \text{ and either } v = 2 \text{ or } b = 1.
\end{split}
\]
Then
\begin{equation} \label{eq:u'-v'}
  u^{\prime} - v^{\prime} = [(v - 1)b + 1] t - bv - [(v - 1)(b - 1) + 1] t + (b - 1)v = (v - 1)t - v \ge v - 2 \ge 0,
\end{equation}
with equality if and only if $t = v = 2$. Then the reduced problem is $u$-weighted and further is strongly $u$-weighted whenever $t > 2$ or $v > 2$.

When $b = 1$, $u^{\prime} = v(t - 1)$ and $v^{\prime} = t$ so $(v^{\prime} - 1)t^{\prime} = (t - 1)u$.
Then the $P^{\prime}$ is weakly $u$-weighted if and only if $P$ is $u$-weighted. Otherwise $b > 1$. In this case, because $t^{\prime} = u \ge 2$, we have
\[
\begin{split}
  (v^{\prime} - 1)t^{\prime} - u^{\prime} &\ge 2(v^{\prime} - 1) - u^{\prime}\\
  &= 2[(v - 1)(b - 1) + 1] t - 2(b - 1)v - 2 - [(v - 1)b + 1] t + bv\\
  &= [(v - 1)(b - 2) + 1] t - (b - 2)v - 2\\
  &= [(v - 1)t - v](b - 2) + t - 2 \ge 0,
\end{split}
\]
because $(v - 1)t - v \ge 0$ and $t \ge 2$. Then $P^{\prime}$ is weakly $u$-weighted.
\hspace{\stretch{1}}$\blacksquare$

\subsection{Reducing a 0-problem of type 2} \label{sec:reducing_0prob_type2}

When $x = x_b$ for some $b \ge 0$, we have previously shown a direct solution to the 3M-DAP $P = (T; U, V)$ that achieves the $x_u/u$ upper bound. However, it is possible to reduce such a problem.

Recall that $r_v = n_u - (t - 2)s_t - (v - 2)s_v = 2s_t - 2(v - 1)s_v$. If we reduce a problem with
\[
  b = b^* = \frac{2s_v}{r_v} = \frac{s_v}{s_t - (v - 1)s_v},
\]
then $s_v^{\prime} = 0$ and
\[
  (v - 1)b + 1 = \frac{s_t}{s_t - (v - 1)s_v},
\]
so
\[
  \frac{x_u^{\prime}}{u^{\prime}} = \frac{s_tx_t - s_vx_v}{s_tt - s_vv} = \frac{s_ux_u}{s_uu} = \frac{x_u}{u} = \frac{x_t^{\prime}}{t^{\prime}}.
\]
Then $x^{\prime} = \gamma^{\prime}u^{\prime} = x_0^{\prime}$. So certainly the reduced problem is a 0-problem of type 2 but perhaps more importantly its solution is immediate: all elements in both $U^{\prime}$ and $T^{\prime}$ can be set to $x_u^{\prime}/u^{\prime} = x_u/u$.

\subsection{The Recursive Algorithm}

Here we present the recursive algorithm for solving any 3M-DAP (Algorithm \ref{alg:recursive}). The algorithm solves any 0-problem directly using either Algorithm \ref{alg:0prob_type1} or Algorithm \ref{alg:0prob_type2}. Any 3M-DAP $P(x) = (T; U, V)$ with $x_b < x < x_{b-1}$ for some $b > 0$, i.e., that is not a 0-problem, is reduced to a smaller problem $P^{\prime} = (T^{\prime}; U^{\prime}, V^{\prime})$, which is itself then solved using the recursive algorithm. Given the solution to the smaller problem, $U = T^{\prime}$ and each row of $U^{\prime}$, respectively $V^{\prime}$, corresponds to a $b$-pair, respectively a $(b - 1)$-pair, in the original problem. For each such row, Algorithm \ref{alg:(T,V)-pair} is used to complete the $(T, V)$-pair, thereby producing a solution to the original problem.

\fontsize{8}{8}
\begin{algorithm}
\SetKwInOut{Input}{input}
\SetKwInOut{Output}{output}
\Input{$T$, $U$, $V$, empty matrices of dimensions $s_t \times t$, $s_u \times u$, and $s_v \times v$, respectively\\
$x_t$, $x_u$, $x_v$, the required row sums for each row of $T$, $U$, and $V$, respectively\\
$s_uu + s_vv = s_tt$\\
$s_ux_u + s_vx_v = s_tx_t$\\
}
\Output{the filled matrices $T$, $U$, $V$\\
the multiset of entries in $T$ is the union of the multiset of entries in $U$ and $V$\\
}
\BlankLine

\uIf{$x \le x_{\infty}$}{
	apply Algorithm \ref{alg:0prob_type1} to solve the $(T; U, V)$ problem\\
}
\uElseIf{$x = x_b$ for some integer $b \ge 0$}{
	apply Algorithm \ref{alg:0prob_type2} to solve the $(T; U, V)$ problem\\
}
\Else{
	$b \leftarrow \lceil s_v/[s_t - (v - 1)s_v] \rceil$\\
	form the reduced problem $(T^{\prime}, U^{\prime}, V^{\prime})$ as described in Section \ref{sec:reduced_problem}\\
	apply Algorithm \ref{alg:recursive} to solve the $(T^{\prime}, U^{\prime}, V^{\prime})$ problem\\
	$U \leftarrow T^{\prime}$\\
	\For{$1 \le i \le s_u^{\prime}$}{
		initialize a $b$-pair $(A, B)$\\
		set $\boldsymbol{u}$ to be the $i^{th}$ row of $U^{\prime}$\\
		apply Algorithm \ref{alg:(T,V)-pair} with inputs $A$, $x_t$, $B$, $x_v$, and $\boldsymbol{u}$ to complete the pair $(A, B)$\\
		insert $A$ into $T$ and $B$ into $V$\\
	}

	\uIf{$b = 1$}{
		insert $V^{\prime}$ into $T$\\
	}
	\uElse{
		\For{$1 \le j \le s_v^{\prime}$}{
			initialize a $(b - 1)$-pair $(A, B)$\\
			set $\boldsymbol{u}$ to be the $j^{th}$ row of $V^{\prime}$\\
			apply Algorithm \ref{alg:(T,V)-pair} with inputs $A$, $x_t$, $B$, $x_v$, and $\boldsymbol{u}$ to complete the pair $(A, B)$\\
			insert $A$ into $T$ and $B$ into $V$\\
		}
	}
}

\caption{Algorithm for solving a 3M-DAP}
\label{alg:recursive}
\end{algorithm}

\normalsize

\subsection{Reducing and transforming to standard form}

Algorithm \ref{alg:recursive} reduces any problem that is not a 0-problem. It will be helpful for analysis of the algorithm to transform the reduced problem to standard form. Recall that the reduced problem can only be transformed to standard form if $u^{\prime} > v^{\prime}$, so per \eqref{eq:u'-v'} this section focuses only on problems with $v \ge 2$ and either $t > 2$ or $v > 2$. 

For $x_b < x < x_{b-1}$, create the reduced problem using $b$- and $(b - 1)$-pairs. Then
\begin{alignat*}{2}
  u^{\prime}(b) &= [(v - 1)b + 1] t - bv\\
  x_u^{\prime}(b) &= [(v - 1)b + 1] x_t - bx_v\\
  t^{\prime}(b) &= u\\
  x_t^{\prime}(b) &= x_u\\
  v^{\prime}(b) &= [(v - 1)(b - 1) + 1] t - (b - 1)v\\
  x_v^{\prime}(b) &= [(v - 1)(b - 1) + 1] x_t - (b - 1)x_v.
\end{alignat*}
Then
\[
\begin{split}
  u^{\prime}(b) - v^{\prime}(b) &= (v - 1)t - v\\
  \lambda^{\prime}(x) = \lambda^{\prime}(x, b) = x_u^{\prime}(b) - x_v^{\prime}(b) &= (v - 1)x_t - x_v = \lambda + \gamma (v - 1)t - x\\
  \gamma^{\prime}(x) = \gamma^{\prime}(x, b) = \frac{x_t^{\prime}(b)}{t^{\prime}(b)} &= \frac{x_u}{u} = \frac{x}{u}.
\end{split}
\]
Because each of these is independent of $b$, it follows from \eqref{eqfunc} that when transforming from reduced to standard form problem the same transformation is applied for all members from the same extended family.

\begin{theorem} \label{thm:map_to_full_family} Given an extended 3M-DAP family $(t, u, v, \lambda, \gamma)$ with $v \ge 2$ and either $t > 2$ or $v > 2$, let $X$ be the subset of the family having $x > x_{\infty}$, and let $H(\cdot)$ be the map that transforms a 3M-DAP in $X$ to a 3M-DAP in standard form through reduction followed by standardization. Then $H(\cdot)$ is a bijection from $X$ to the complete 3M-DAP family $(u, \{u^{\prime}(b)\}_{b \ge 1}, \{v^{\prime}(b)\}_{b \ge 1}, 0, 1/2)$ of standard form problems.

In particular, for any integer $b \ge 1$, let $X_b$ be the subset of $X$ of problems having $x \in (x_b, x_{b - 1}]$. Then $H(\cdot)$ is a bijection from $X_b$ to the extended 3M-DAP family $(u, u^{\prime}(b), v^{\prime}(b), 0, 1/2)$ of standard form problems.
\end{theorem}
\textbf{Proof.} It is sufficient to prove the second claim. It is immediate that the given subset of problems will be mapped into the same extended family. As noted above, the transformation from reduced to standard form problem is independent of $b$:
\[
  h(y) = \frac{1}{2} \frac{[u^{\prime}(b) - v^{\prime}(b)]y - [x_u^{\prime}(b) - x_v^{\prime}(b)]}{[u^{\prime}(b) - v^{\prime}(b)]x_t^{\prime}(b)/t^{\prime}(b) - [x_u^{\prime}(b) - x_v^{\prime}(b)]} = \frac{1}{2} \frac{[(v - 1)t - v]y - [(v - 1)x_t - x_v]}{[(v - 1)t - v]x_u/u - [(v - 1)x_t - x_v]}.
\]

Let $\hat{x}(x)$ represent the value of $x_u$ in the standardized reduced problem. If $y_1 + \cdots + y_{u^{\prime}(b)} = x_u^{\prime}(b)$, then $\hat{x}(x) = h(y_1) + \cdots + h(y_{u^{\prime}(b)})$. Observe that
\[
\begin{split}
  &[u^{\prime}(b) - v^{\prime}(b)]x_u^{\prime}(b) - [x_u^{\prime}(b) - x_v^{\prime}(b)]u^{\prime}(b)\\
  &= x_v^{\prime}(b)u^{\prime}(b) - x_u^{\prime}(b)v^{\prime}(b)\\
  &= \{[(v - 1)(b - 1) + 1]x_t - (b - 1)x_v\}\{[(v - 1)b + 1]t - bv\}\\
  &\quad \quad \quad \quad \quad - \{[(v - 1)b + 1]x_t - bx_v\}\{[(v - 1)(b - 1) + 1]t - (b - 1)v\}\\
  &= - [(v - 1)(b - 1) + 1]bvx_t - [(v - 1)b + 1](b - 1)tx_v\\
  &\quad \quad \quad \quad \quad + [(v - 1)b + 1](b - 1)vx_t + [(v - 1)(b - 1) + 1]btx_v\\
  &= tx_v - vx_t\\
  &= t(x - \lambda - \gamma v),
\end{split}
\]
which is independent of $b$. Then $\hat{x}(x)$ will also be independent of $b$:
\begin{equation} \label{eq:xhat}
\begin{split}
 \hat{x}(x) &= \frac{1}{2} \frac{[u^{\prime}(b) - v^{\prime}(b)]x_u^{\prime}(b) - [x_u^{\prime}(b) - x_v^{\prime}(b)]u^{\prime}(b)}{[u^{\prime}(b) - v^{\prime}(b)]x_t^{\prime}(b)/t^{\prime}(b) - [x_u^{\prime}(b) - x_v^{\prime}(b)]}\\
 &= \frac{t}{2} \frac{x - \lambda - \gamma v}{[(v - 1)t - v]x/u - [\lambda + \gamma(v - 1)t - x]}\\
 &= \frac{t}{2} \frac{u(x - \lambda - \gamma v)}{[(v - 1)t + u - v]x - [\lambda + \gamma(v - 1)t]u}\\
 &= \frac{t}{2} \frac{u(x - \lambda - \gamma v)}{[(v - 1)t + u - v](x - x_{\infty})}\\
 &= \frac{t}{2} + \frac{t}{2} \frac{[(v - 1)t - v](\gamma u - x)}{[(v - 1)t + u - v](x - x_{\infty})}.
\end{split}
\end{equation}
As $x$ increases from $x_{\infty}$ to $\gamma u$, $\hat{x}(x)$ decreases from $+\infty$ to $t/2$. The following will be useful later in the proof:
\begin{equation} \label{eq:xhat_frac}
  \frac{\hat{x}(x)}{\hat{x}(x) - t/2}[(v - 1)t - v] = \frac{u(x - \lambda - \gamma v)}{\gamma u - x}.
\end{equation}
Now we are interested in the value taken by $\hat{x}(x_b)$:
\[
\begin{split}
  \gamma - \frac{x_b}{u} &= \gamma - \frac{[\lambda + \gamma (v - 1)t] b + \gamma t}{[(v - 1)t + u - v]b + t}\\
  &= \frac{[\gamma (u - v) - \lambda] b}{[(v - 1)t + u - v]b + t},
\end{split}
\]
and
\[
\begin{split}
  \frac{1}{u}(x_b - x_{\infty}) &= \frac{[\lambda + \gamma (v - 1)t] b + \gamma t}{[(v - 1)t + u - v]b + t} - \frac{\lambda + \gamma (v - 1)t}{(v - 1)t + u - v}\\
  &= \frac{[(v - 1)t + u - v][\lambda + \gamma (v - 1)t] b + [(v - 1)t + u - v] \gamma t}{\{[(v - 1)t + u - v]b + t\}[(v - 1)t + u - v]}\\
  &\qquad \qquad \qquad \qquad \qquad \qquad - \frac{\{[(v - 1)t + u - v]b + t\}[\lambda + \gamma (v - 1)t]}{\{[(v - 1)t + u - v]b + t\}[(v - 1)t + u - v]}\\
  &= \frac{[\gamma (u - v) - \lambda] t}{\{[(v - 1)t + u - v]b + t\}[(v - 1)t + u - v]}.
\end{split}
\]
Then
\[
  \hat{x}(x_b) = \frac{t}{2}  + \frac{t}{2} \frac{b}{t} [(v - 1)t - v] = \frac{1}{2} \{[(v - 1)t - v]b + t\} = \frac{u^{\prime}(b)}{2} = \frac{v^{\prime}(b + 1)}{2}.
\]
Thus for any problem in the original extended family with $x \in (x_b, x_{b - 1}]$, the corresponding standardized reduced problem is a member of the extended family $(u, u^{\prime}(b), v^{\prime}(b), 0, 1/2)$, as claimed.

For the converse, we must show that every 3M-DAP $\hat{P} = (\hat{T}; \hat{U}, \hat{V})$ a member of the extended family $(u, u^{\prime}(b), v^{\prime}(b), 0, 1/2)$, there is a unique $P = (T; U, V) \in X_b$ such that $H(P) = \hat{P}$.

Let $\hat{s}_t$, $\hat{s}_u$, and $\hat{s}_v$ be the number of rows in the matrices $\hat{T}$, $\hat{U}$, and $\hat{V}$, respectively, and let the required rowsums be $\hat{x}_t = u/2$ and $\hat{x}_u = \hat{x}_v = \hat{x}$, with $\hat{x} \in (v^{\prime}(b)/2, u^{\prime}(b)/2]$ and rational. Then we must set $s_u = \hat{s}_t$ and from \eqref{eq:st} and \eqref{eq:sv} we must set:
\begin{align*}
  s_t &= [(v - 1)b + 1]\hat{s}_u + [(v - 1)(b - 1) + 1]\hat{s}_v\\
  s_v &= b\hat{s}_u + (b - 1)\hat{s}_v.
\end{align*}
Then
\[
  s_tt - s_vv = \hat{s}_uu^{\prime}(b) + \hat{s}_vv^{\prime}(b) = \hat{s}_tt^{\prime}(b) = s_uu.
\]
Further, we must set $x_u = x$ to be the result of the inverse transformation to \eqref{eq:xhat} applied to $\hat{x}$. This will result in $x \in (x_b, x_{b-1}]$. Finally, we must set $x_v = x - \lambda$ and $x_t = \gamma t$. Now
\[
  \hat{s}_t u = 2\hat{x}[\hat{s}_u + \hat{s}_v] = u^{\prime}(b)\hat{s}_u + v^{\prime}(b)\hat{s}_v,
\]
so we can write
\[
\begin{split}
  \hat{s}_t u &= \frac{2\hat{x}}{2\hat{x} - t}[u^{\prime}(b)\hat{s}_u + v^{\prime}(b)\hat{s}_v] - \frac{t}{2\hat{x} - t}2\hat{x}[\hat{s}_u + \hat{s}_v]\\
  &= \frac{2\hat{x}}{2\hat{x} - t}\{[u^{\prime}(b) - t]\hat{s}_u + [v^{\prime}(b) - t]\hat{s}_v\}\\
  &= \frac{2\hat{x}}{2\hat{x} - t}[(v - 1)t - v][b\hat{s}_u + (b - 1)\hat{s}_v].
\end{split}
\]
Then, using \eqref{eq:xhat_frac}:
\[
\begin{split}
  \hat{s}_t u &= \frac{x - \lambda - \gamma v}{\gamma u - x}u[b\hat{s}_u + (b - 1)\hat{s}_v]\\
  \Leftrightarrow \hat{s}_t &= \frac{x - \lambda - \gamma v}{\gamma u - x}s_v\\
  \Leftrightarrow (\gamma u - x)s_u &= (x - \lambda - \gamma v)s_v = (x_v - \gamma v)s_v\\
  \Leftrightarrow \gamma s_u u + \gamma s_v v &= s_ux_u + s_vx_v\\
  \Leftrightarrow s_tx_t &= s_ux_u + s_vx_v.
\end{split}
\]
Then the problem $P = (T; U, V)$ is a 3M-DAP whose standardized reduction is $\hat{P} = (\hat{T}, \hat{U}, \hat{V})$ and is the unique such problem.
\hspace{\stretch{1}}$\blacksquare$

\section{Analysis of The Recursive Algorithm}

Let $P = (T; U, V)$ be a problem from the extended 3M-DAP family $(t, u, v, \lambda, \gamma)$. When appropriate, we indicate that $x_u = x$ for the problem by writing $P(x)$. Given problem $P$, applying Algorithm \ref{alg:recursive} results in a feasible solution. Let $g(P)$ be the size of the smallest piece in this solution.

%\begin{figure}
%\centering
%\begin{minipage}{.48\linewidth}
%  \includegraphics[width=\linewidth]{fig1.pdf}
%%  \captionof{figure}{(a)}
%%  \label{img1}
%\end{minipage}
%\hspace{.01\linewidth}
%\begin{minipage}{.48\linewidth}
%  \includegraphics[width=\linewidth]{Gchart.pdf}
%%  \captionof{figure}{Second caption}
%%  \label{img2}
%\end{minipage}
%    \caption{The left chart shows three standard reward functions. All three pass through the point $(1,1)$ and  at $u=1$ have gradient $0.5$. The right chart shows the corresponding tangent functions, which can be seen to satisfy the properties outlined in Lemmas \ref{Gconvex} and \ref{Gdiff}. (a) Budgeted linear reward function: $M(u) = \min\{u,1\}$. Here $y_{\max} = 1$ and $G(y) = 1 - y$. Note that $v(y) = 1$ for all $y$. (b) $M(u) = \min\{3u,0.5u+0.5,1.15\}$. Here $y_{\max} = 1.15$ and $G(y) = \max\{3 - 5y, (1.15 - y)/1.3\}$. Note that while $v(0.5) = 1.3$, $\bar{v}(0.5) = 0.2$. (c) $M(u) = \min\{u\sqrt{10},\sqrt{u}\}$. Here $y_{\max} = +\infty$ and $G(y) = \max\{\sqrt{10} - 10y, 1/4y\}$.}
%    \label{fig:one}
%\end{figure}

\begin{figure}[!htb]
  \center{\includegraphics[width=\textwidth]{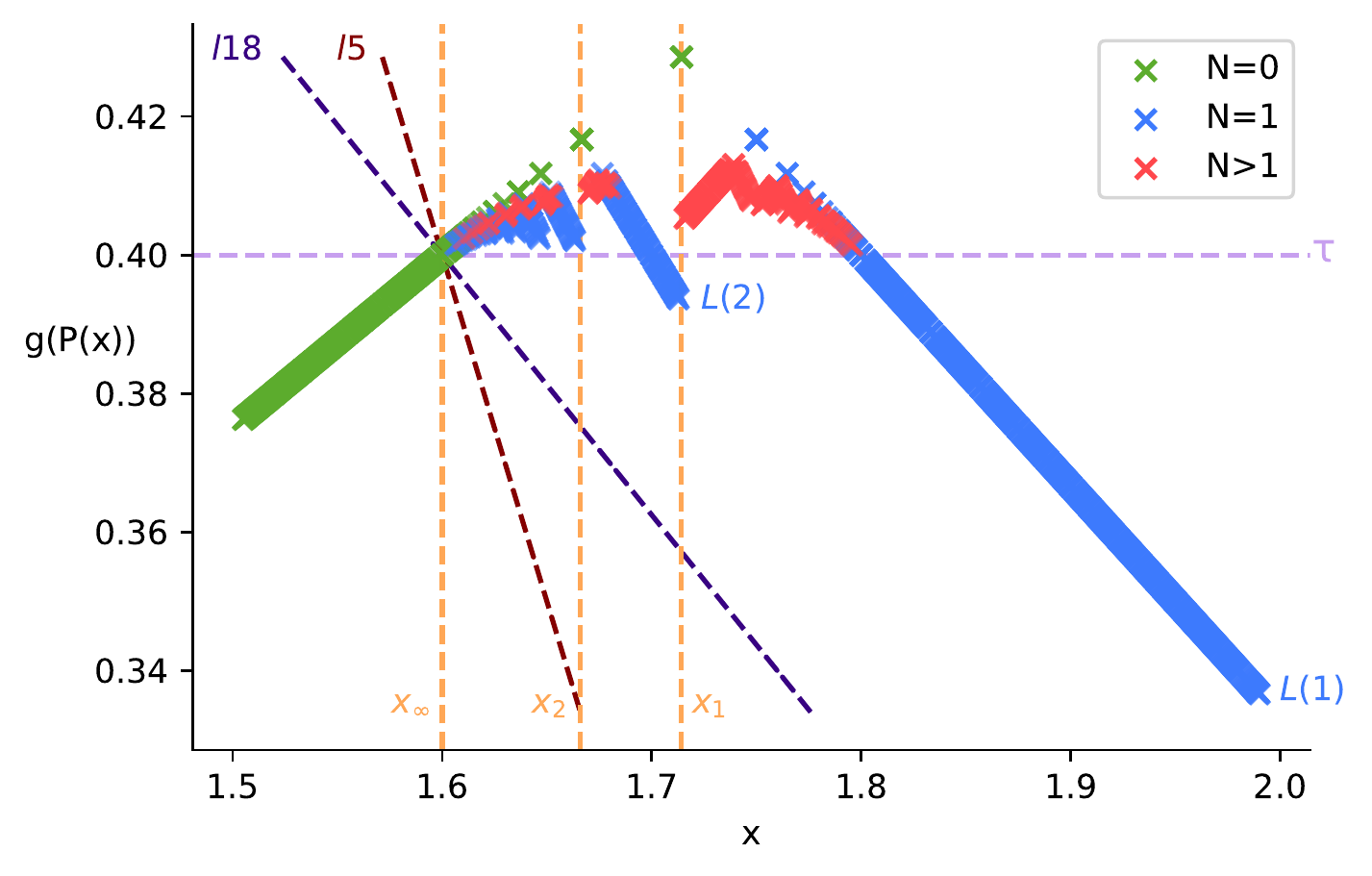}}
  \caption{Solution values produced by Algorithm \ref{alg:recursive} for problems in the extended fully-constrained-muffin-problem family with $n=3$, $F(t = 2, u = 4, v = 3, \lambda = 0, \gamma = \sfrac{1}{2})$. As we will show, these are in fact optimal values. $L(1)$ and $L(2)$ are the limits of $g(P(x))$ as $x$ converges to $x_0$ and $x_1$, respectively. The lines $l5$ and $l18$ represent the lower bounds described in Conjecture \ref{NprobLB_conj} and Proposition \ref{prop:Vsub>x_u/u}, respectively.}
  \label{fig:one}
\end{figure}

\begin{figure}[!htb]
  \center{\includegraphics[width=\textwidth]{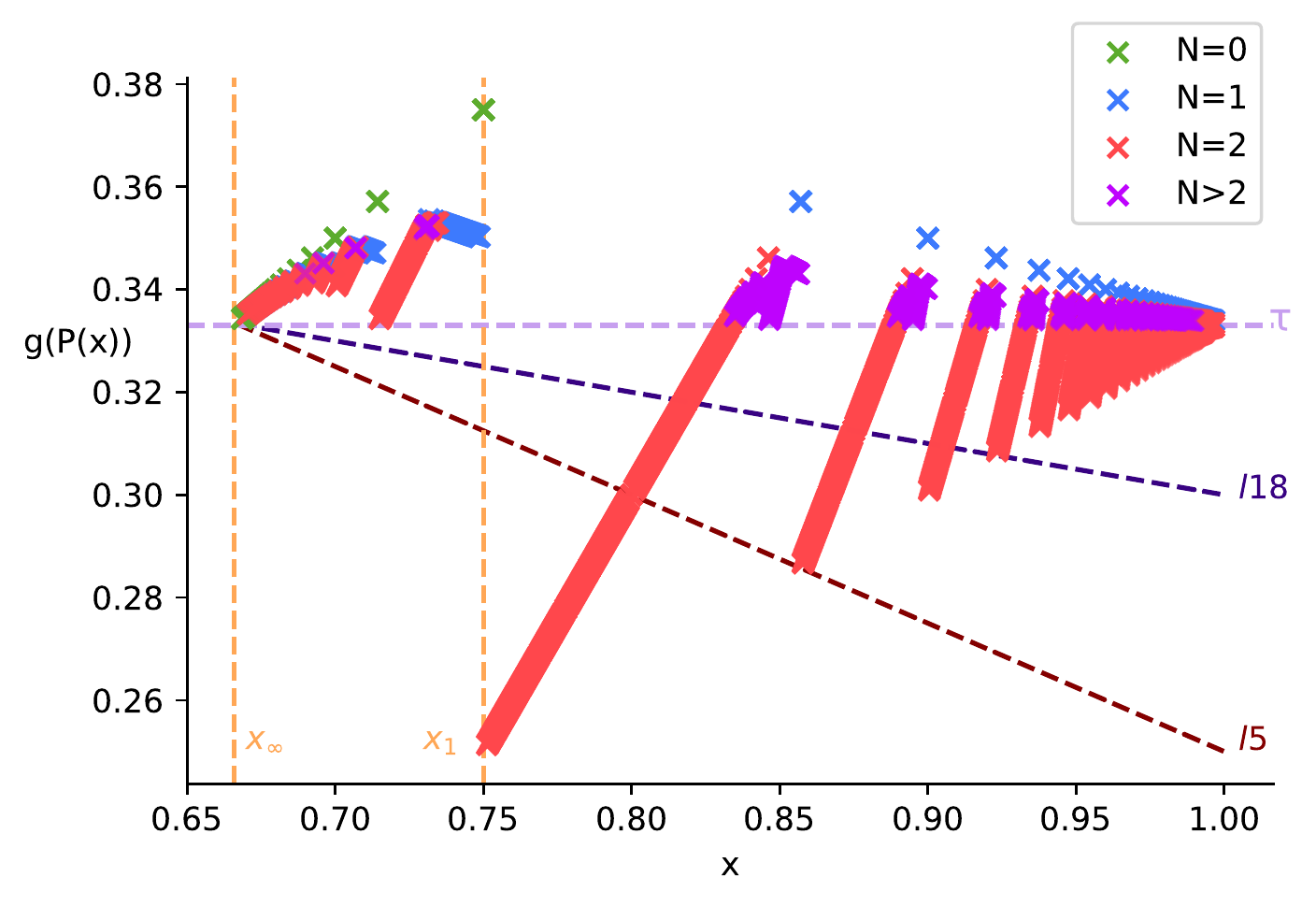}}
  \caption{Solution values produced by Algorithm \ref{alg:recursive} for problems in the extended 3M-DAP family $F(t = 2, u = 2, v = 6, \lambda = -3, \gamma = \sfrac{1}{2})$. As we will show, these are in fact optimal values. The lines $l5$ and $l18$ represent the lower bounds described in Conjecture \ref{NprobLB_conj} and Proposition \ref{prop:Vsub>x_u/u}, respectively.}
  \label{fig:two}
\end{figure}

\subsection{$N$-problems}

\begin{definition}
Algorithm \ref{alg:recursive} solves any 3M-DAP $P$ that is not a 0-problem by reducing it to a smaller 3M-DAP $P^{\prime}$. If this reduced problem is a 0-problem, then call the original problem a \emph{$1$-problem}. Then inductively define \emph{$N$-problems} for any $N \ge 2$ to be those 3M-DAPs for which it takes $N$ reductions in Algorithm \ref{alg:recursive} to reach a 0-problem. Any 3M-DAP must be an $N$-problem for some finite $N$.
\end{definition}

In this section we prove the following conjecture and establish lower bounds on the value of many 3M-DAPs. 

\begin{conjecture} \label{g=g'_conj} For any 3M-DAP $P$ that is not a 0-problem, $g(P) = g(P^{\prime})$.
\end{conjecture}

\subsection{Illustrative problem families}

Define
\begin{align*}
  \tau &= \frac{x_{\infty}}{u} = \frac{(v - 1)x_t + x_u - x_v}{(v - 1)t + u - v}, \quad \quad \text{and}\\
  L(b) &= \gamma - b\frac{x_{b - 1} - \lambda - \gamma v}{[(v - 1)t - v]b + t}.
\end{align*}
We will show below that $L(b)$ is the limit of $g(P(x))$ as $x$ converges to $x_{b - 1}$ from below. 

By Proposition \ref{prop:reduced_problem_properties} most $N$-problems with $N > 1$ reduce to $u$-weighted problems. Because fully-constrained muffin problems are $u$-weighted, we can use such a problem family to illustrate the typical properties of $u$-weighted problem families. 

Figure \ref{fig:one} shows the problem values for solutions computed by Algorithm \ref{alg:recursive} for the extended fully-constrained-muffin-problem family with $n=3$. This illustrates that for $u$-weighted problems, all $N$-problems with $N > 1$ have $g(P) > \tau$ and this is also true for many $1$-problems. However, there are some $1$-problems for which this does not hold: for these we see that $g(P) > L(1)$. 

The figure includes the lines $l5$ and $l18$ that represent the lower bounds on the problem value described in Conjecture \ref{NprobLB_conj} and Proposition \ref{prop:Vsub>x_u/u}, respectively. The former was used to show that any optimal solution must be constructed from $k$-pairs. The latter was used to show that the elements of $V$ used in completing a $k$-pair must be larger than $x_u/u$. These are critical to establishing the optimality of the solutions produced by Algorithm \ref{alg:recursive}. Figure \ref{fig:one} illustrates that these lower bounds apply for any $u$-weighted problem.

For an $N$-problem with $N > 1$ and $v = 1$, the reduced problem is $v$-weighted. Figure \ref{fig:two}
shows the problem values for solutions computed by Algorithm \ref{alg:recursive} for the extended 3M-DAP family $F(2, 2, 6, -3, \sfrac{1}{2})$ on the range $x > x_{\infty}$. This illustrates that the lower bounds of Conjecture \ref{NprobLB_conj} and Proposition \ref{prop:Vsub>x_u/u} don't necessarily apply for a $v$-weighted problem with $N = 2$ and $b = 1$. Such cases require special treatment to demonstrate the optimality of the solutions produced by Algorithm \ref{alg:recursive}.

We will prove Conjecture \ref{g=g'_conj} using an inductive argument. The preceding observations are reflected in the inductive hypothesis.

As shown earlier, for $P$ a $v$-weighted problem with $N = 2$ and $b = 1$, it may be the case that
\[
  f(P) \le \frac{(v - 1)x_t - x_v}{(v - 1)t - v}.
\]
In this case, we are not guaranteed that an optimal solution must consist of $k$-pairs. Consider the problem $(T, U, V) = ((7, 4, 1), (4, 2, 0.8), (1, 6, 3.8))$. Table \ref{tab:ex_soln} shows the solution produced by Algorithm \ref{alg:recursive}. Because $b = 1$, this is comprised of $1$- and $0$-pairs, in this case one of each. The optimal solution is
\[
  f(P) = 0.3 = \frac{(v - 1)x_t - x_v}{(v - 1)t - v},
\]
i.e., the strict lower bound of Conjecture \ref{NprobLB_conj} does not apply. Then it may be possible for there to be an optimal solution in which the $(T, V)$-pair that contains the one row of $V$ contains less than $v = 6$ rows of $T$. Table \ref{tab:alt_soln} shows just such an alternative optimal solution: the $(T, V)$-pair containing the row of $V$ contains only 5 rows of $T$.

\begin{table}
\centering
\begin{tabular}{ccc}
$T = $
\begin{tabular}{cc}
  \hline
  \cellcolor{yellow}0.3& \cellcolor[rgb]{0,1,1}0.7\\
  \cellcolor{yellow}0.3& \cellcolor[rgb]{0,1,1}0.7\\
  \cellcolor[rgb]{1,0.647,0}0.4& \cellcolor[rgb]{0,1,1}0.6\\
  \cellcolor[rgb]{1,0.647,0}0.4& \cellcolor[rgb]{0,1,1}0.6\\
  \cellcolor[rgb]{1,0.647,0}0.4& \cellcolor[rgb]{0,1,1}0.6\\
  \cellcolor[rgb]{1,0.647,0}0.4& \cellcolor[rgb]{0,1,1}0.6\\
  \cellcolor{green}0.5& \cellcolor{green}0.5\\
  \hline
\end{tabular}
&
$T^{\prime} = U = $
\begin{tabular}{cc}
  \hline
  \cellcolor{yellow}0.3& \cellcolor{green}0.5\\
  \cellcolor{yellow}0.3& \cellcolor{green}0.5\\
  \rowcolor[rgb]{1,0.647,0}0.4& 0.4\\
  \rowcolor[rgb]{1,0.647,0}0.4& 0.4\\
  \hline
\end{tabular}
&
\begin{tabular}{c}
$V = $ \\
\begin{tabular}{cccccc}
  \hline
  \rowcolor[rgb]{0,1,1} 0.7& 0.7& 0.6& 0.6& 0.6& 0.6\\
  \hline
\end{tabular}
\end{tabular}
\\
\rule{0pt}{10ex}
&
\begin{tabular}{c}
$T^{\prime \prime} = U^{\prime} = $ \\
\begin{tabular}{cccccc}
  \hline
  \rowcolor[rgb]{1,0.647,0} \cellcolor{yellow}0.3& \cellcolor{yellow}0.3& 0.4& 0.4& 0.4& 0.4\\
  \hline
\end{tabular}
\end{tabular}
&
$V^{\prime} = $
\begin{tabular}{cc}
  \hline
  \rowcolor{green} 0.5& 0.5\\
  \hline
\end{tabular}
\\
\rule{0pt}{10ex}
&
$U^{\prime \prime} = $
\begin{tabular}{cc}
  \hline
  \rowcolor{yellow} 0.3& 0.3\\
  \hline
\end{tabular}
&
$V^{\prime \prime} = $
\begin{tabular}{cc}
  \hline
  \rowcolor[rgb]{1,0.647,0} 0.4& 0.4\\
  \rowcolor[rgb]{1,0.647,0} 0.4& 0.4\\
  \hline
\end{tabular}
\\
\end{tabular}
\caption{Solution to problem $(T, U, V) = ((7, 4, 1), (4, 2, 0.8), (1, 6, 3.8))$ by Algorithm \ref{alg:recursive}}
\label{tab:ex_soln}
\end{table}

\begin{table}
\centering
\begin{tabular}{ccc}
$T = $
\begin{tabular}{cc}
  \hline
  \cellcolor[rgb]{0,1,1}0.5& \cellcolor[rgb]{0,1,1}0.5\\
  \cellcolor{yellow}0.3& \cellcolor[rgb]{0,1,1}0.7\\
  \cellcolor{yellow}0.3& \cellcolor[rgb]{0,1,1}0.7\\
  \cellcolor{yellow}0.3& \cellcolor[rgb]{0,1,1}0.7\\
  \cellcolor{yellow}0.3& \cellcolor[rgb]{0,1,1}0.7\\
  \cellcolor{green}0.5& \cellcolor{green}0.5\\
  \cellcolor{green}0.5& \cellcolor{green}0.5\\
  \hline
\end{tabular}
&
$T^{\prime} = U = $
\begin{tabular}{cc}
  \hline
  \cellcolor{yellow}0.3& \cellcolor{green}0.5\\
  \cellcolor{yellow}0.3& \cellcolor{green}0.5\\
  \cellcolor{yellow}0.3& \cellcolor{green}0.5\\
  \cellcolor{yellow}0.3& \cellcolor{green}0.5\\
  \hline
\end{tabular}
&
$V = $
\begin{tabular}{cccccc}
  \hline
  \rowcolor[rgb]{0,1,1} 0.5& 0.5& 0.7& 0.7& 0.7& 0.7\\
  \hline
\end{tabular}
\\
\rule{0pt}{10ex}
&
\begin{tabular}{c}
$U^{\prime} = $ \\
\begin{tabular}{cccc}
  \hline
  \rowcolor{yellow} 0.3& 0.3& 0.3& 0.3\\
  \hline
\end{tabular}
\end{tabular}
&
$V^{\prime} = $
\begin{tabular}{cc}
  \hline
  \rowcolor{green} 0.5& 0.5\\
  \rowcolor{green} 0.5& 0.5\\
  \hline
\end{tabular}
\\
\end{tabular}
\caption{Alternative solution to problem $(T, U, V) = ((7, 4, 1), (4, 2, 0.8), (1, 6, 3.8))$}
\label{tab:alt_soln}
\end{table}

\subsection{Inductive hypothesis}
 
We establish Conjecture \ref{g=g'_conj} by induction on $N$. The inductive hypothesis holds for $N$ if for any $N$-problem $P$, the solution produced by Algorithm \ref{alg:recursive} has the following properties:
\begin{enumerate}
\item if $N > 0$, then $g(P) = g(P^{\prime})$;
\item[2.1] if $N = 0$, then $g(P) = x_u/u$;
\item[2.2] if $N = 1$ and $v > 1$, then $g(P) > \min\left\{\tau, L(1)\right\}$;
\item[2.3] if $N > 1$ and $v > 1$ (except when $N = 2$, $u < v$, and $b = 1$), then $g(P) > \tau$;
\item[3.1] all elements of $U$ have size less than $x_v/v$;
\item[3.2] all elements of $V$ have size greater than $x_u/u$.
\end{enumerate}

If the inductive hypothesis holds, conclude that the smallest element of $T$ is in $U$, i.e., that $g(P) \in U$; and the largest element of $T$ is in $V$.

\begin{proposition} \label{prop:IH_N=0} The inductive hypothesis holds when $N = 0$, i.e., for any 0-problem.
\end{proposition}
\textbf{Proof.} For any 0-problem, items 1., 2.2, and 2.3 of the hypothesis do not apply. The remaining items hold by Theorem \ref{thm:0type1} for 0-problems of type 1 and by Theorem \ref{thm:0type2} for 0-problems of type 2.
\hspace{\stretch{1}}$\blacksquare$

Given an $N$-problem $P$ with $N \ge 1$, we can assume that the inductive hypothesis holds for any $(N - 1)$-problem and, in particular, for the problem $P^{\prime}$ to which $P$ is reduced by Algorithm \ref{alg:recursive} Our plan of attack is:
\begin{enumerate}
\item[1.] When $v > 1$ and $b > 1$, use the lower bound on $g(P^{\prime})$ assumed in the inductive hypothesis to establish a lower bound on $g(P^{\prime})$ in terms of the parameters of $P$;
\item[2.1] When $v > 1$ and $b > 1$, use this lower bound to establish item 3.2 of the inductive hypothesis for $P$; when $v = 1$ or $b = 1$, use the inductive hypothesis for $P^{\prime}$ to directly establish item 3.2 of the inductive hypothesis for $P$;
\item[2.2] Use the inductive hypothesis for $P^{\prime}$ to directly establish item 3.1 of the inductive hypothesis for $P$;
\item[3.] Use these items to establish item 1 of the inductive hypothesis for $P$;
\item[4.] Use this item to establish item 2.2 or 2.3 as appropriate of the inductive hypothesis for $P$.
\end{enumerate}

\subsection{Lower bounds on the reduced problem $P^{\prime}$ in terms of the parameters of the original problem $P$}

\begin{lemma} \label{lem:Y'_1prob} Let $P(x)$ be a 1-problem with $v > 1$ and $x_b < x < x_{b-1}$ for some $b > 0$. Then
\[
  g(P^{\prime}) > \frac{x_u/u + (v - 1)x_t - x_v}{(v - 1)(t - 1)} = \frac{(v - 1)t\gamma + \lambda - x + x/u}{(v - 1)t + u - v - u + 1}.
\]
If further $(v - 1)t - v > 0$, then
\[
  g(P^{\prime}) > \frac{(v - 1)x_t - x_v}{(v - 1)t - v}.
\]
\end{lemma}
\textbf{Proof.} Form the reduced problem $P^{\prime}$ using $b$- and $(b - 1)$-pairs. Then this is a 0-problem, so
\[
  g(P^{\prime}) = \frac{x_u^{\prime}}{u^{\prime}} = \frac{[(v - 1)b + 1]x_t - bx_v}{[(v - 1)b + 1]t - bv} = \frac{[(v - 1)t\gamma + \lambda]b + t\gamma - bx}{[(v - 1)t + u - v]b + t - bu}.
\]
This is a mediant of $x_b/u$ and $(bx)/(bu) = x/u$ with a negative weight on the latter, so because $x > x_b$, it follows that $g(P^{\prime}) < x_b/u$. Now consider
\[
  [(v - 1)t + u - v](b - 1) + t - (b - 1)u - 1 = [(v - 1)t - v](b - 1) + t - 1.
\]
This is positive: $v \ge 2$ and $b \ge 1$, so this is at least $t - 1 > 0$. Then
\[
  \frac{[(v - 1)t\gamma + \lambda](b - 1) + t\gamma - (b - 1)x - x/u}{[(v - 1)t + u - v](b - 1) + t - (b - 1)u - 1}
\]
is a mediant of $x_{b - 1}/u$ and $x/u$ with a negative weight on the latter, so on noting that the denominator is positive and $x < x_{b - 1}$, it follows that this is greater than $x_{b - 1}/u > x_b/u > g(P^{\prime})$. Now taking the mediant of this with a negative weight and $g(P^{\prime})$ results in the desired lower bound. The final result then follows directly on noting that $x > x_{\infty}$.
\hspace{\stretch{1}}$\blacksquare$

\begin{lemma} \label{lem:Y'_2prob} Let $P(x)$ be a 2-problem with $x_b < x < x_{b-1}$ for some $b > 0$. Suppose that for the reduced problem $P^{\prime}$, $g(P^{\prime}) > L^{\prime}(1)$. If $v = 1$, or $v \ge 2$ and $b \ge 2$, or $u \ge v \ge 2$ and $b = 1$, then $g(P^{\prime}) > \tau$.
\end{lemma}
\textbf{Proof.} The reduced problem $P^{\prime}$ is formed using $b$- and $(b - 1)$-pairs. By assumption,
\begin{align*}
  g(P^{\prime}) > L^{\prime}(1) = \gamma^{\prime} - \frac{\gamma^{\prime}(u^{\prime} - v^{\prime}) - \lambda^{\prime}}{(t^{\prime} - 1)v^{\prime}} &= \frac{x}{u} - \frac{[(v - 1)t + u - v]x/u - \lambda - \gamma(v - 1)t}{(u - 1)v^{\prime}}\\
 &= \frac{[(u - 1)v^{\prime} - (v - 1)t - u + v]x + \lambda u + \gamma(v - 1)tu}{u(u - 1)v^{\prime}}.
\end{align*}
Now $v^{\prime} = [(v - 1)(b - 1) + 1]t - (b - 1)v = [(v - 1)t - v](b - 1) + t$. If $v = 1$, then $v^{\prime} = t + 1 - b > 1$ by Lemma \ref{lem:v1} and the coefficient of $x$ in the numerator above is $(u - 1)(v^{\prime} - 1) > 0$. Otherwise, $v \ge 2$, so if $b \ge 2$, the coefficient of $x$ in the numerator above is greater than or equal to:
\[
  (u - 1)[(v - 1)t - v] + (u - 1)t - (v - 1)t - u + v = (u - 2)[(v - 1)t - v] + (u - 1)t - u \ge 0.
\]
If $b = 1$, this coefficient is $(u - v)(t - 1)$, which is nonnegative when $u \ge v$. In these cases the lower bound on $g(P^{\prime})$ is nondecreasing in $x$ and $x > x_{\infty}$, so
\[
  g(P^{\prime}) > \frac{x_{\infty}}{u} - \frac{[(v - 1)t + u - v]x_{\infty}/u - \lambda - \gamma(v - 1)t}{(u - 1)v^{\prime}} = \frac{x_{\infty}}{u} = \tau.
\]
\hspace{\stretch{1}}$\blacksquare$

\begin{lemma} \label{lem:Y'_Nprob} Let $P$ be an $N$-problem with $N \ge 2$. Suppose that for the reduced problem $P^{\prime}$, $g(P^{\prime}) > \tau^{\prime}$. Then $g(P^{\prime}) > \tau$.
\end{lemma}
\textbf{Proof.} We have
\[
\begin{split}
  g(P^{\prime}) > \frac{(v^{\prime} - 1)x_t^{\prime} + x_u^{\prime} - x_v^{\prime}}{(v^{\prime} - 1)t^{\prime} + u^{\prime} - v^{\prime}} &= \frac{(v^{\prime} - 1)x_u + (v - 1)x_t - x_v}{(v^{\prime} - 1)u + (v - 1)t - v}\\
  &\ge \frac{(v - 1)x_t + x_u - x_v}{(v - 1)t + u - v},
\end{split}
\]
because $v^{\prime} - 1 \ge 1$ and one of the following holds:
\begin{itemize}
\item $x_u/u > x_{\infty}/u > [(v - 1)x_t - x_v]/[(v - 1)t - v]$, if $(v - 1)t - v > 0$ (when $v > 2$, or $v = 2$ and $t > 2$);
\item $(v - 1)x_t - x_v < 0$, if $(v - 1)t - v = 0$ (when $t = v = 2$);
\item $x_u/u < x_v/v$, if $v = 1$.
\end{itemize}
\hspace{\stretch{1}}$\blacksquare$

The previous three lemmas more than suffice to establish the following proposition.

\begin{proposition} \label{prop:Y'_LB} Let $P(x)$ be an $N$-problem with $N > 0$, $v > 1$, and $x_b < x < x_{b-1}$ for some $b > 1$. Assume the inductive hypothesis holds for any $(N - 1)$-problem. Then for the reduced problem $P^{\prime}$:
\[
  g(P^{\prime}) > \frac{x_u/u + (v - 1)x_t - x_v}{(v - 1)(t - 1)}.
\]
\end{proposition}
\textbf{Proof.} First, note that
\[
  \tau = \frac{x_{\infty}}{u} = \frac{(v - 1)x_t + x_u - x_v}{(v - 1)t + u - v} > \frac{x_u/u + (v - 1)x_t - x_v}{(v - 1)(t - 1)},
\]
because $u > 1$ and one of the following holds:
\begin{itemize}
\item $x_u/u > x_{\infty}/u > [(v - 1)x_t - x_v]/[(v - 1)t - v]$, if $(v - 1)t - v > 0$ (when $v > 2$, or $v = 2$ and $t > 2$);
\item $(v - 1)x_t - x_v < 0$, if $(v - 1)t - v = 0$ (when $t = v = 2$).
\end{itemize}
	
If $N = 1$, then apply Lemma \ref{lem:Y'_1prob}. If $N \ge 2$, then by Proposition \ref{prop:reduced_problem_properties} and the inductive hypothesis, the assumptions of at least one of Lemma \ref{lem:Y'_2prob} or Lemma \ref{lem:Y'_Nprob} hold, and hence so does the conclusion.
\hspace{\stretch{1}}$\blacksquare$

\subsection{Bounds on the sizes of the elements of $U$ and $V$}

\begin{proposition} \label{lem:U<x_v/v} Let $P(x)$ be an $N$-problem with $x_b < x < x_{b-1}$ for some $N > 0$ and $b > 0$. Assume the inductive hypothesis holds for any $(N - 1)$-problem. Then, in a solution provided by Algorithm \ref{alg:recursive}, all elements of $U$ have size less than $x_v/v$. If further $t \le b$, then all elements of $U$ have size less than $x_t/t$.
\end{proposition}
\textbf{Proof.} Form the reduced problem $P^{\prime}$. Then $U = T^{\prime}$ and $T^{\prime}$ is comprised of the elements of $U^{\prime}$ and $V^{\prime}$. By the inductive hypothesis, all elements of $U^{\prime}$ are less than $x_v^{\prime}/v^{\prime} \le x_t/t < x_v/v$.

Also by the inductive hypothesis, all elements of $V^{\prime}$ are greater than $x_u^{\prime}/u^{\prime}$, so by considering a row of $V^{\prime}$ in which all elements but one have size $x_u^{\prime}/u^{\prime}$, we can conclude that all elements of $V^{\prime}$ are less than
\[
  x_v^{\prime} - (v^{\prime} - 1)\frac{x_u^{\prime}}{u^{\prime}} = \frac{u^{\prime}x_v^{\prime} - (v^{\prime} - 1)x_u^{\prime}}{u^{\prime}}.
\]
Now $v^{\prime} - 1 = u^{\prime} - (v - 1)t + v - 1 = u^{\prime} - (t - 1)(v - 1)$. Then the numerator of our strict upper bound on the elements of $V^{\prime}$ is
\[
\begin{split}
  u^{\prime}(x_v^{\prime} - x_u^{\prime}) + (t - 1)(v - 1)x_u^{\prime} &= u^{\prime}x_v - u^{\prime}(v - 1)x_t + (t - 1)(v - 1)x_u^{\prime}\\
  &= (v - 1)[(t - 1)(v - 1)b + t - 1 - (v - 1)bt - t + bv]x_t\\
  &\quad \quad \quad + [(v - 1)bt + t - bv - (t - 1)(v - 1)b]x_v\\
  &= (v - 1)(b - 1)x_t + (t - b)x_v.
\end{split}
\]
Then the strict upper bound on the elements of $V^{\prime}$ is
\[
  \frac{(v - 1)(b - 1)x_t + (t - b)x_v}{(v - 1)(b - 1)t + (t - b)v},
\]
and this is at most $x_v/v$ because $x_t/t < x_v/v$, $b \ge 1$, and if $v = 1$, then $b \le t - 2$ by Lemma \ref{lem:v1}. Indeed, if $t \le b$, then this is at most $x_t/t$.\hspace{\stretch{1}}$\blacksquare$

\begin{proposition} \label{lem:V>x_u/u} Let $P(x)$ be an $N$-problem with $x_b < x < x_{b-1}$ for some $N > 0$ and $b > 0$. Assume the inductive hypothesis holds for any $(N - 1)$-problem. Then, in a solution provided by Algorithm \ref{alg:recursive}, all elements of $V$ have size greater than $x_u/u$. If $v = 1$, then all elements of $V$ have size $x_v/v$. If $v > 1$ and $b = 1$, then all elements of $V$ have size greater than $x_t/t$.
\end{proposition}
\textbf{Proof.} If $v = 1$, the result is immediate. If $v > 1$ and $b = 1$, then elements of $V$ appear only in 1-pairs in $T$. Indeed each element of $V$ appears in a row of $T$ with $t - 1$ elements of $U^{\prime}$. By the inductive hypothesis, all elements of $U^{\prime}$ are less than $x_v^{\prime}/v^{\prime} = x_t/t$. Then all elements of $V$ have size greater than
\[
  x_t - (t - 1)\frac{x_t}{t} = \frac{x_t}{t}.
\]

Otherwise, $v > 1$ and $b > 1$. Recall that Algorithm \ref{alg:recursive} divides $(T, V)$ into maximal, inseparable $(T, V)$-pairs, with one pair for each row in $U^{\prime}$ and one pair for each row in $V^{\prime}$. Because $v > 1$, by Proposition \ref{prop:Y'_LB}, the smallest element in any row of $U^{\prime}$ or $V^{\prime}$ is
\[
  g(P^{\prime}) > \frac{x_u/u + (v - 1)x_t - x_v}{(v - 1)(t - 1)}.
\]
Then when Algorithm \ref{alg:(T,V)-pair} is used to fill the remaining elements of any of the $(T, V)$-pairs, Proposition \ref{prop:Vsub>x_u/u} applies, and the conclusion follows.
\hspace{\stretch{1}}$\blacksquare$

\subsection{The value of a problem is the value of its reduced problem}

\begin{proposition} \label{prop:g=g'} Let $P$ be an $N$-problem with $N > 0$. Assume the inductive hypothesis holds for any $(N - 1)$-problem, and in particular for the reduced problem $P^{\prime}$. Then, in a solution provided by Algorithm \ref{alg:recursive}, $g(P) = g(P^{\prime})$.
\end{proposition}
\textbf{Proof.} By Proposition \ref{lem:V>x_u/u}, any element of $V$ has size greater than $x_u/u > x_u^{\prime}/u^{\prime} \ge g(P^{\prime})$. Further, all elements of $U = T^{\prime}$ have size at least $g(P^{\prime})$, and at least one such element has size $g(P^{\prime})$. The conclusion follows.
\hspace{\stretch{1}}$\blacksquare$

\subsection{Lower bounds}

The majority of this section focuses on lower bounds for 1-problems because this is the difficult case. So let $P$ be a $1$-problem and let $P^{\prime}$ be the corresponding reduced problem, which must be a 0-problem. By Propositions \ref{prop:IH_N=0} and \ref{prop:g=g'}, $g(P) = g(P^{\prime})$.

Write $d = (v - 1)t - v$, and for any $b \ge 1$ define
\begin{align*}
  L(b) &= \gamma - b\frac{x_{b-1} - \lambda - \gamma v}{[(v - 1)t - v]b + t} = \gamma - b\frac{x_{b-1} - \lambda - \gamma v}{db + t}\\
  M(b) &= \frac{b[d(b - 1) + t]}{[(d + u)(b - 1) + t][db + t]}.
\end{align*}
Now
\[
  x_{b-1} = \frac{[\lambda + \gamma (v - 1)t] u(b - 1) + \gamma tu}{[(v - 1)t + u - v](b - 1) + t} = \frac{[\lambda + \gamma (d + v)] u(b - 1) + \gamma tu}{(d + u)(b - 1) + t},
\]
so
\[
  x_{b - 1} - \lambda - \gamma v = [\gamma (u - v) - \lambda] \frac{d(b - 1) + t}{(d + u)(b - 1) + t}.
\]
Then
\[
  L(b) = \gamma - [\gamma(u - v) - \lambda]\frac{b[d(b - 1) + t]}{[(d + u)(b - 1) + t][db + t]} = \gamma - [\gamma(u - v) - \lambda]M(b).
\]
The limiting behavior is given by:
\begin{align*}
  \lim_{b \rightarrow \infty} M(b) &= \frac{1}{d + u}\\
  \lim_{b \rightarrow \infty} L(b) &= \gamma - \frac{\gamma(u - v) - \lambda}{d + u} = \tau.
\end{align*}

\begin{lemma} \label{lem:L(b)_LB} Let $P(x)$ be a 1-problem with $x_b < x < x_{b - 1}$ for some $b > 0$. Then $g(P(x)) > L(b)$ and $g(P(x)) \searrow L(b)$ as $x \nearrow x_{b - 1}$.
\end{lemma}
\textbf{Proof.} We have:
\[
\begin{split}
  g(P(x)) = g(P^{\prime}(x^{\prime}(x, b)) = \frac{x^{\prime}(x, b)}{u^{\prime}(b)} &= \frac{[(v - 1)b + 1]x_t - bx_v}{[(v - 1)b + 1]t - bv}\\
  &= \frac{\gamma [(v - 1)b + 1]t + \lambda b - bx}{[(v - 1)b + 1]t - bv}\\
  &= \gamma - b\frac{x - \lambda - \gamma v}{[(v - 1)b + 1]t - bv}\\
  &= \gamma - b\frac{x - \lambda - \gamma v}{[(v - 1)t - v]b + t}.
\end{split}
\]
The denominator $[(v - 1)t - v]b + t$ is positive; this is immediate when $v \ge 2$ and follows from Lemma \ref{lem:v1} when $v = 1$. Then $g(P(x))$ is a strictly decreasing function of $x$ with limit $L(b)$ as $x$ approaches $x_{b-1}$ from below.
\hspace{\stretch{1}}$\blacksquare$

Then $L(b)$ is a strict lower bound on any 1-problem $P(x)$ with $x_b < x < x_{b-1}$. We are interested in the behavior of $L(b)$ as a function of $b$ and, in particular, how it compares to $\tau$.

\begin{lemma} \label{lem:L(b)} Suppose that $v \ge 2$. (i) If $t = v = u = 2$, then $L(b) = \tau$ for all $b \ge 1$; (ii) If $u = t$ and either $t > 2$ or $v > 2$, then $L(1) = \tau$ and $L(b) > \tau$ for all $b \ge 2$; (iii) If $u < t$, then $L(b) > \tau$ for all $b \ge 1$; (iv) If $u > t$, then $L(1) = \min_{b \ge 1} L(b) < \tau$. In summary, when $v \ge 2$, $L(b) \ge \min\{\tau, L(1)\}$.
\end{lemma}
\textbf{Proof.} When $t = v = 2$ so that $d = 0$, we have
\[
  L(b) = \gamma - [\gamma(u - 2) - \lambda]\frac{b}{u(b - 1) + 2} = \gamma - \frac{\gamma(u - 2) - \lambda}{u - (u - 2)/b}  
\]
and
\[
  \tau = \gamma - \frac{\gamma(u - 2) - \lambda}{u}.
\]
Then if $u = 2$, $L(b) = \tau$ for all values of $b$. This establishes case (i). If $u > 2$, then $L(b)$ is increasing in $b$ and converges to $\tau$ as $b$ tends to $+\infty$, so $L(b) < \tau$ for all values of $b$. This establishes case (iv) when $t = v = 2$.

Otherwise, $t \ge 2$, $v \ge 2$, and either $t > 2$ or $v > 2$. Then $d > 0$. Then
\begin{equation}
\begin{aligned}
  L(b) &\ge \tau\\
  \Leftrightarrow \gamma - [\gamma(u - v) - \lambda]M(b) &\ge \gamma - \frac{\gamma(u - v) - \lambda}{d + u}\\
  \Leftrightarrow \frac{1}{d + u} &\ge M(b)\\
  \Leftrightarrow \frac{1}{d + u} &\ge \frac{b[d(b - 1) + t]}{[(d + u)(b - 1) + t][db + t]}\\
  \Leftrightarrow [(d + u)(b - 1) + t][db + t] &\ge (d + u)b[d(b - 1) + t]\\
  \Leftrightarrow d(d + u)b(b - 1) + (d + u)t(b - 1) + (db + t)t &\ge d(d + u)b(b - 1) + (d + u)tb\\
  \Leftrightarrow (db + t)t - (d + u)t &\ge 0\\
  \Leftrightarrow db &\ge d + u - t\\
  \Leftrightarrow b &\ge 1 + \frac{u - t}{d}.
\end{aligned}
\label{L(b)>tau}
\end{equation}
Then if $u = t$, $L(1) = \tau$ and $L(b) > \tau$ for all $b \ge 2$. This establishes case (ii). If $u < t$, then $L(b) > \tau$ for all values of $b$. This establishes case (iii).

Otherwise $u > t$. Let us investigate how $M(b)$ varies with $b$. Clearly, $M(b)$ is positive and finite for all $b \ge 1$. Write
\[
  M(b) = \frac{\alpha b^2 - \beta b}{\rho b^2 - \delta b - \varepsilon},
\]
where
\[
\begin{split}
  \alpha &= d\\
  \beta &= d - t\\
  \rho &= d(d + u)\\
  \delta &= d(d + u - t) - (d + u)t\\
  \varepsilon &= (d + u - t)t.
\end{split}
\]
Then its derivative is
\[
\begin{split}
  \frac{dM(b)}{db} &= \frac{(2\alpha b - \beta)(\rho b^2 - \delta b - \varepsilon) - (\alpha b^2 - \beta b)(2\rho b - \delta)}{(\rho b^2 - \delta b - \varepsilon)^2}\\
  &= \frac{(\beta \rho - \alpha \delta)b^2 - 2\alpha \varepsilon b + \beta \varepsilon}{(\rho b^2 - \delta b - \varepsilon)^2} \equiv \frac{m(b)t}{(\rho b^2 - \delta b - \varepsilon)^2}.
\end{split}
\]
Then the sign of the derivative is determined by the sign of the numerator $m(b)t$, which is quadratic in $b$. Now
\[
  \beta \rho - \alpha \delta = d(d + u)(d - t) - d^2(d + u - t) + d(d + u)t = d^2t,
\]
so
\begin{align*}
  m(b) &= d^2 b^2 - 2d(d + u - t) b + (d - t)(d + u - t)\\
  &= [db - (d + u - t)]^2 - u(d + u - t).
\end{align*}
This is a quadratic that goes from positive to negative to positive. 

But $m(1) = (u - t)^2 - u(d + u - t) = -t(u - t) - ud < 0$ because we are assuming $u > t$. Then on $b \ge 1$, $m(b)$ goes from negative to positive, so $M(b)$ is initially decreasing and then increasing. Conversely, $L(b)$ is initially increasing and then decreasing on $b \ge 1$. Then the minimum value taken by $L(b)$ on $b \ge 1$ is $L(1) < \tau$. This completes the proof for case (iv).
\hspace{\stretch{1}}$\blacksquare$

\begin{lemma} \label{lem:L(b)_a=0} Let $P(x)$ be a 1-problem with $x_b < x < x_{b - 1}$ for some $b > 0$, and with $t = v = 2$. Then if $u = 2$, $g(P(x)) > \tau$, and if $u > 2$, then $g(P(x)) < \tau$ for $x$ sufficiently close to but less than $x_{b - 1}$.
\end{lemma}

\begin{lemma} \label{lem:L(b)_a>0} Let $P(x)$ be a weakly $u$-weighted 1-problem with $x_b < x < x_{b - 1}$ for some $b > 0$, and with either $t > 2$ or $v > 2$. Then
\begin{enumerate}
\item[(i)] if $b \ge 3$, then $g(P(x)) > \tau$;
\item[(ii)] if $b = 2$, then $g(P(x)) > \tau$ if $u \le v(t - 1)$, else $g(P(x)) < \tau$ for $x$ sufficiently close to but less than $x_1$; \text{ and}
\item[(iii)] if $b = 1$, then $g(P(x)) > \tau$ if $u \le t$, else $g(P(x)) < \tau$ for $x$ sufficiently close to but less than $x_0$.
\end{enumerate}
\end{lemma}
\textbf{Proof.} Recall \eqref{L(b)>tau}:
\[
  L(b) \ge \tau \Leftrightarrow b \ge 1 + \frac{u - t}{d}.
\]
If $b = 1$, this condition is satisfied provided $u \le t$. If $b = 2$, the condition requires
\[
  \frac{u - t}{(v - 1)t - v} \le 1 \Leftrightarrow u \le v(t - 1).
\]
If $b \ge 3$, the condition requires
\[
  \frac{u - t}{(v - 1)t - v} \le 1 \Leftrightarrow u \le 2(v - 1)t - 2v + t = (v - 1)t + v(t - 2),
\]
which is always satisfied because by assumption $u \le (v - 1)t$.

By Lemma \ref{lem:L(b)_LB}, $L(b)$ is the limit of $g(P(x))$ as $x$ approaches $x_{b-1}$ from below, so the conclusions follow.
\hspace{\stretch{1}}$\blacksquare$

\begin{proposition} \label{prop:L(1)_or_tau} Let $P$ be a 1-problem with $v \ge 2$. Then $g(P) > \min\{\tau, L(1)\}$.
\end{proposition}
\textbf{Proof.} This is a direct consequence of Lemmas \ref{lem:L(b)_LB} and \ref{lem:L(b)}.
\hspace{\stretch{1}}$\blacksquare$

Finally, we can now easily prove the desired lower bounds on $N$-problems with $N \ge 2$.

\begin{proposition} \label{lem:Y'_LB_indhyp} Let $P$ be an $N$-problem with $N \ge 2$ and $v \ge 2$ such that if $N = 2$ and $u < v$, then $b > 1$. Assume the inductive hypothesis holds for any $(N - 1)$-problem, and in particular for the reduced problem $P^{\prime}$. Then, in a solution provided by Algorithm \ref{alg:recursive}, 
\[
  g(P) > \frac{(v - 1)x_t + x_u - x_v}{(v - 1)t + u - v} = \tau.
\]
\end{proposition}
\textbf{Proof.} Because $v \ge 2$, the reduced problem $P^{\prime}$ is $u$-weighted by Proposition \ref{prop:reduced_problem_properties}, i.e., $u^{\prime} \ge v^{\prime} \ge 2$.

If $N = 2$, then $P^{\prime}$ is a 1-problem with $v^{\prime} \ge 2$. By Proposition \ref{prop:L(1)_or_tau}, $g(P^{\prime}) > \min\{\tau^{\prime}, L^{\prime}(1)\}$. If $\tau^{\prime} \ge g(P^{\prime}) > L^{\prime}(1)$, then because by assumption either $v \ge 2$ and $b \ge 2$, or $u \ge v \ge 2$ and $b = 1$, Lemma \ref{lem:Y'_2prob} applies and we can conclude that $g(P^{\prime}) > \tau$. Otherwise, $g(P^{\prime}) > \tau^{\prime}$ and Lemma \ref{lem:Y'_Nprob} applies: once again we can conclude that $g(P^{\prime}) > \tau$.

If $N > 2$, then because $u^{\prime} \ge v^{\prime} \ge 2$, by the inductive hypothesis, $g(P^{\prime}) > \tau^{\prime}$. Then by Lemma \ref{lem:Y'_Nprob}, $g(P^{\prime}) > \tau$.

Finally, we appeal to Proposition \ref{prop:g=g'}: $g(P) = g(P^{\prime}) > \tau$.
\hspace{\stretch{1}}$\blacksquare$

\subsection{Proof of the Inductive Hypothesis}

The sequence of results in the prior subsections completes the inductive proof. The following theorem summarizes the results.

\begin{theorem} \label{thm:g} Let $P = (T; U, V)$ be a 3M-DAP. Suppose that $P$ is an $N$-problem, and if $N > 0$, let $P^{\prime}$ be the corresponding reduced problem. Then in a solution provided by Algorithm \ref{alg:recursive}:
\begin{enumerate}
\item if $N > 0$, then $g(P) = g(P^{\prime})$;
\item[2.1] if $N = 0$, then $g(P) = x_u/u$;
\item[2.2] if $N = 1$ and $v > 1$, then $g(P) > \min\left\{\tau, L(1)\right\}$;
\item[2.3] if $N > 1$ and $v > 1$ (except when $N = 2$, $u < v$, and $b = 1$), then $g(P) > \tau$;
\item[3.1] all elements of $U$ have size less than $x_v/v$;
\item[3.2] all elements of $V$ have size greater than $x_u/u$.
\end{enumerate}
\end{theorem}

Finally, we note that Conjecture \ref{NprobLB_conj}, which we relied on to motivate the solution approach, does indeed hold for most 3M-DAP with $(v - 1)t - v > 0$ that are not 0-problems (the exceptions are when $N = 2$, $u < v$, and $b = 1$, cf. Figure \ref{fig:two}).

\begin{corollary} \label{cor:Conj_NprobLB_conj} Let $P = (T; U, V)$ be a 3M-DAP with $(v - 1)t - v > 0$ and suppose that $P$ is an $N$-problem with $N > 0$ such that if $N = 2$ and $u < v$, then $b > 1$. Then
\[
  f(P) \ge g(P) > \frac{(v - 1)x_t - x_v}{(v - 1)t - v}.
\]
\end{corollary}
\textbf{Proof.} $N > 0$ so $x_u > x_{\infty}$, so
\[
  \frac{x_u}{u} > \tau > \frac{(v - 1)x_t - x_v}{(v - 1)t - v}.
\]
Then the conclusion follows from Lemma \ref{lem:Y'_1prob} when $N = 1$ and from Theorem \ref{thm:g} when $N > 1$.
\hspace{\stretch{1}}$\blacksquare$

\section{Algorithm \ref{alg:recursive} produces an optimal solution for any 3M-DAP}

We now prove that Algorithm \ref{alg:recursive} produces an optimal solution for any 3M-DAP $P = (T; U, V)$. Clearly $f(P) \ge g(P)$. The proof that equality pertains will be direct for a subset of problems and by induction on $N$ for the remainder of problems. 

Suppose $P(x)$ is an $N$-problem with $N > 0$ and $x_b < x < x_{b - 1}$ for some $b > 0$. Any solution to $P$ consists of maximal, inseparable $(T, V)$-pairs. If such a pair is a $k$-pair $(A, B)$, then an upper bound on $f(P)$ is given by $z(k)$, the average value of the $U$-elements in $A$:
\[
  z(k) = \frac{[(v - 1)k + 1]x_t - kx_v}{[(v - 1)k + 1]t - kv} = \frac{(v - 1 + 1/k)x_t - x_v}{(v - 1 + 1/k)t - v},
\]
where the second equality applies only when $k > 0$. This is a decreasing function of $k$. Any solution must contain a $k$-pair with $k \ge b$ and also a $k$-pair with $k \le b - 1$. We will demonstrate that the optimal solution to
\begin{itemize}
\item any 1-problem contains a $b$-pair and no pair larger;
\item any 2-problem for which the reduced problem $P^{\prime}$ has $x^{\prime} < x_1^{\prime}$ contains a $(b - 1)$-pair and no pair smaller; and
\item any other problem consists of only $b$- or $(b - 1)$-pairs.
\end{itemize}

\subsection{Lower bound to exclude pairs larger than a $b$-pair}

The following lemma will only be applicable for reduced problems so its assumptions reflect the conclusions of Proposition \ref{prop:reduced_problem_properties}.

\begin{lemma} \label{lem:g>2xu-xv_N=1} Let $P(x) = (T; U, V)$ be a 3M-DAP with $x_b < x < x_{b-1}$ for some $b > 0$ that is either $u$-weighted or is $v$-weighted with $v = u + 1$. Suppose that either: (i) $P$ is an $N$-problem with $N \ge 2$; or (ii) $P$ is a weakly $u$-weighted 1-problem and $b > 1$. Then, in a solution provided by Algorithm \ref{alg:recursive},
\[
  g(P) > \frac{2x_u - x_v}{2u - v}.
\]
\end{lemma}
\textbf{Proof.} Observe that
\[
  \frac{2x_u - x_v}{2u - v} = \frac{\lambda + x}{2u - v}.
\]
This is an increasing function of $x$. 

\noindent \textbf{Case 1.} $g(P) \ge \tau$. By Theorem \ref{thm:g}, this is certainly the case when $N \ge 2$ (except when $N = 2$, $v = u + 1$, and $b = 1$). Recall that
\[
  \tau = \frac{\lambda + (v - 1)t\gamma}{u - v + (v - 1)t}.
\]

\textbf{Case 1.1} $u \le (v - 1)t$. Then
\[
  \tau = \frac{\lambda + \gamma(v - 1)t}{(v - 1)t + u - v} \ge \frac{\lambda + \gamma u}{2u - v} > \frac{\lambda + x}{2u - v} = \frac{2x_u - x_v}{2u - v},
\]
because $\lambda < \gamma(u - v)$ and $x < \gamma u$.

\textbf{Case 1.2} $u > (v - 1)t$. An upper bound on $x$ can be obtained on noting that when $b = 1$, so that $x_1 < x < x_0$, $P$ is not a 1-problem. Then $x$ must be sufficiently small to ensure that $x_{\infty}^{\prime} < x^{\prime}$. When $b = 1$:
\[
\begin{split}
  x^{\prime} = x_u^{\prime} &= vx_t - x_v = \lambda + vt\gamma - x\\
  u^{\prime} &= v(t - 1)\\
  x_v^{\prime} &= x_t\\
  v^{\prime} &= t\\
  x_t^{\prime} &= x_u = x\\
  t^{\prime} &= u,
\end{split}
\]
and
\[
  x_{\infty}^{\prime} = \frac{[(v^{\prime} - 1)x_t^{\prime} + x_u^{\prime} - x_v^{\prime}]u^{\prime}}{(v^{\prime} - 1)t^{\prime} + u^{\prime} - v^{\prime}} = (t - 1)v\frac{(t - 2)x + \lambda + (v - 1)t \gamma}{(t - 1)(u + v) - t}.
\]
Then the upper bound on $x$ results:
\[
\begin{split}
  x_{\infty}^{\prime} &< x^{\prime}\\
  \Leftrightarrow (t - 1)v[(t - 2)x + \lambda + (v - 1)t\gamma] &< [(t - 1)(u + v) - t](\lambda + vt\gamma - x)\\
  \Leftrightarrow [(t - 1)u + (t - 1)^2v - t]x &< [(t - 1)u - t]\lambda + [(t - 1)u - 1]vt\gamma\\
  \Leftrightarrow \frac{x}{u} &< \frac{[(t - 1)u - t]\lambda + [(t - 1)u - 1]vt\gamma}{[(t - 1)u - t](u - v) + [(t - 1)u - 1]vt}.
\end{split}
\]
Now
\[
  \{[(u - v)(v - 2) - v](t - 1) + v\}[u - (v - 1)t] + (2u - v)(v - 1)t[(v - 1)t - v] \ge 0,
\]
because if $v = 2$, the expression reduces to $2(t - 2)[(2t - 1)u - t] \ge 0$, and if $v > 2$, then $v - 2 \ge 1$ so $(u - v)(v - 2) - v \ge u$, whence all terms are positive. Then
\[
\begin{split}
  \{[(u - v)(v - 2) - v](t - 1) + v\}[u - (v - 1)t] \quad \quad \quad \quad \quad \quad &\\
  + (2u - v)(v - 1)t[(v - 1)t - v] &\ge 0\\
  \Leftrightarrow (2u - v)(v - 1)(t - 1)[u - (v - 1)t] \quad \quad \quad \quad \quad \quad &\\
  + (2u - v)(v - 1)t[(v - 1)t - v] &\ge [(t - 1)u - 1]v[u - (v - 1)t]\\
  \Leftrightarrow (2u - v)(v - 1)[(t - 1)u - t] &\ge [(t - 1)u - 1]v[u - (v - 1)t]\\
  \Leftrightarrow \frac{(2u - v)(v - 1)}{u - (v - 1)t} &\ge \frac{[(t - 1)u - 1]v}{(t - 1)u - t}\\
  \Leftrightarrow \frac{(2u - v)(v - 1)t}{[u - (v - 1)t](u - v)} &\ge \frac{[(t - 1)u - 1]vt}{[(t - 1)u - t](u - v)}.
\end{split}
\]
It follows that
\[
\begin{split}
  \frac{x}{u} < \frac{[(t - 1)u - t]\lambda + [(t - 1)u - 1]vt\gamma}{[(t - 1)u - t](u - v) + [(t - 1)u - 1]vt} &\le \frac{[u - (v - 1)t]\lambda + (2u - v)(v - 1)t\gamma}{[u - (v - 1)t](u - v) + (2u - v)(v - 1)t}\\
  &= \frac{[u - (v - 1)t]\lambda + (2u - v)(v - 1)t\gamma}{u[(u - v) + (v - 1)t]},
\end{split}
\]
from which it follows that
\[
\begin{split}
  [u - v + (v - 1)t]x &< [u - (v - 1)t]\lambda + (2u - v)(v - 1)t\gamma\\
  \Leftrightarrow (v - 1)t\lambda + [u - v + (v - 1)t]x &< u\lambda + (2u - v)(v - 1)t\gamma\\
  \Leftrightarrow [u - v + (v - 1)t](\lambda + x) &< (2u - v)[\lambda + (v - 1)t\gamma]\\
  \Leftrightarrow \frac{\lambda + x}{2u - v} &< \frac{\lambda + (v - 1)t\gamma}{u - v + (v - 1)t}\\
  \Leftrightarrow \frac{2x_u - x_v}{2u - v} &< \tau \le g(P).
\end{split}
\]

\noindent \textbf{Case 2.} $g(P) < \tau$. Then we can restrict attention to two subcases. 

\textbf{Case 2.1} $P$ a weakly $u$-weighted 1-problem with $b > 1$. We cannot have $t = v = 2$: the problem is weakly $u$-weighted so this would require $u = 2$, but then by Lemma \ref{lem:L(b)_a=0}, $g(P) > \tau$. Then either $t > 2$ or $v > 2$. By Lemma \ref{lem:L(b)_a>0}, we must have $b = 2$ and $u > v(t - 1)$. Then $L(2) < g(P) < \tau$.

Then $v(t - 1) < u \le (v - 1)t$, from which it follows that $t < v$. Recall that
\[
  x_1 = \frac{\lambda + v\gamma t}{u - v + vt}u,
\]
so
\[
  \lambda = (u - v + vt)\frac{x_1}{u} - v\gamma t.
\]
Then
\[
  L(2) = \frac{2\lambda + \gamma(2v - 1)t - 2x_1}{2(u - v) + (2v - 1)t - 2u} = \frac{2v(t - 1)x_1/u - \gamma t}{2v(t - 1) - t} = \frac{2v^2(t - 1)x_1/u - v\gamma t}{2v^2(t - 1) - vt},
\]
and
\[
\begin{split}
  \frac{\lambda + x_1}{2u - v}
  &= \frac{(2u - v + vt)x_1/u - v\gamma t}{(2u - v + vt) - vt}.
\end{split}
\]
Now $2 \le t < v$, so $2v - 1 > 2(v - 1)$ and $v(t - 1) \ge v > t$, so
\[
  v(t - 1)(2v - 1) > 2(v - 1)t \ge 2u \quad \quad \Leftrightarrow \quad \quad 2v^2(t - 1) > 2u + v(t - 1) = 2u - v + vt.
\]
Then
\[
  g(P) > L(2) = \frac{2v^2(t - 1)x_1/u - v\gamma t}{2v^2(t - 1) - vt} > \frac{(2u - v + vt)x_1/u - v\gamma t}{(2u - v + vt) - vt} = \frac{\lambda + x_1}{2u - v} > \frac{2x_u - x_v}{2u - v}.
\]
%Then we certainly have $x < x_1$ and
%\[
%  g(P) > L(2) = \gamma - 2 \frac{x_1 - \lambda - \gamma v}{2v(t - 1) - t} = \frac{2\lambda + \gamma(2v - 1)t - 2x_1}{2v(t - 1) - t}.
%\]
%Using the previous lemma, we have
%\[
%\begin{split}
%  (2v + 1) x_1 &\ge \lambda(2v - 1) + \gamma [2u + (2v - 1)v]\\
%  \Leftrightarrow 4\lambda u + 4\gamma vtu - \lambda(2v - 1)t - \gamma [2u + (2v - 1)v]t &\ge 4(vt + u - v)x_1 - (2v + 1)t x_1\\
%  \Leftrightarrow \lambda(4u - 2vt + t) + \gamma(2v - 1)t(2u - v) &\ge (4u - 4v + 2vt - t)x_1\\
%  \Leftrightarrow [2\lambda + \gamma(2v - 1)t - 2x_1](2u - v) &\ge (\lambda + x_1)[2v(t - 1) - t]\\
%  \Leftrightarrow L(2) = \frac{2\lambda + \gamma(2v - 1)t - 2x_1}{2v(t - 1) - t} &\ge \frac{\lambda + x_1}{2u - v}\\
%  \Leftrightarrow g(P) &> \frac{\lambda + x}{2u - v} = \frac{2x_u - x_v}{2u - v}.
%\end{split}
%\]

\textbf{Case 2.2} $P$ is a 2-problem with $v = u + 1$ and $b = 1$. Then
\begin{align*}
  x_u^{\prime \prime} &= [(v^{\prime} - 1)b^{\prime} + 1]x_t^{\prime} - b^{\prime}x_v^{\prime}\\ 
  u^{\prime \prime} &= [(v^{\prime} - 1)b^{\prime} + 1]t^{\prime} - b^{\prime}v^{\prime}\\ 
  x_t^{\prime} &= x_u = x\\
  t^{\prime} &= u\\
  x_v^{\prime} &= x_t\\
  v^{\prime} &= t
\end{align*}
so
\[
  g(P) = \frac{x_u^{\prime \prime}}{u^{\prime \prime}} = \frac{[(t - 1)b^{\prime} + 1]x - b^{\prime}t\gamma}{[(t - 1)b^{\prime} + 1]u - b^{\prime}t}.
\]

If $b^{\prime} = 1$, then
\[
  g(P) = \frac{x - \gamma}{u - 1} > \frac{\lambda + x}{u - 1} = \frac{2x_u - x_v}{2u - v},
\]
because $-\gamma = \gamma (u - v) > \lambda$.

If $b^{\prime} > 1$ and either $t > 2$ or $u > 2$, then
\begin{align*}
  \frac{2u}{u + 1}b^{\prime} &< (t - 1)(u - 1)b^{\prime} + u\\
  \Leftrightarrow \frac{u - 1}{u + 1}b^{\prime} &< [(t - 1)b^{\prime} + 1]u - b^{\prime}t\\
  \Leftrightarrow \frac{b^{\prime}(u - 1)t}{[(t - 1)b^{\prime} + 1]u - b^{\prime}t} &< (u + 1)t\\
  \Leftrightarrow \frac{\{[(t - 1)b^{\prime} + 1]u - b^{\prime}t\}\lambda + b^{\prime}(u - 1)t\gamma}{b^{\prime}t - [(t - 1)b^{\prime} + 1]} &< \frac{\lambda + (u + 1)t\gamma}{(u + 1)t - 1}u = x_1 < x.
\end{align*}
Then
\begin{align*}
  \{b^{\prime}t - [(t - 1)b^{\prime} + 1]\}x &> \{[(t - 1)b^{\prime} + 1]u - b^{\prime}t\}\lambda + b^{\prime}(u - 1)t\gamma\\
  \Leftrightarrow [(t - 1)b^{\prime} + 1]x(u - 1) - b^{\prime}(u - 1)t\gamma &> \{[(t - 1)b^{\prime} + 1]u - b^{\prime}t\}x + \{[(t - 1)b^{\prime} + 1]u - b^{\prime}t\}\lambda\\
  \Leftrightarrow \{[(t - 1)b^{\prime} + 1]x - b^{\prime}t\gamma\}(u - 1) &> \{[(t - 1)b^{\prime} + 1]u - b^{\prime}t\}(x + \lambda)\\
  \Leftrightarrow g(P) = \frac{[(t - 1)b^{\prime} + 1]x - b^{\prime}t\gamma}{[(t - 1)b^{\prime} + 1]u - b^{\prime}t} &> \frac{x + \lambda}{u - 1}.
\end{align*}

If $b^{\prime} > 1$ and $t = u = 2$, then we need to argue more carefully. We have $x^{\prime} < x_{b^{\prime} - 1}^{\prime}$, so
\begin{align*}
  x^{\prime} &< x_{b^{\prime} - 1}^{\prime}\\
  \Leftrightarrow vx_t - x_v &< \frac{[\lambda^{\prime} + \gamma^{\prime}(v^{\prime} - 1)t^{\prime}](b^{\prime} - 1) + \gamma^{\prime}t^{\prime}}{[(v^{\prime} - 1)t^{\prime} + u^{\prime} - v^{\prime}](b^{\prime} - 1) + t^{\prime}}u^{\prime}\\
  \Leftrightarrow vx_t - x_v &< \frac{[(v - 1)x_t - x_v + (t - 1)x](b^{\prime} - 1) + x}{[(t - 1)u + (t - 1)v - t](b^{\prime} - 1) + u}(t - 1)v\\
  \Leftrightarrow \lambda + 6\gamma - x &< \frac{3(\lambda + 4\gamma)(b^{\prime} - 1) + x}{3b^{\prime} - 1}\\
  \Leftrightarrow (3b^{\prime} - 1)\lambda + 6(3b^{\prime} - 1)\gamma - (3b^{\prime} - 1)x &< 3(b^{\prime} - 1)\lambda + 12(b^{\prime} - 1)\gamma + 3x\\  
  \Leftrightarrow 2\lambda + 6(b^{\prime} + 1)\gamma &< (3b^{\prime} + 2)x.
\end{align*}
Now $\lambda + \gamma < 0$, so
\begin{align*}
  0 &> 2(2b^{\prime} + 3)(\lambda + \gamma)\\
  \Leftrightarrow (b^{\prime} - 1)[2\lambda + 6(b^{\prime} + 1)\gamma] &> (3b^{\prime} + 2)[2\lambda + 2b^{\prime}\gamma]\\
  \Leftrightarrow \frac{2\lambda + 6(b^{\prime} + 1)\gamma}{3b^{\prime} + 2} &> \frac{2\lambda + 2b^{\prime}\gamma}{b^{\prime} - 1}.
\end{align*}
Then
\begin{align*}
  x&> \frac{2\lambda + 6(b^{\prime} + 1)\gamma}{3b^{\prime} + 2}\\
  \Leftrightarrow x&> \frac{2\lambda + 2b^{\prime}\gamma}{b^{\prime} - 1}\\
  \Leftrightarrow (b^{\prime} - 1)x &> 2\lambda + 2b^{\prime}\gamma\\
  \Leftrightarrow (b^{\prime} + 1)x - 2b^{\prime}\gamma &> 2\lambda + 2x\\
  \Leftrightarrow g(P) = \frac{1}{2}[(b^{\prime} + 1)x - 2b^{\prime}\gamma] &> \lambda + x.
\end{align*}

\hspace{\stretch{1}}$\blacksquare$

\begin{proposition} \label{prop:g>z(b+1)} Let $P(x) = (T; U, V)$ be a 3M-DAP with $x_b < x < x_{b-1}$ for some $b > 0$. Suppose that $P$ is one of: (i) an $N$-problem with $N \ge 3$; (ii) a 2-problem with $b > 1$ and the reduced problem $P^{\prime}$ has $x^{\prime} < x_1^{\prime}$; or (iii) a 1-problem. Then in a solution provided by Algorithm \ref{alg:recursive}, $g(P) > z(k)$ for all $k > b$.
\end{proposition}
\textbf{Proof.} Because $z(k)$ is a decreasing function of $k$, it suffices to show the result for $k = b + 1$. Observe that
\[
  z(b + 1) = \frac{[(v - 1)(b + 1) + 1]x_t - (b + 1)x_v}{[(v - 1)(b + 1) + 1]t - (b + 1)v}.
\]
Because $P$ is not a 0-problem, we can form the reduced problem $P^{\prime}$. Now
\[
  2u^{\prime} - v^{\prime} = 2[(v - 1)b + 1]t - 2bv - [(v - 1)(b - 1) + 1]t + (b - 1)v = [(v - 1)(b + 1) + 1]t - (b + 1)v,
\]
and
\[
  2x_u^{\prime} - x_v^{\prime} = 2[(v - 1)b + 1]x_t - 2bx_v - [(v - 1)(b - 1) + 1]x_t + (b - 1)x_v = [(v - 1)(b + 1) + 1]x_t - (b + 1)x_v.
\]
Then
\[
  z(b + 1) = \frac{2x_u^{\prime} - x_v^{\prime}}{2u^{\prime} - v^{\prime}}.
\]
Because we have established that $g(P) = g(P^{\prime})$, it suffices to show that
\[
  g(P^{\prime}) > \frac{2x_u^{\prime} - x_v^{\prime}}{2u^{\prime} - v^{\prime}}.
\]
This follows from the previous lemma for the first two cases, on noting that in the second case, because $b > 1$, by Proposition \ref{prop:reduced_problem_properties} $P^{\prime}$ is weakly $u$-weighted, and because $x^{\prime} < x_1^{\prime}$, $P^{\prime}$ has $b^{\prime} > 1$. 

For the third case, $P^{\prime}$ is a 0-problem. Then
\[
  g(P^{\prime}) = \frac{x_u^{\prime}}{u^{\prime}} > \frac{2x_u^{\prime} - x_v^{\prime}}{2u^{\prime} - v^{\prime}}.
\]
\hspace{\stretch{1}}$\blacksquare$

\subsection{Lower bound to exclude pairs smaller than a $(b - 1)$-pair}

Again, the following lemma will only be applicable for reduced problems so its assumptions reflect the conclusions of Proposition \ref{prop:reduced_problem_properties}.

\begin{lemma} \label{lem:g>(xt - (2xv-xu)/(2v - u))/(t - 1)_N=2} Let $P(x) = (T; U, V)$ be a 3M-DAP with $x_b < x < x_{b-1}$ for some $b > 1$ and $2v > u$ that is either $u$-weighted or is $v$-weighted with $v = u + 1$. Suppose that either: (i) $P$ is an $N$-problem with $N \ge 2$; or (ii) $P$ is a weakly $u$-weighted 1-problem. Then, in a solution provided by Algorithm \ref{alg:recursive},
\[
  g(P) > \frac{1}{t - 1}\left(x_t - \frac{2x_v - x_u}{2v - u}\right).
\]
\end{lemma}
\textbf{Proof.} Observe that
\[
  \frac{1}{t - 1}\left(x_t - \frac{2x_v - x_u}{2v - u}\right) = \frac{1}{t - 1}\left(\gamma t - \frac{x - 2\lambda}{2v - u}\right) = \frac{\gamma (2v - u)t + 2\lambda - x}{(t - 1)(2v - u)}.
\]
This is a decreasing function of $x$. 

\noindent \textbf{Case 1.} $g(P) \ge \tau$. By Theorem \ref{thm:g}, this is certainly the case when $N \ge 2$ (because $b > 1$).

Because $P$ is not a 0-problem, we must have $x > x_{\infty}$. Then it suffices to show that
\[
  g(P) \ge \frac{\gamma (2v - u)t + 2\lambda - x_{\infty}}{(t - 1)(2v - u)}.
\]
Now $u \ge 2 \Leftrightarrow u/2 \ge 1 \Leftrightarrow v - 1 \ge v - u/2$. Then, because $\lambda < \gamma(u - v)$,
\[
  \frac{x_{\infty}}{u} = \frac{\lambda + \gamma(v - 1)t}{(v - 1)t + u - v} \ge \frac{\lambda + \gamma (v - u/2)t}{(v - u/2)t + u - v} = \frac{2\lambda + \gamma (2v - u)t}{(2v - u)t + 2u - 2v},
\]
with equality if and only if $u = 2$. Taking the mediant of the right-hand side with $x_{\infty}/u$ gives
\begin{equation} \label{ub}
  \frac{x_{\infty}}{u} \ge \frac{2\lambda + \gamma (2v - u)t}{(2v - u)t + 2u - 2v} \ge \frac{\gamma (2v - u)t + 2\lambda - x_{\infty}}{(t - 1)(2v - u)},
\end{equation}
with the two inequalities being equalities if and only if $u = 2$. Then if $g(P) \ge \tau$, we have
\[
  g(P) \ge \tau = \frac{x_{\infty}}{u} \ge \frac{\gamma (2v - u)t + 2\lambda - x_{\infty}}{(t - 1)(2v - u)}.
\]

\noindent \textbf{Case 2.} $g(P) < \tau$. Then we can restrict attention to $P$ a weakly $u$-weighted 1-problem with $b > 1$. We cannot have $t = v = 2$: the problem is weakly $u$-weighted so this would require $u = 2$, but then by Lemma \ref{lem:L(b)_a=0}, $g(P) > \tau$. Then either $t > 2$ or $v > 2$. By Lemma \ref{lem:L(b)_a>0}, we must have $b = 2$ and $u > v(t - 1)$. But by assumption $u < 2v$, so $t = 2$.

%Because $b > 1$ by assumption, certainly $x < x_1$ and
%\[
%  g(P) > L(2) = \gamma - 2 \frac{x_1 - \lambda - \gamma v}{2v(t - 1) - t} = \frac{2\lambda + \gamma(2v - 1)t - 2x_1}{(2v - 1)t - 2v},
%\]
%by Lemma \ref{lem:L(b)_a=0} or Lemma \ref{lem:L(b)_a>0}. Recall that
%\[
%  x_1 = \frac{\lambda + \gamma vt}{vt + u - v}u \Leftrightarrow (vt + u - v)x_1 = \lambda u + \gamma vtu.
%\]
%Then one of two cases pertains.
%
%\noindent \textbf{Case 1.} $t = v = 2$. Then $u = 3$ because $P$ is a 3M-DAP, and this does satisfy the assumption that $u < 2v$. Then $x_{\infty} = \lambda + 2\gamma$ and
%\[
%  x_1 = \frac{3\lambda + 12\gamma}{5}.
%\]
%Because $u - v = 1$, $\gamma = (u - v)\gamma > \lambda$. Then
%\[
%\begin{split}
%  g(P) > L(2) &= \frac{2\lambda + \gamma(2v - 1)t - 2x_1}{(2v - 1)t - 2v}\\
%  &= \lambda + 3\gamma - x_1\\
%  &= \frac{1}{5}(2\lambda + 3\gamma)\\
%  &> \lambda\\
%  &= 2\lambda + 2\gamma - x_{\infty}\\
%  &= \frac{\gamma (2v - u)t + 2\lambda - x_{\infty}}{(t - 1)(2v - u)}.
%\end{split}
%\]
%
%\noindent \textbf{Case 2.} Either $t > 2$ or $v > 2$ and $u > v(t - 1)$. 

Because $P$ is a weakly $u$-weighted 3M-DAP, $u \le (v - 1)t$. Then $2 = t < v < u \le 2v - 2$. Because $b = 2$, $x_2 < x < x_1$ and
\[
  g(P) = \frac{(2v - 1)x_t - 2x_v}{(2v - 1)t - 2v} = \frac{2\lambda + \gamma (2v - 1)t - 2x}{(2v - 1)t - 2v} = \frac{\lambda + \gamma (2v - 1) - x}{v - 1}.
\]
We want to contrast this with
\[
  \frac{1}{t - 1}\left(x_t - \frac{2x_v - x_u}{2v - u}\right) = \frac{2\lambda + 2(2v - u)\gamma - x}{2v - u}
\]
Now
\[
  \frac{x_{\infty}}{u} = \frac{\lambda + 2(v - 1)\gamma}{u + v - 2},
\]
and taking the mediant of this with $\gamma > x_{\infty}/u$ gives
\[
  \frac{\lambda + \gamma (2v - 1)}{u + v - 1} > \frac{x_{\infty}}{u} > \frac{2\lambda + 2(2v - u)\gamma}{2v},
\]
where the second inequality is taken from \eqref{ub}. Taking mediants on both sides with $x/u$ (with a coefficient of $-1$) and noting that $2v \le u + v - 1$ gives
\[
  g(P) = \frac{\lambda + \gamma (2v - 1) - x}{v - 1} > \frac{2\lambda + 2(2v - u)\gamma - x}{2v - u} = \frac{1}{t - 1}\left(x_t - \frac{2x_v - x_u}{2v - u}\right),
\]
as desired.

%This is a mediant of $x/u$ (with a coefficient of $-1$) and
%\[
%  \frac{\lambda + \gamma (2v - 1)}{u + v - 1} > \frac{x_{\infty}}{u},
%\]
%where the inequality follows because the left-hand side is itself a mediant of $\gamma$ and
%and $\gamma > x_{\infty}/u$. On the other hand,
%\[
%  \frac{1}{t - 1}\left(x_t - \frac{2x_v - x_u}{2v - u}\right) = \frac{2\lambda + 2(2v - u)\gamma - x}{2v - u}
%\]
%is a mediant of $x/u$ (again with a coefficient of $-1$) and
%\[
%  \frac{2\lambda + 2(2v - u)\gamma}{2v} < \frac{x_{\infty}}{u},
%\]
%where the inequality was shown above. Because $2v \le u + v - 1$, the conclusion follows.
\hspace{\stretch{1}}$\blacksquare$

\begin{proposition} \label{prop:g>(x_u - z(b-2))/(u - 1)} Let $P(x) = (T; U, V)$ be a 3M-DAP with $x_b < x < x_{b-1}$ for some $b > 0$. Suppose that either: (i) $P$ is an $N$-problem with $N \ge 3$; or (ii) $P$ is a 2-problem with $b > 1$ and the reduced problem $P^{\prime}$ has $x^{\prime} < x_1^{\prime}$. Then in a solution provided by Algorithm \ref{alg:recursive}, $g(P) > (x_u - z(k))/(u - 1)$ for all $k < b - 1$.
\end{proposition}
\textbf{Proof.} Because $z(k)$ is a decreasing function of $k$, it suffices to show the result for $k = b - 2$. Observe that
\[
  z(b - 2) = \frac{[v - 1 + 1/(b - 2)]x_t - x_v}{[v - 1 + 1/(b - 2)]t - v} = \frac{[(v - 1)(b - 2) + 1]x_t - (b - 2)x_v}{[(v - 1)(b - 2) + 1]t - (b - 2)v}.
\]
Because $P$ is not a 0-problem, we can form the reduced problem $P^{\prime}$. Now
\[
\begin{split}
  2v^{\prime} - u^{\prime} &= 2[(v - 1)(b - 1) + 1]t - 2(b - 1)v - [(v - 1)b + 1]t + bv\\
  &= [(v - 1)(b - 2) + 1]t - (b - 2)v > 0,
\end{split}
\]
and
\[
\begin{split}
  2x_v^{\prime} - x_u^{\prime} &= 2[(v - 1)(b - 1) + 1]x_t - 2(b - 1)x_v - [(v - 1)b + 1]x_t + bx_v\\
  &= [(v - 1)(b - 2) + 1]x_t - (b - 2)x_v.
\end{split}
\]
Then
\[
  z(b - 2) = \frac{2x_v^{\prime} - x_u^{\prime}}{2v^{\prime} - u^{\prime}}.
\]

In the first case $P^{\prime}$ is an $N$-problem with $N \ge 2$. In the second case, because $b > 1$, by Proposition \ref{prop:reduced_problem_properties} $P^{\prime}$ is weakly $u$-weighted, and because $x^{\prime} < x_1^{\prime}$, $P^{\prime}$ has $b^{\prime} > 1$. In both cases, $2v^{\prime} > u^{\prime}$. So Lemma \ref{lem:g>(xt - (2xv-xu)/(2v - u))/(t - 1)_N=2} applies. Then
\[
  g(P) = g(P^{\prime}) > \frac{1}{t^{\prime} - 1}\left(x_t^{\prime} - \frac{2x_v^{\prime} - x_u^{\prime}}{2v^{\prime} - u^{\prime}}\right) = \frac{1}{u - 1}\left(x_u - \frac{2x_v^{\prime} - x_u^{\prime}}{2v^{\prime} - u^{\prime}}\right) = \frac{x_u - z(b - 2)}{u - 1}.
\]
\hspace{\stretch{1}}$\blacksquare$

\subsection{Proof of Optimality}

For any problem that is not a 0-problem, any solution consists of maximal, inseparable pairs. Suppose a solution includes a maximal inseparable $k$-pair. The average size of a $U$-element in this pair is $z(k)$, so there is a $U$-element with at most this size, so $z(k)$ is an upper bound on the solution. Further, there is a $U$-element with at least size $z(k)$; considering a row of $U$ containing such an element leads us to conclude that there is a $U$-element with size at most
\[
  \frac{x_u - z(k)}{u - 1},
\]
and this provides an alternative upper bound on the solution. We use these facts to establish that Algorithm \ref{alg:recursive} always produces an optimal solution. First, we consider two cases that require special consideration.

\begin{lemma} \label{lem:f=g_N=2_b'=1} Consider a 3M-DAP $P(x) = (T; U, V)$ that is a 2-problem with $x_b < x < x_{b-1}$ for some $b \ge 1$, for which the reduced problem $P^{\prime}(x^{\prime})$ has $x^{\prime}_1 < x^{\prime} < x^{\prime}_0$. Then $f(P) = g(P)$, i.e., Algorithm \ref{alg:recursive} produces an optimal solution.
\end{lemma}
\textbf{Proof.} We have $x_u^{\prime \prime} = v^{\prime}x_t^{\prime} - x_v^{\prime}$, $u^{\prime \prime} = v^{\prime}t^{\prime} - v^{\prime}$, and $g(P^{\prime}) = x_u^{\prime \prime}/u^{\prime \prime}$. Then
\[
\begin{split}
  g(P) = g(P^{\prime}) = \frac{v^{\prime}x_t^{\prime} - x_v^{\prime}}{v^{\prime}t^{\prime} - v^{\prime}} = \frac{v^{\prime}x_u - x_v^{\prime}}{v^{\prime}u - v^{\prime}} = \frac{x_u - x_v^{\prime}/v^{\prime}}{u - 1} = \frac{x_u - z(b - 1)}{u - 1}.
\end{split}
\]
Because any solution must contain a maximal inseparable $k$-pair with $k \le b - 1$ and $z(k)$ is decreasing in $k$, it follows that $[x_u - z(b - 1)]/(u - 1)$ is an upper bound on $f(P)$. Then the conclusion follows.
\hspace{\stretch{1}}$\blacksquare$

\begin{lemma} \label{lem:f=g_N=2_b=1} Consider a 3M-DAP $P(x) = (T; U, V)$ that is a 2-problem with $x_1 < x < x_0$. Then $f(P) = g(P)$, i.e., Algorithm \ref{alg:recursive} produces an optimal solution.
\end{lemma}
\textbf{Proof.} Suppose that the reduced problem $P^{\prime}(x^{\prime})$ has $x_{b^{\prime}}^{\prime} < x^{\prime} < x_{b^{\prime} - 1}^{\prime}$ for some $b^{\prime} \ge 1$. Then
\[
  g(P) = g(P^{\prime}) = z^{\prime}(b^{\prime}) = \frac{[(v^{\prime} - 1)b^{\prime} + 1]x_t^{\prime} - b^{\prime}x_v^{\prime}}{[(v^{\prime} - 1)b^{\prime} + 1]t^{\prime} - b^{\prime}v^{\prime}}.
\]
%and note that this is a function only of the dimensions and rowsums of the matrices $T^{\prime}$ and $V^{\prime}$. 
To solve $P$ Algorithm \ref{alg:recursive} employs only 1- and 0-pairs.
%: the reduced problem has:
%\begin{align*}
%  s_t^{\prime} &= s_u\\
%  t^{\prime} &= u \ge 2\\
%  x_t^{\prime} &= x_u\\
%  s_u^{\prime} &= s_v\\
%  u^{\prime} &= v(t - 1)\\
%  x^{\prime} = x_u^{\prime} &= vx_t - x_v\\
%  s_v^{\prime} &= s_t - vs_v = s_t - n_v\\
%  v^{\prime} &= t\\
%  x_v^{\prime} &= x_t.
%\end{align*}
%Suppose that
%\begin{align*}
%  v &> \frac{t(2t - 1)(u - 1)}{2(t - 1)}\\
%  \Leftrightarrow 2v(t - 1) &> t(2t - 1)(u - 1)\\
%  \Leftrightarrow 2u^{\prime} &> v^{\prime}(2v^{\prime} - 1)(t^{\prime} - 1),
%\end{align*}
%and note that this assumption requires that $v > u$ (the right-hand side of the first line is strictly greater than $u$ because $t \ge 2$ and $u \ge 2$). Then by the proof of Lemma \ref{lem:g>2xu-xv_N=1}, for $x$ sufficiently small that $x^{\prime}$ is sufficiently close to $x_{b^{\prime} - 1}^{\prime}$, 
%\[
%  g(P) < \frac{2x_u^{\prime} - x_v^{\prime}}{2u^{\prime} - v^{\prime}} = z(2).
%\]
%Then it seems plausible that a solution to $P$ could be found that includes 2-pairs and whose value is greater than $g(P)$.
%
%Any solution to $P$ is formed from $k$-pairs. Consider such a solution consisting of $(T, V)$-pairs $\{(T_1, V_1), \ldots, (T_m, V_m)\}$, where each pair $(T_i, V_i)$ is a $k_i$-pair with $k_i \ge 1$, together with a group of 0-pairs that we represent by the matrix pair $(T_0, V_0)$. This leads to a reduced problem $\hat{P} = (\hat{T}; \{\boldsymbol{\hat{u}}_1, \ldots, \boldsymbol{\hat{u}}_m\}, \hat{V})$ that is a DAP. Here $\hat{T} = U$ and $\hat{v} = t$. For this reduced problem
By Lemma \ref{lem:maximal_inseparable_pairs} any solution to $P$ is formed from maximal, inseparable $(T, V)$-pairs; at least $s_t - n_v$ of these must be 0-pairs. Consider such a solution consisting of $(T, V)$-pairs $\{(T_1, V_1), \ldots, (T_m, V_m)\}$, where $V_i$ has $k_i \ge 1$ rows, together with a group of 0-pairs that we represent by the matrix pair $(T_0, V_0)$; $V_0$ has no rows. This leads to a reduced problem $\hat{P} = (\hat{T}; \{\boldsymbol{\hat{u}}_1, \ldots, \boldsymbol{\hat{u}}_m\}, \hat{V})$ that is a DAP. Here $\hat{T} = U$ and $\hat{V} = T_0$, so that $\hat{v} = t$. For this reduced problem
\[
  \frac{1}{2}\hat{r}_v(\hat{s}_v) = \hat{s}_t - (\hat{v} - 1)\hat{s}_v = s_u - (t - 1)\hat{s}_v.
\]
Then
\[
  \hat{b}^*(\hat{s}_v) = \frac{2\hat{s}_v}{\hat{r}_v(\hat{s}_v)} = \frac{\hat{s}_v}{s_u - (t - 1)\hat{s}_v} = \frac{1}{s_u/\hat{s}_v - (t - 1)}.
\]
Here we write $\hat{b}^*(\hat{s}_v)$ as a function of $\hat{s}_v$, the number of 0-pairs in the solution to $P$, and note that it is an increasing function thereof. The smallest that $\hat{s}_v$ can be is $s_t - n_v$ when all other pairs are 1-pairs. But by assumption $\hat{b}^*(s_t - n_v) > b^{\prime} - 1$, so in any solution to $P$, the reduced problem $\hat{P}$ has $\hat{b}^* > b^{\prime} - 1$. This means that any solution to $\hat{P}$ must contain a maximal, inseparable $(\hat{T}, \hat{V})$-pair with $\hat{V}$ having $k \ge b^{\prime}$ rows. An upper bound on the value of $\hat{P}$ is given by the average value of the elements in $\hat{T}$ that are from the $\{\boldsymbol{\hat{u}}_1, \ldots, \boldsymbol{\hat{u}}_m\}$. If $\hat{T}$ has $a$ rows, then this average value is
\[
  \frac{ax_t^{\prime} - kx_v^{\prime}}{at^{\prime} - kv^{\prime}} \le \frac{[(v^{\prime} - 1)k + 1]x_t^{\prime} - kx_v^{\prime}}{[(v^{\prime} - 1)k + 1]t^{\prime} - kv^{\prime}} = z^{\prime}(k) \le z^{\prime}(b^{\prime}),
\]
because $a \le (v^{\prime} - 1)k + 1$ by Lemma \ref{lem:inseparable_pair}. Hence $z^{\prime}(b^{\prime})$ is an upper bound on value of both $\hat{P}$ and the solution to $P$. But Algorithm 8 produces a solution that achieves this upper bound so this solution must be optimal.
\hspace{\stretch{1}}$\blacksquare$

\begin{theorem} \label{thm:f=g} For any 3M-DAP $P = (T; U, V)$, $f(P) = g(P)$, i.e., Algorithm \ref{alg:recursive} produces an optimal solution.
\end{theorem}
\textbf{Proof.} The proof is by induction on $N$; the inductive hypothesis is that $f(P) = g(P)$ for any problem $P$ that is an $M$-problem with $M < N$. The induction is initialized for $P$ a 0-problem: then $g(P) = x_u/u$, which is an upper bound on $f(P)$ so the conclusion follows.

If $P(x)$ is a 1-problem with $x_b < x < x_{b-1}$, then $g(P) = x_u^{\prime}/u^{\prime} = z(b)$. Because any solution must contain a maximal inseparable $k$-pair with $k \ge b$ and $z(k)$ is decreasing in $k$, it follows that $z(b)$ is an upper bound on $f(P)$. Then the conclusion follows.

If $P(x)$ is a 2-problem with $x_b < x < x_{b-1}$ and $x^{\prime}_1 < x^{\prime} < x^{\prime}_0$, appeal to Lemma \ref{lem:f=g_N=2_b'=1}.

If $P(x)$ is a 2-problem with $x_1 < x < x_0$, appeal to Lemma \ref{lem:f=g_N=2_b=1}.

Otherwise, $P$ is either an $N$-problem with $N \ge 3$ or is a 2-problem with $b > 1$ and for which the reduced problem $P^{\prime}$ has $b^{\prime} > 1$. Then both Propositions \ref{prop:g>z(b+1)} and \ref{prop:g>(x_u - z(b-2))/(u - 1)} apply. Conclude from the former that $f(P) \ge g(P) > z(k)$ for all $k > b$ and from the latter that $f(P) \ge g(P) > (x_u - z(k))/(u - 1)$ for all $k < k - 1$. Together, these two results imply that an optimal solution can contain only $b$-pairs and $(b - 1)$-pairs. But the optimal such solution is given by Algorithm \ref{alg:recursive} (applying the inductive hypothesis). The conclusion follows.
\hspace{\stretch{1}}$\blacksquare$

\section{Complexity of Algorithm \ref{alg:recursive}}

Here we show that the running time of Algorithm \ref{alg:recursive} is linear in the size of the problem, i.e., in $n_t$, the number of elements in the matrix $T$.

To begin, consider $P$ an $N$-problem with $N > 1$. Focus on using Algorithm \ref{alg:(T,V)-pair} to complete the solution to $P$ given the solution to the reduced problem.

\begin{lemma} \label{lem:P_from_P'_complexity} Consider $P$ an $N$-problem with $N > 1$. Constructing the solution to $P$ given the solution to the reduced problem through application of Algorithm \ref{alg:(T,V)-pair} for each $(T, V)$-pair in the solution has complexity $\Theta(n_t)$.
\end{lemma}
\textbf{Proof.} Suppose that the solution to $P$ comprises the set of maximal, inseparable $(T, V)$-pairs $\{(A_1, B_1), \ldots, (A_p, B_p)\}$, where $n_{a_i}$ is the number of elements in $A_i$ for each $i = 1, \ldots, p$. Then $n_t = \sum_{i=1}^p n_{a_i}$. Completing each $(T, V)$-pair $(A_i, B_i)$ requires between $n_{a_i}$ and $3n_{a_i}$ operations by Lemma \ref{lem:alg_(T,V)-pair_complexity}. (Note that any $(T, V)$-pair with $v = 1$ will have $b \le t - 2$ by Lemma \ref{lem:v1}.)

Then constructing the solution to $P$ given the solution to the reduced problem through multiple applications of Algorithm \ref{alg:(T,V)-pair} requires between $n_t$ and $3n_t$ operations and the conclusion follows.
\hspace{\stretch{1}}$\blacksquare$

\begin{proposition} \label{prop:u>2_complexity} Consider a 3M-DAP $P$ with $t \ge 2$, $u > 2$, and $v \ge 2$. Solving $P$ via Algorithm \ref{alg:recursive} has complexity $O(n_t \log_2 \log_2 n_t)$.
\end{proposition}
\textbf{Proof.} $P$ is either solved with less than three reductions, which requires linear time, or in three or more reductions. In the latter case we find an upper bound on the number of reductions. Consider:
\begin{align*}
  t^{\prime} &= u\\
  u^{\prime} &= [(v - 1)b + 1]t - bv = [(v - 1)t - v]b + t \ge v(t - 1)\\
  v^{\prime} &= [(v - 1)(b - 1) + 1]t - (b - 1)v = [(v - 1)t - v](b - 1) + t \ge t.
\end{align*}
Write $t^1 = t^{\prime}$, $t^2 = t^{\prime \prime}$, etc. Then
\begin{align*}
  t^1 &= u\\
  u^1 &\ge v(t - 1)\\
  v^1 &\ge t\\
  t^2 &= u^1 \ge v(t - 1)\\
  u^2 &\ge v^1(t^1 - 1) \ge t(u - 1)\\
  v^2 &\ge t^1 = u\\
  t^3 &= u^2 \ge t(u - 1)\\
  u^3 &\ge v^2(t^2 - 1) \ge uv(t - 1)\\
  v^3 &\ge t^2 \ge v(t - 1).
\end{align*}
Because $u > 2$, $t^3 \ge t(u - 1) \ge 2t$. Then after three reductions, the new $T$ matrix is guaranteed to have at least twice as many columns as the original $T$ matrix. Because the new $T$ matrix has fewer elements than the original $T$ matrix, we can conclude that the number of rows in the new $T$ matrix is less than half the number of rows in the original $T$ matrix.

In fact, we can be more aggressive. Write $w_0 = w = tu/2$, $t_1 = t^3 \ge t(u - 1) \ge tu/2 = w$, and $u_1 = u^3 \ge uv(t - 1) \ge tu = 2w$. Then $w_1 = t_1u_1/2 \ge w^2$, $w_2 \ge w_1^2 \ge w^4$, $w_3 \ge w_2^2 \ge w^8$, and $w_i \ge w^{2^i}$ for all $i \ge 0$. Then
\[
  t_i \ge w_{i-1} \ge w^{2^{i-1}} = \left(\frac{tu}{2}\right)^{2^{i-1}}.
\]
Then
\[
  n_t > t_is_{t_i} > \left(\frac{tu}{2}\right)^{2^{i - 1}}s_{t_i} \quad \Leftrightarrow \quad s_{t_i} < \left(\frac{2}{tu}\right)^{2^{i - 1}} n_t < t^{-2^{i - 1}}n_t.
\]
Here $i$ counts the number of steps, with each step corresponding to three reductions. The maximum number of steps is bounded above by the value $i^*$ that ensures that $s_{t_{i^*}}$ is at most 1 when we set $t=2$:
\[
  2^{-2^{{i^* - 1}}}n_t \le 1 \quad \Leftrightarrow \quad n_t \le 2^{2^{{i^* - 1}}} \quad \Leftrightarrow \quad \log_2 \log_2 n_t \le i^* - 1.
\]
Then number of reductions is $O(\log_2 \log_2 n_t)$. Because solving the final 0-problem and constructing the solution to any problem from its reduced problem can done in time linear in the number of elements in the problem(by Lemmas \ref{lem:alg_0prob_type1_complexity} and \ref{lem:P_from_P'_complexity}, respectively), which is bounded above by $n_t$, the complexity of the algorithm is at worst $O(n_t \log_2 \log_2 n_t)$.
\hspace{\stretch{1}}$\blacksquare$

\begin{lemma} \label{lem:t>2_or_v>2_complexity} Consider a 3M-DAP $P$ with $u \ge 2$ and either $t > 2$ and $v \ge 2$, or $t \ge 2$ and $v > 2$. Solving $P$ via Algorithm \ref{alg:recursive} has complexity $O(n_t \log_2 \log_2 n_t)$.
\end{lemma}
\textbf{Proof.} $P$ is either a 0-problem, which can be solved in linear time by Lemma \ref{lem:alg_0prob_type1_complexity}, or can be reduced. The reduced problem has $t^{\prime} = u \ge 2$, $u^{\prime} \ge v(t - 1) > 2$, and $v^{\prime} \ge t \ge 2$. By Proposition \ref{prop:u>2_complexity}, solving the reduced problem has time complexity at most $O(n_t \log_2 \log_2 n_t)$, while constructing the solution to $P$ from the solution to $P^{\prime}$ has linear time complexity by Lemma \ref{lem:P_from_P'_complexity}.
\hspace{\stretch{1}}$\blacksquare$

\begin{lemma} \label{lem:v=1_complexity} Consider a 3M-DAP $P$ with $v = 1$. Solving $P$ via Algorithm \ref{alg:recursive} has complexity $O(n_t \log_2 \log_2 n_t)$.
\end{lemma}
\textbf{Proof.} By assumption $P$ has $t > 2$ and $u \ge 2$. Now $P$ is either a 0-problem, which can be solved in linear time by Lemma \ref{lem:alg_0prob_type1_complexity}, or can be reduced. By Proposition \ref{prop:reduced_problem_properties}, the reduced problem has $t^{\prime} = u \ge 2$, $u^{\prime} \ge v(t - 1) \ge 2$, and $v^{\prime} = u^{\prime} + 1 > 2$. By Lemma \ref{lem:t>2_or_v>2_complexity}, solving the reduced problem has time complexity at most $O(n_t \log_2 \log_2 n_t)$, while constructing the solution to $P$ from the solution to $P^{\prime}$ has linear time complexity by Lemma \ref{lem:P_from_P'_complexity}.
\hspace{\stretch{1}}$\blacksquare$

\begin{lemma} \label{lem:t=u=v=2_complexity} Consider a 3M-DAP $P$ with $t=u=v = 2$. Solving $P$ via Algorithm \ref{alg:recursive} has complexity $\Theta(n_t)$.
\end{lemma}
\textbf{Proof.} If $P$ is a 0-problem, it can be solved in linear time by Lemma \ref{lem:alg_0prob_type1_complexity}. If $P$ is a 1-problem, its reduced problem $P^{\prime}$ is a 0-problem, which can be solved in linear time by Lemma \ref{lem:alg_0prob_type1_complexity}, and the solution to $P$ can be constructed in linear time from the solution to $P^{\prime}$ by Lemma \ref{lem:P_from_P'_complexity}.

Otherwise, $P$ is an $N$-problem with $N \ge 2$. Then $s_t^{\prime \prime} = s_u^{\prime} < s_t^{\prime} = s_u$. Because $s_u + s_v = s_t$, if $s_u < s_v$, then $s_t^{\prime \prime} < s_u < s_t/2$. If $s_u > s_v$, then
\[
  \hat{b} = \left \lceil \frac{s_v}{s_t - s_v} \right \rceil = \left \lceil \frac{s_v}{s_u} \right \rceil = 1.
\]
Then $s_t^{\prime \prime} = s_u^{\prime} = s_v - (\hat{b} - 1)(s_t - s_v) = s_v < s_t/2$. Then
\[
  n_t^{\prime \prime} = 2s_t^{\prime \prime} < s_t = \frac{1}{2}n_t.
\]
Let $T(n_t)$ be the time complexity of a 3M-DAP with $t=u=v=2$ and $n_t$ elements in matrix $T$, and let $f(n_t)$ be the time complexity to construct the solution to the original problem from the doubly-reduced problem. Then we obtain the recursion
\[
  T(n_t) \le T(n_t/2) + f(n_t).
\]
By Lemma \ref{lem:P_from_P'_complexity}, $f(n_t) = \Theta(n_t)$. Then it follows from the Master Theorem (Cormen, et al. 2001) that $T(n_t) = \Theta(n_t)$.
\hspace{\stretch{1}}$\blacksquare$

We summarize the results of this section in the following theorem.
\begin{theorem} \label{thm:alg_recursive_complexity} Consider a 3M-DAP $P$. Solving $P$ via Algorithm \ref{alg:recursive} has complexity $\Theta(n_t)$ when $t=u=v = 2$, otherwise has complexity $O(n_t \log_2 \log_2 n_t)$.
\end{theorem}

\section{Scott Huddleston's Algorithm for solving a 3M-DAP} \label{sec:scott}

In 2010 Scott Huddleston developed an algorithm for solving muffin problems. The algorithm is described in Chapter 13 of Gasarch et al. (2020) using a rather different framework than is employed herein. However, the two algorithms are essentially identical (as such, Huddleston's algorithm optimally solves any 3M-DAP). Algorithm \ref{alg:recursive_scott} below presents Huddleston's algorithm in the terminology of this paper. There two situations in which Algorithms \ref{alg:recursive} and \ref{alg:recursive_scott} employ different strategies, which do however lead to identical results:
\begin{enumerate}
\item When applied to a 0-problem of type 1, Algorithm \ref{alg:recursive} creates a subproblem that is a 3M-DAP and then applies either Algorithm \ref{alg:u=v+1} or Algorithm \ref{alg:t_even} to solve the subproblem. These two algorithms have lower bound guarantees that demonstrate that the solution produced for the subproblem has value exceeding $x_u/u$. But given that we have now established that Algorithm \ref{alg:recursive} produces an optimal solution (even when $q_t = t - 2$ and $r_t > 0$, so that subproblem created has $v^{\prime} = 1$), we can simply use this algorithm to solve the subproblem. This revision is reflected in Algorithm \ref{alg:recursive_scott}. 
\item When applied to a 0-problem of type 2, Algorithm \ref{alg:recursive} employs Algorithm \ref{alg:0prob_type2} to solve the problem directly. However, as described in Section \ref{sec:reducing_0prob_type2}, it is possible to reduce such a problem. The reduced problem is a 0-problem of type 1 and has a trivial solution. As with any other problem that is reduced, the solution to the original problem can be constructed from the solution to the reduced problem by applying Algorithm \ref{alg:(T,V)-pair}. This is the approach employed by Algorithm \ref{alg:recursive_scott}.
\end{enumerate}

\fontsize{8}{8}
\begin{algorithm}
\SetKwInOut{Input}{input}
\SetKwInOut{Output}{output}
\Input{$T$, $U$, $V$, empty matrices of dimensions $s_t \times t$, $s_u \times u$, and $s_v \times v$, respectively\\
$x_t$, $x_u$, $x_v$, the required row sums for each row of $T$, $U$, and $V$, respectively\\
$s_uu + s_vv = s_tt$\\
$s_ux_u + s_vx_v = s_tx_t$\\
}
\Output{the filled matrices $T$, $U$, $V$\\
the multiset of entries in $T$ is the union of the multiset of entries in $U$ and $V$\\
}
\BlankLine

\uIf{$x \le x_{\infty}$}{
	set all elements of $U$ to value $x_u/u$\\
	$q_t \leftarrow \lfloor n_u/s_t \rfloor$\\
	$r_t \leftarrow n_u - q_ts_t$\\
	\uIf{$r_t = 0$}{
		set all elements of $V$ to value $x_v/v$\\
		fill each row of $T$ with $q_t$ elements of value $x_u/u$ and $t - q_t$ elements of value $x_v/v$\\
	}
	\uElse{
		for $r_t$ rows of $T$, set $q_t + 1$ elements to value $x_u/u$; call the unfilled submatrix $V^{\prime}$\\
		set $x_v^{\prime} = x_t - (q_t + 1)x_u/u$\\
		for the remaining $s_t - r_t$ rows of $T$, set $q_t$ elements to value $x_u/u$; call the unfilled submatrix $U^{\prime}$\\
		set $x_u^{\prime} = x_t - q_tx_u/u$\\
		set $T^{\prime} = V$\\
		apply Algorithm \ref{alg:recursive_v2} to solve the $(T^{\prime}, U^{\prime}, V^{\prime})$ 3M-DAP\\
	}
}
%\uElseIf{$x = x_b$ for some integer $b \ge 0$}{
%	apply Algorithm \ref{alg:0prob_type2} to solve the $(T; U, V)$ problem\\
%}
\Else{
	$b \leftarrow \left \lceil s_v/[s_t - (v - 1)s_v] \right \rceil$\\
	form the reduced problem $(T^{\prime}, U^{\prime}, V^{\prime})$ as described in Section \ref{sec:reduced_problem}\\
	apply Algorithm \ref{alg:recursive_v2} to solve the $(T^{\prime}, U^{\prime}, V^{\prime})$ problem\\
	$U \leftarrow T^{\prime}$\\
	\For{$1 \le i \le s_u^{\prime}$}{
		initialize a $b$-pair $(A, B)$\\
		set $\boldsymbol{u}$ to be the $i^{th}$ row of $U^{\prime}$\\
		apply Algorithm \ref{alg:(T,V)-pair} with inputs $A$, $x_t$, $B$, $x_v$, and $\boldsymbol{u}$ to complete the pair $(A, B)$\\
		insert $A$ into $T$ and $B$ into $V$\\
	}

	\uIf{$b = 1$}{
		insert $V^{\prime}$ into $T$\\
	}
	\uElse{
		\For{$1 \le j \le s_v^{\prime}$}{
			initialize a $(b - 1)$-pair $(A, B)$\\
			set $\boldsymbol{u}$ to be the $j^{th}$ row of $V^{\prime}$\\
			apply Algorithm \ref{alg:(T,V)-pair} with inputs $A$, $x_t$, $B$, $x_v$, and $\boldsymbol{u}$ to complete the pair $(A, B)$\\
			insert $A$ into $T$ and $B$ into $V$\\
		}
	}
}

\caption{Scott Huddleston's Algorithm for solving a 3M-DAP}
\label{alg:recursive_scott}
\end{algorithm}

\normalsize

\section{Further Results for Muffin Problems}

For the remainder of the paper, we provide further analysis of muffin problems. 

In Section \ref{sec:n_or_n+1_conj}, we prove that Conjecture \ref{n_or_n+1_conj} does indeed hold, namely that for any supply-constrained muffin problem, there exists an optimal solution in which each student receives either $n$ or $n + 1$ pieces. 

In Section \ref{sec:1/3}, we investigate when $f(m, s) = \sfrac{1}{3}$. This is the case when $g(m, s) \le \sfrac{1}{3}$, which occurs only when $n = 2$ and the problem is a 1-problem of type 1, i.e., the reduced problem is a 0-problem of type 1.

Any other fully-constrained muffin problem of order 2 can be doubly reduced and then standardized to form a related fully-constrained muffin problem of order greater than 2. In Section \ref{sec:relationships} we investigate the nature of these relationships between different muffin problems. Section \ref{sec:relationship_2_to_>2} shows that for any muffin problem of order greater than 2 and for every $b \ge 1$, there is a muffin problem $P(x)$ of order 2 with $x_b < x < x_{b - 1}$ that is related to it. It follows that there are relationships between muffin problems of order 2: we begin investigation of these relationships in Section \ref{sec:relationship_2_to_2}. We demonstrate that it is possible to construct a solution to a muffin problem of order 2 from the solution to a related muffin problem of order 2 by reducing (Section \ref{sec:alternative_reduction}) or expanding (Section \ref{sec:equivalent_expansion}), as appropriate.

Finally, in Section \ref{conj1_proof}, we prove Conjecture \ref{conj:f(km,ks)=f(m,s)}, that the value of a muffin problem depends only on $x = \sfrac{m}{s}$, the ratio of the numbers of muffins to students.

\subsection{Proof of Conjecture \ref{n_or_n+1_conj}} \label{sec:n_or_n+1_conj}

Assume $x$ not an integer or half-integer. Consider a supply-constrained muffin problem (so that each muffin must be cut into two pieces). Then some student must receive at least $n + 1$ pieces, giving an upper bound of $x/(n + 1)$. Equally, some student must receive at most $n$ pieces, giving an upper bound of $1 - x/n$.

If we require that each student receive either $n$ or $n + 1$ pieces, then the problem is a fully-constrained muffin problem. In this section we prove Conjecture \ref{n_or_n+1_conj}, which claims that any optimal solution to a fully-constrained muffin problem is also optimal for the supply-constrained version of the problem.

A fully-constrained muffin problem can be formulated as a 3M-DAP, with $t = 2$, $s_t = m$, $u = n + 1$, $v = n$, $x_t = 1$, $x_u = x$, and $x_v = x$. We must have
\[
\begin{split}
  s_u + s_v &= s, \quad \text{and}\\
  (n + 1)s_u + ns_v &= 2m,
\end{split}
\]
so
\[
\begin{split}
  s_u &= 2m - ns, \quad \text{and}\\
  s_v &= (n + 1)s - 2m.
\end{split}
\]
We have $\lambda = x_u - x_v = 0$ and $\gamma = x_t/t = 1/2$. Further,
\[
  x_{\infty} = \frac{[\lambda + \gamma(v - 1)t]u}{(v - 1)t + u - v} = \frac{(n - 1)(n + 1)}{2n - 1} = \frac{n^2 - 1}{2n - 1},
\]
and
\[
  x_1 = \frac{[\lambda + \gamma(v - 1)t]u + \gamma tu}{[(v - 1)t + u - v] + t} = \frac{n(n + 1)}{2n + 1}.
\]

\begin{lemma} \label{lem:x_for_N=1_b'=1} A fully-constrained muffin problem $P(x)$ is a 1-problem of type 1 (i.e., the reduced problem is a 0-problem of type 1) with $x_1 < x < x_0$ if and only if 
\[
  \frac{n^2}{2n - 1} \le x < x_0.
\]
\end{lemma}
\textbf{Proof.} Suppose $P(x)$ is a 1-problem of type 1 with $x_1 < x < x_0$. Because $b=1$, the reduced problem has:
\begin{align*}
  x^{\prime} = x_u^{\prime} &= vx_t - x_v = n - x\\
  u^{\prime} &= vt - v = n\\
  x_v^{\prime} &= x_t = 1\\
  v^{\prime} &= t = 2\\
  x_t^{\prime} &= x_u = x\\
  t^{\prime} &= u = n + 1.
\end{align*}
Then
\[
  x_{\infty}^{\prime} = \frac{[(v^{\prime} - 1)x_t^{\prime} + x_u^{\prime} - x_v^{\prime}]u^{\prime}}{(v^{\prime} - 1)t^{\prime} + u^{\prime} - v^{\prime}} = \frac{n(n - 1)}{2n - 1}.
\]
Because $P$ is a 1-problem of type 1, the reduced problem is a 0-problem of type 1 so has $x^{\prime} \le x_{\infty}^{\prime}$. Then
\[
  x^{\prime} \le x_{\infty}^{\prime} \Leftrightarrow n - x \le \frac{n(n - 1)}{2n - 1} \Leftrightarrow x \ge \frac{n^2}{2n - 1}.
\]
The converse follows immediately.
\hspace{\stretch{1}}$\blacksquare$

A solution to the 3M-DAP gives a feasible solution and hence a lower bound to the supply-constrained muffin problem. We will demonstrate in this section that the solution is in fact optimal. Indeed, in most cases the only optimal solution is the solution to the 3M-DAP, i.e., for most supply-constrained muffin problems in an optimal solution each student must receive either $n$ or $n + 1$ pieces.

\begin{lemma} \label{lem:3mat_opt_N=0} If a fully-constrained muffin problem $P$ is a 0-problem, then the only optimal solutions to $Q$, the supply-constrained version of the problem, are the optimal solutions to $P$.
\end{lemma}
\textbf{Proof.} The optimal solution to $P$ is $x_u/u = x/(n + 1)$. Because this is an upper bound on the solution to $Q$, the solution is optimal for $Q$. 

Suppose a student receives $n + 2$ or more pieces. Then an upper bound on the solution is $x/(n + 2) < x/(n + 1)$, so such a solution cannot be optimal.

Suppose a student receives $n - 1$ or fewer pieces. Then an upper bound on the solution is $1 - x/(n - 1)$. But
\[
\begin{split}
  x &> \frac{n}{2} > \frac{n^2 - 1}{2n}\\
  \Leftrightarrow \frac{2n}{n^2 - 1}x &> 1\\
  \Leftrightarrow \frac{x}{n - 1} + \frac{x}{n + 1} &> 1\\
  \Leftrightarrow \frac{x}{n + 1} &> 1 - \frac{x}{n - 1},
\end{split}
\]
so such a solution cannot be optimal.
\hspace{\stretch{1}}$\blacksquare$

\begin{lemma} \label{lem:3mat_opt_N=1_b=2} If a fully-constrained muffin problem $P(x)$ is a 1-problem with $x_2 < x < x_1$, then the only optimal solutions to $Q(x)$, the supply-constrained version of the problem, are the optimal solutions to $P(x)$.
\end{lemma}
\textbf{Proof.} To start, we have
\[
  x < x_1 = \frac{n(n + 1)}{2n + 1}.
\]
Because $P(x)$ is a 1-problem with $x_2 < x < x_1$, the reduced problem is formed using only $2$- and $1$-pairs. Then
\[
  g(P) = \frac{x_u^{\prime}}{u^{\prime}} = \frac{(v - 1/2)x_t - x_v}{(v - 1/2)t - v} = \frac{n - 1/2 - x}{n - 1} = 1 - \frac{x - 1/2}{n - 1}.
\]
First, suppose a student receives $n - 1$ or fewer pieces. Then an upper bound on the solution is $1 - x/(n - 1)$. But we have
\[
  x - 1/2 < x \Leftrightarrow \frac{x - 1/2}{n - 1} < \frac{x}{n - 1} \Leftrightarrow 1 - \frac{x}{n - 1} < 1 - \frac{x - 1/2}{n - 1},
\]
so such a solution cannot be optimal. Second, suppose a student receives $n + 2$ or more pieces. Then an upper bound on the solution is $x/(n + 2)$. Because $n \ge 2$,
\[
\begin{split}
  x < x_1 =\frac{n(n + 1)}{2n + 1} &\le \frac{(2n - 1)(n + 2)}{2(2n + 1)}\\
  \Leftrightarrow (4n + 2)x &< (2n - 1)(n + 2)\\
  \Leftrightarrow 2(n - 1)x &< 2(n - 1)(n + 2) - (n + 2)(2x - 1)\\
  \Leftrightarrow \frac{x}{n + 2} &< 1 - \frac{x - 1/2}{n - 1},
\end{split}
\]
so such a solution cannot be optimal. 

Then in any optimal solution each student must be assigned either $n$ or $n + 1$ pieces. The conclusion follows.
\hspace{\stretch{1}}$\blacksquare$

\begin{lemma} \label{lem:3mat_opt_n=2_N=1_b>2} If a fully-constrained muffin problem $P(x)$ with $n = 2$ is a 1-problem with $x_b < x < x_{b-1}$ where $b \ge 2$, then the only optimal solutions to $Q(x)$, the supply-constrained version of the problem, are the optimal solutions to $P(x)$.
\end{lemma}
\textbf{Proof.} Because $n = 2$, we have $t = v = 2$ and $u = 3$. Then
\[
  x < x_{b - 1} = \frac{[\lambda + \gamma(v - 1)t]u(b - 1) + \gamma tu}{[(v - 1)t + u - v](b - 1) + t} = \frac{3b}{3b - 1}.
\]
Because $P$ is a 1-problem, we have
\[
  g(P) = \frac{x_u^{\prime}}{u^{\prime}} = \frac{(v - 1 + 1/b)x_t - x_v}{(v - 1 + 1/)t - v} = \frac{(1 + 1/b) - x}{2/b} = \frac{1}{2}[1 - b(x - 1)].
\]

First, a student cannot receive only one piece because the size of that piece is at most 1, which is less than $x$. Second, suppose a student receives four or more pieces. Then an upper bound on the solution is $x/4$. Because $b \ge 2$,
\[
\begin{split}
  x &< \frac{3b}{3b - 1} = \frac{2b + 2 + b - 2}{2b + 1 + b - 2} \le \frac{2(b + 1)}{2b + 1}\\
  \Leftrightarrow (2b + 1)x &< 2(b + 1)\\
  \Leftrightarrow x &< 2[1 - b(x - 1)]\\
  \Leftrightarrow \frac{x}{4} &< \frac{1}{2}[1 - b(x - 1)],
\end{split}
\]
so such a solution cannot be optimal. 

Then in any optimal solution each student must be assigned either two or three pieces. The conclusion follows.
\hspace{\stretch{1}}$\blacksquare$

\begin{lemma} \label{lem:3mat_opt_N>0} If a fully-constrained muffin problem $P(x)$ is either (i) a 1-problem with $n > 2$ and $x_b < x < x_{b-1}$ where $b \ge 3$; or (ii) an $N$-problem with $N \ge 2$, then the only optimal solutions to $Q(x)$, the supply-constrained version of the problem, are the optimal solutions to $P(x)$.
\end{lemma}
\textbf{Proof.} If $P$ is a 1-problem with $n > 2$ and $x_b < x < x_{b-1}$ where $b \ge 3$, then $g(P) > \tau$ by Lemma \ref{lem:L(b)_a>0}. Alternatively, if $P$ is an $N$-problem with $N \ge 2$, then $g(P) > \tau$ by Theorem \ref{thm:g}. So
\[
  g(P) > \tau = \frac{x_{\infty}}{u} = \frac{n - 1}{2n - 1}.
\]

Suppose a student receives $n - 1$ or fewer pieces. Then an upper bound on the solution is $1 - x/(n - 1)$. Because $P$ is not a 0-problem, $x > x_{\infty} = (n^2 - 1)/(2n - 1)$, so we have
\[
  1 - \frac{x}{n - 1} < 1 - \frac{n^2 - 1}{(2n - 1)(n - 1)} = 1 - \frac{n + 1}{2n - 1} = \frac{n - 2}{2n - 1} < \frac{n - 1}{2n - 1},
\]
so such a solution cannot be optimal. 

Alternatively, suppose a student receives $n + 2$ or more pieces. Then an upper bound on the solution is $x/(n + 2)$. Because $P$ is not a 1-problem with $x_1 < x < x_0$, from Lemma \ref{lem:x_for_N=1_b'=1}, $x < n^2/(2n - 1)$, so we have
\[
\begin{split}
  n &\ge 2\\
  \Leftrightarrow n^2 + n - 2 = (n - 1)(n + 2) &\ge n^2\\
  \Leftrightarrow \frac{n - 1}{2n - 1} &\ge \frac{n^2}{(2n - 1)(n + 2)}\\
  \Rightarrow \frac{n - 1}{2n - 1} &> \frac{x}{n + 2},
\end{split}
\]
so such a solution cannot be optimal. 

Then in any optimal solution each student must be assigned either $n$ or $n + 1$ pieces. The conclusion follows.
\hspace{\stretch{1}}$\blacksquare$

Finally, we address 1-problems with $x_1 < x < x_0$. As we will see, for this type of supply-constrained muffin problem, it is possible to find optimal solutions that are not solutions to the fully-constrained muffin problem. 

\begin{lemma} \label{lem:3mat_opt_N=1_b=1_part1} If a fully-constrained muffin problem $P(x)$ is a 1-problem with $x_1 < x < x_0$, then any optimal solution to $P(x)$ is optimal for $Q(x)$, the supply-constrained version of the problem.
\end{lemma}
\textbf{Proof.} The optimal solution to $P$ is
\[
  \frac{x_u^{\prime}}{u^{\prime}} = \frac{vx_t - x_v}{vt - v} = \frac{n - x}{n} = 1 - \frac{x}{n}.
\]
Because this is an upper bound on the solution to the $Q$, the solution is optimal for $Q$.
\hspace{\stretch{1}}$\blacksquare$

\begin{lemma} \label{lem:claim_num_pieces} For a muffin problem $(m, s)$ with
\[
  \frac{n + 1}{2} - \frac{1}{2(n + 1)} = \frac{n(n + 2)}{2(n + 1)} \le x < \frac{n + 1}{2},
\]
and either $n \ge 3$, or $n = 2$ and $x > 4/3$, we have
\[
  x \ge \frac{n^2}{2n - 1} + \frac{n - 1}{(2n - 1)s} > \frac{n}{2} + \frac{3}{2s}.
\]
\end{lemma}
\textbf{Proof.} First, suppose $n \ge 3$. When $n$ is odd, we must have $s \ge 2(n + 1)$, and when $n$ is even, we must have $s \ge n + 1$. Note that $2(n - 1)/[n(n - 2)]$ is less than 2 for all $n \ge 3$ and is less than 1 for all $n \ge 4$. Then it follows that
\[
\begin{split}
  s &> \frac{2(n - 1)}{n(n - 2)}(n + 1)\\
  \Leftrightarrow \frac{1}{s} &< \frac{n(n - 2)}{2(n - 1)(n + 1)}\\
  \Leftrightarrow \frac{1}{s} &< \frac{n(n + 2)(2n - 1)}{2(n - 1)(n + 1)} - \frac{n^2}{n - 1}\\
  \Leftrightarrow \frac{n^2}{2n - 1} + \frac{n - 1}{(2n - 1)s} &< \frac{n(n + 2)}{2(n + 1)} \le x.
\end{split}
\]
Second, suppose that $n = 2$ and $x > 4/3$. We want to show that $x \ge 4/3 + 1/(3s)$. Write $d = m - s$, $m = 3d + a$, and $s = 2d + a$. Because $4/3 < x < 3/2$, we must have $1 \le a < d$. For a fixed value of $d$, $x$ is decreasing in $a$, so it suffices to consider $a = d - 1$. Then $m = 4d - 1$ and $s = 3d - 1$, so $x = 4/3 + 1/(3s)$, as required.

For the second part, first note that $s \ge 5$ (when $n \ge 4$, $s \ge n + 1$, when $n = 3$, $s \ge 8$, and when $n = 2$, because $4/3 < x < 3/2$, $s \ge 5$). Then
\[
\begin{split}
  (s - 4)n &> -1\\
  \Leftrightarrow ns - 1 &> 4n - 2\\
  \Leftrightarrow \frac{ns - 1}{2n - 1} &> 2\\
  \Leftrightarrow \frac{ns + 2n - 2}{2n - 1} &> 3\\
  \Leftrightarrow \frac{n}{2} + \frac{ns + 2n - 2}{2(2n - 1)s} &> \frac{n}{2} + \frac{3}{2s}\\
  \Leftrightarrow \frac{n^2}{2n - 1} + \frac{n - 1}{(2n - 1)s} > \frac{n}{2} + \frac{3}{2s}.
\end{split}
\]
\hspace{\stretch{1}}$\blacksquare$

\fontsize{8}{8}
\begin{algorithm}
\SetKwInOut{Input}{input}
\SetKwInOut{Output}{output}
\Input{$T$, an empty matrix of dimensions $m \times 2$\\
$x_t = 1$, the required row sum for each row of $T$\\
$\boldsymbol{u}$, an empty demand vector with length $n + 2$\\
$U$, an empty matrix of dimensions $(2m - ns - 2) \times (n + 1)$\\
$V$, an empty matrix of dimensions $[(n + 1)s - 2m + 1] \times n$\\
$x$, the required row sum for all demand vectors\\
\quad \quad \quad \quad (the vector $\boldsymbol{u}^*$ and all rows of matrices $U$ and $V$)\\
$m = sx$\\
}
\Output{the filled matrices $T$, $U$, and $V$ and vector $\boldsymbol{u}$\\
the multiset of entries in $T$ is the multiset of entries in $\boldsymbol{u}$, $U$, and $V$\\
}
\BlankLine

divide $n[(n + 1)s - 2m + 1]$ rows of $T$ $[x/n; 1 - x/n]$\\
divide one row of $T$ $[(2 + 1/n)x - (n + 1), (n + 2) - (2 + 1/n)x]$\\
set all elements of $V$ to $x/n$\\
set $n + 1$ elements of $\boldsymbol{u}$ to $1 - x/n$ and the remaining element to $(2 + 1/n)x - (n + 1)$\\
$q \leftarrow \lfloor \{n[(n + 1) - 2m] - 1\}/(2m - ns - 2) \rfloor$\\
$r \leftarrow n[(n + 1) - 2m] - 1 - q(2m - ns - 2)$\\
assign the remaining $n[(n + 1)s - 2m + 1] - (n + 1)$ pieces of size $1 - x/n$ and the piece of size $(n + 2) - (2 + 1/n)x$ to the rows of $U$ so that:\\
\quad $2m - ns - 3 - r$ rows each receive $q$ pieces of size $1 - x/n$\\
\quad \quad let $U_1^{\prime}$ be the unfilled submatrix of dimensions $(2m - ns - 3 - r) \times (n + 1 - q)$\\
\quad \quad set $x_1^{\prime} = x - q(1 - x/n)$\\
\quad one row receives $q$ pieces of size $1 - x/n$ and one piece of size $(n + 2) - (2 + 1/n)x$\\
\quad \quad let $\boldsymbol{u}_2^{\prime}$ be the vector consisting of the $n - q$ unfilled elements\\
\quad \quad set $x_2^{\prime} = x - q(1 - x/n) - [(n + 2) - (2 + 1/n)x]$\\

\quad $r$ rows each receive $q + 1$ pieces of size $1 - x/n$\\
\quad \quad let $U_3^{\prime}$ be the unfilled submatrix of dimensions $r \times (n - q)$\\
\quad \quad set $x_3^{\prime} = x - (q + 1)(1 - x/n)$\\
set $T^{\prime}$ to be the matrix consisting of the unfilled rows of $T$\\
use Algorithm \ref{alg:t_even_gen} to solve the $(T^{\prime}; U_1^{\prime}, \boldsymbol{u}_2^{\prime}, U_3^{\prime})$ DAP\\
\caption{Algorithm to construct an optimal solution for $Q(x)$ that is not feasible for $P(x)$}
\label{alg:opt_for_Q_not_P}
\end{algorithm}

\normalsize

\begin{lemma} \label{lem:3mat_opt_N=1_b=1_part2} If a fully-constrained muffin problem $P(x)$ is a 1-problem with $x_1 < x < x_0$, then:
\begin{enumerate}
\item[(1)] when $x < n(n + 2)/[2(n + 1)]$, the only optimal solutions to $Q(x)$, the supply-constrained version of the problem, are the optimal solutions $P(x)$;
\item[(2)] when $n(n + 2)/[2(n + 1)] \le x < (n + 1)/2$, there exist optimal solutions to $Q(x)$ that are not feasible for $P(x)$.
\end{enumerate} 
\end{lemma}
\textbf{Proof.} The optimal solution to $P$ is
\[
  \frac{x_u^{\prime}}{u^{\prime}} = \frac{vx_t - x_v}{vt - v} = \frac{n - x}{n} = 1 - \frac{x}{n}.
\]
A solution to $Q$ in which a student has less than $n$ pieces has an upper bound of $1 - x/(n - 1) < 1 - x/n$, so cannot be optimal. A solution in which a student has $n + 2$ or more pieces has an upper bound of $x/(n + 2)$. Now
\[
  \frac{n + 1}{2} - \frac{1}{2(n + 1)} = \frac{n(n + 2)}{2(n + 1)} > x \Leftrightarrow n(n + 2) > 2(n + 1)x \Leftrightarrow 1 - \frac{x}{n} > \frac{x}{n + 2}.
\]
Thus, on this range the solution cannot be optimal. This demonstrates the first result.

For the second result, suppose that
\[
  \frac{n + 1}{2} - \frac{1}{2(n + 1)} = \frac{n(n + 2)}{2(n + 1)} \le x < \frac{n + 1}{2}.
\]
First, assume that $x \neq 4/3$, so that Lemma \ref{lem:claim_num_pieces} applies. 

Algorithm \ref{alg:opt_for_Q_not_P} constructs an optimal solution to the supply-constrained muffin problem $Q(x)$ in which one student has $n + 2$ pieces, $2m - ns - 2$ students have $n + 1$ pieces, and $(n + 1)s - 2m + 1$ students have $n$ pieces. By Lemma \ref{lem:claim_num_pieces}
\[
  x > \frac{n}{2} + \frac{3}{2s} \quad \Leftrightarrow \quad 2m - ns - 2 > 1,
\]
so there are at least two students who receive $n + 1$ pieces.

The initial division of muffins results in pieces that are all at least as large as $1 - x/n$, because
\[
\begin{split}
  x \ge \frac{n(n + 2)}{2(n + 1)} \quad &\Leftrightarrow \quad \left(2 + \frac{1}{n}\right)x - (n + 1) \ge 1 - \frac{x}{n}\\
  x < \frac{n + 1}{2} \quad &\Leftrightarrow \quad n + 2 - \left(2 + \frac{1}{n}\right)x > 1 - \frac{x}{n}.
\end{split}
\]

After assigning pieces to the $(n + 2)$-student and all the $n$-students, there are $n[(n + 1)s - 2m + 1] - (n + 1)$ pieces of size $1 - x/n$ and one piece of size $(n + 2) - (2 + 1/n)x$ to be assigned to the $(n + 1)$-students, making $n[(n + 1)s - 2m]$ pieces in total. This can be done so that each such student still needs at least two pieces, because by Lemma \ref{lem:claim_num_pieces}
\[
\begin{split}
  x &\ge \frac{n^2}{2n - 1} + \frac{n - 1}{(2n - 1)s}\\
  \Leftrightarrow (2n - 1)m &\ge n^2s + n - 1\\
  \Leftrightarrow (n - 1)[2m - ns - 2] &\ge n[(n + 1)s - 2m].
\end{split}
\]
This allows us to use Algorithm \ref{alg:t_even_gen} to solve the remaining subproblem.

It remains to show that the subproblem has a solution with smallest piece at least $1 - x/n$. Following the application of Algorithm \ref{alg:t_even_gen}, write
\[
\begin{split}
  x_1^{\prime \prime} &= x_1^{\prime} - (n - 1 - q)\frac{1}{2} = x - q\left(1 - \frac{x}{n}\right) - (n - 1 - q)\frac{1}{2}\\
  x_2^{\prime \prime} &= x_2^{\prime} - (n - 2 - q)\frac{1}{2} = x - q\left(1 - \frac{x}{n}\right) - [(n + 2) - (2 + 1/n)x] - (n - 2 - q)\frac{1}{2}\\
  x_3^{\prime \prime} &= x_3^{\prime} - (n - 2 - q)\frac{1}{2} = x - (q + 1)\left(1 - \frac{x}{n}\right) - (n - 2 - q)\frac{1}{2}.
\end{split}
\]
Now
\[
\begin{split}
  x_2^{\prime \prime} - x_1^{\prime \prime} &= \frac{1}{2} - [(n + 2) - (2 + 1/n)x]\\
  x_3^{\prime \prime} - x_1^{\prime \prime} &= \frac{x}{n} - \frac{1}{2} > 0\\
  x_3^{\prime \prime} - x_2^{\prime \prime} &= [(n + 2) - (2 + 1/n)x] - \left(1 - \frac{x}{n}\right) = (n + 1) - 2x > 0.
\end{split}
\]
It is possible that either $U_1^{\prime}$ or $U_3^{\prime}$ may not be present and which of $x_1^{\prime \prime}$ and $x_2^{\prime \prime}$ is larger depends on the value of $x$. Then there are several cases to consider. Recall that $x_{\min}^{\prime \prime} = \min\{x_1^{\prime \prime}, x_2^{\prime \prime}, x_3^{\prime \prime}\}$ and $x_{\max}^{\prime \prime} = \max\{x_1^{\prime \prime}, x_2^{\prime \prime}, x_3^{\prime \prime}\}$.

\textbf{Case 1.} Both $U_1^{\prime}$ and $U_3^{\prime}$ are present and $x_1^{\prime \prime} < x_2^{\prime \prime}$. Then $x_{\min}^{\prime \prime} = x_1^{\prime \prime}$ and $x_{\max}^{\prime \prime} = x_3^{\prime \prime}$, so $\lambda^{\prime \prime} = x_{\min}^{\prime \prime} - x_{\max}^{\prime \prime} = x_1^{\prime \prime} - x_3^{\prime \prime} = \frac{1}{2} - \frac{x}{n}$. Then by Lemma \ref{lem:LB_on_f;t_even_gen}, there is a solution to the problem of assigning the remaining pieces with smallest piece having size greater than
\[
  \frac{1}{2}\lambda^{\prime \prime} + \frac{1}{2} = \frac{3}{4} - \frac{x}{2n} > 1 - \frac{x}{n},
\]
because $x > n/2$.

\textbf{Case 2.} $U_3^{\prime}$ is present (but maybe not $U_1^{\prime}$) and $x_2^{\prime \prime} \le x_1^{\prime \prime}$. Then $x_{\min}^{\prime \prime} = x_2^{\prime \prime}$ and $x_{\max}^{\prime \prime} = x_3^{\prime \prime}$, so $\lambda^{\prime \prime} = x_{\min}^{\prime \prime} - x_{\max}^{\prime \prime} = x_2^{\prime \prime} - x_3^{\prime \prime} = 2x - (n + 1)$. Then by Lemma \ref{lem:LB_on_f;t_even_gen}, there is a solution to the problem of assigning the remaining pieces with smallest piece having size greater than
\[
  \frac{1}{2}\lambda^{\prime \prime} + \frac{1}{2} = x - \frac{n}{2},
\]
and observe that
\[
  x \ge \frac{n(n + 2)}{2(n + 1)} \quad \Leftrightarrow \quad \frac{n + 1}{n}x \ge \frac{n + 2}{2} \quad \Leftrightarrow \quad x - \frac{n}{2} \ge 1 - \frac{x}{n}.
\]

\textbf{Case 3.} $U_1^{\prime}$ is present but not $U_3^{\prime}$ and $x_1^{\prime \prime} < x_2^{\prime \prime}$. Then $x_{\min}^{\prime \prime} = x_1^{\prime \prime}$ and $x_{\max}^{\prime \prime} = x_2^{\prime \prime}$, so $\lambda^{\prime \prime} = x_{\min}^{\prime \prime} - x_{\max}^{\prime \prime} = x_1^{\prime \prime} - x_2^{\prime \prime} = [(n + 2) - (2 + 1/n)x] - 1/2$. Then by Lemma \ref{lem:LB_on_f;t_even_gen}, there is a solution to the problem of assigning the remaining pieces with smallest piece having size greater than
\[
  \frac{1}{2}\lambda^{\prime \prime} + \frac{1}{2} = \frac{1}{2}[(n + 2) - (2 + 1/n)x] + \frac{1}{4} > 1 - \frac{x}{n},
\]
because this is the average of two numbers that are largely than $1 - x/n$.

\textbf{Case 4.} $U_1^{\prime}$ is present but not $U_3^{\prime}$ and $x_2^{\prime \prime} \le x_1^{\prime \prime}$. Then $x_{\min}^{\prime \prime} = x_2^{\prime \prime}$ and $x_{\max}^{\prime \prime} = x_1^{\prime \prime}$, so $\lambda^{\prime \prime} = x_{\min}^{\prime \prime} - x_{\max}^{\prime \prime} = x_2^{\prime \prime} - x_1^{\prime \prime} = 1/2 - [(n + 2) - (2 + 1/n)x] = [(2 + 1/n)x - (n + 1)] - 1/2$. Then by Lemma \ref{lem:LB_on_f;t_even_gen}, there is a solution to the problem of assigning the remaining pieces with smallest piece having size greater than
\[
  \frac{1}{2}\lambda^{\prime \prime} + \frac{1}{2} = \frac{1}{2}[(2 + 1/n)x - (n + 1)] + \frac{1}{4} > 1 - \frac{x}{n},
\]
because this is the average of two numbers that are largely than $1 - x/n$.

Finally, we must consider when $x = 4/3$, so $m = 4d$ and $s = 3d$ for some $d \ge 1$. Then we can give $d$ students 4 pieces of size $1 - x/2$ and $2d$ students 2 pieces of size $x/2$ and this solution is optimal.
\hspace{\stretch{1}}$\blacksquare$

This last construction for $x = 4/3$ can be extended to any $n \ge 2$, i.e., when $x = n(n + 2)/[2(n + 1)]$ there is a straightforward optimal solution to the supply-constrained muffin problem that is not feasible for the fully-constrained muffin problem. When $n$ is even, let $d$ be any positive integer; when $n$ is odd, require that $d$ be even. Then:
\begin{enumerate}
\item Cut all the muffins into two pieces, one of size $n/[2(n + 1)]$ and one of size $(n + 2)/[2(n + 1)]$
\item Assign $nd/2$ students $n + 2$ pieces of size $n/[2(n + 1)]$
\item Assign $(n + 2)d/2$ students $n$ pieces of size $(n + 2)/[2(n + 1)]$.
\end{enumerate}
For example, when $x = 15/8$, we have $n = 3$, so the solution with $d = 2$ is:
\begin{enumerate}
\item Cut all 15 muffins into two pieces, one of size $3/8$ and one of size $5/8$
\item Assign $3$ students $5$ pieces of size $3/8$
\item Assign $5$ students $3$ pieces of size $5/8$.
\end{enumerate}

We summarize the conclusions of this section in the following theorem.

\begin{theorem} \label{thm:fully->semi} Any optimal solution for a fully-constrained muffin problem is optimal for the corresponding supply-constrained muffin problem. There exist optimal solutions for the supply-constrained muffin problem that are not feasible for the corresponding fully-constrained muffin problem if and only if $n(n + 2)/[2(n + 1) \le x < (n + 1)/2$ for all $n \ge 2$.
\end{theorem}

\subsection{When does $f(m, s) = 1/3$?} \label{sec:1/3}

We have previously established Conjecture \ref{f>=1/3_conj}: $f(m, s) \ge 1/3$. In this section we investigate when this bound pertains.

\begin{lemma} \label{lem:g>1/3_spec_cases} Let $P$ be a fully-constrained muffin problem. If $P$ is a 0-problem or if $g(P) > \tau$, then $g(P) > 1/3$.
\end{lemma}
\textbf{Proof.} We have $\lambda = 0$, $\gamma = 1/2$, $t = 2$, $u = n + 1$, and $v = n$.

If $P$ is a 0-problem, then
\[
  g(P) = \frac{x}{u} > \frac{n}{2(n + 1)} \ge \frac{1}{3},
\]
because $n \ge 2$.

If $g(P) > \tau$, then
\[
  g(P) > \tau = \frac{x_{\infty}}{u} = \frac{n^2 - 1}{(n + 1)(2n - 1)} = \frac{n - 1}{2n - 1}\ge \frac{1}{3},
\]
because $n \ge 2$.
\hspace{\stretch{1}}$\blacksquare$

\begin{theorem} \label{thm:f>1/3_n>=3} For a muffin problem with $n \ge 3$, $f(m, s) > 1/3$.
\end{theorem}
\textbf{Proof.} It suffices to show that if $P(x)$ is a fully-constrained muffin problem with $\lfloor 2x \rfloor = n \ge 3$, then $g(P) > 1/3$. For such a problem, $\lambda = 0$, $\gamma = 1/2$, $t = 2$, $u = n + 1$, and $v = n$.

If $P$ is a 0-problem, then we can appeal to the previous lemma. Otherwise, $x_b < x < x_{b - 1}$ for some $b > 0$. If $P$ is a 1-problem with $b \ge 3$, then because $n \ge 3$, Lemma \ref{lem:L(b)_a>0} applies so $g(P) > \tau$. Alternatively, if $P$ is an $N$-problem with $N \ge 2$, then $g(P) > \tau$ by Theorem \ref{thm:g}. In either case, we can appeal to the previous lemma.

If $P$ is a 1-problem with either $b = 1$ or $b = 2$, then by Lemma \ref{lem:L(b)_LB}, $g(P) > L(b)$. Now
\[
  L(b) = \gamma - b\frac{x_{b - 1} - \lambda - \gamma v}{[(v - 1)t - v]b + t} = \frac{1}{2} - b\frac{x_{b - 1} - n/2}{(n - 2)b + 2}.
\]
If $b = 2$, recall that $x_1 = n(n + 1)/(2n + 1)$, so
\[
  g(P) > L(2) = \frac{1}{2} - \frac{2n(n + 1) - n(2n + 1)}{2(n - 1)(2n + 1)} = \frac{1}{2} - \frac{n}{2(n - 1)(2n + 1)}
\]
The right-hand side is smallest when $n$ is smallest, so (setting $n = 3$) is at least $11/28 > 1/3$. Alternatively, $b = 1$, and with $x_0 = \gamma u = (n + 1)/2$:
\[
  g(P) > L(1) = \frac{1}{2} - \frac{(n + 1)/2 - n/2}{n} = \frac{1}{2} - \frac{1}{2n} = \frac{n - 1}{2n}.
\]
Again, the right-hand side is smallest when $n$ is smallest, so (setting $n = 3$) is at least 1/3. 
\hspace{\stretch{1}}$\blacksquare$

\begin{theorem} \label{thm:f<1/3_n=2_1problem} For a muffin problem with $n = 2$, $f(m, s) = 1/3$ if and only if
\[
  \frac{3b + 1}{3b} \le x < \frac{3b}{3b - 1},
\]
for all integer $b > 0$. This is precisely when the fully-constrained version of the problem is a 1-problem of type 1. Otherwise $f(m, s) > 1/3$.
\end{theorem}
\textbf{Proof.} Let $P(x)$ be the fully-constrained version of the problem. Then $f(m, s) > 1/3$ if and only if $g(P) > 1/3$ and $f(m, s) = 1/3$ if and only if $g(P) \le 1/3$. Recall that
\[
  \tau = \frac{x_{\infty}}{u} = \frac{n^2 - 1}{(n + 1)(2n - 1)} = \frac{n - 1}{2n - 1} = \frac{1}{3},
\]
because $n = 2$.

If $P$ is a 0-problem, then by Lemma \ref{lem:g>1/3_spec_cases}, $g(P) > 1/3 \Leftrightarrow f(m, s) > 1/3$. Otherwise, $x_b < x < x_{b - 1}$ for some $b > 0$. If $P$ is an $N$-problem with $N \ge 2$, then $g(P) > \tau = 1/3$ by Theorem \ref{thm:g}, so once again $f(m, s) > 1/3$.

Otherwise, $P$ is a 1-problem: its solution is $x^{\prime}/u^{\prime}$. The reduced problem has
\[
\begin{split}
  x^{\prime} &= [(v - 1)b + 1]x_t - bx_v = 1 - b(x - 1)\\
  \lambda^{\prime} &= \lambda + \gamma(v - 1)t - x = 1 - x\\
  \gamma^{\prime} &= \frac{x}{u} = \frac{x}{3}\\
  v^{\prime} &= [(v - 1)(b - 1) + 1]t - (b - 1)v = 2\\
  t^{\prime} &= u = 3\\
  u^{\prime} &= [(v - 1)b + 1]t - bv = 2,
\end{split}
\]
and
\[
  x^{\prime}_{\infty} = \frac{[\lambda^{\prime} + \gamma^{\prime}(v^{\prime} - 1)t^{\prime}]u^{\prime}}{(v^{\prime} - 1)t^{\prime} + u^{\prime} - v^{\prime}} = \frac{2}{3}.
\]
Then $g(P) > 1/3$ when $x^{\prime} > x^{\prime}_{\infty} = 2/3$ and $g(P) \le 1/3$ when $x^{\prime} \le x^{\prime}_{\infty} = 2/3$ (and note that the reduced problem is a 0-problem of type 1 when $x^{\prime} \le x^{\prime}_{\infty} = 2/3$). Now
\[
\begin{split}
  x^{\prime} &\le x^{\prime}_{\infty} = \frac{2}{3}\\
  \Leftrightarrow 1 - b(x - 1) &\le \frac{2}{3}\\
  \Leftrightarrow b(x - 1) &\ge \frac{1}{3}\\
  \Leftrightarrow x - 1 &\ge \frac{1}{3b}\\
  \Leftrightarrow x &\ge 1 + \frac{1}{3b} = \frac{3b + 1}{3b}.
\end{split}
\]
But recall that $x_{b} < x < x_{b - 1}$ and
\[
  x_b = \frac{[\lambda + (v - 1)\gamma t]ub + \gamma tu}{[(v - 1)t + u - v]b + t} = \frac{3b + 3}{3b + 2}, \quad \text{so} \quad x_{b - 1} = \frac{3b}{3b - 1}.
\]
Then $g(P) \le 1/3$ when
\[
  \frac{3b + 1}{3b} \le x < \frac{3b}{3b - 1},
\]
for all integer $b > 0$; otherwise $g(P) > 1/3$.
\hspace{\stretch{1}}$\blacksquare$

\subsection{Relationships between muffin problems} \label{sec:relationships}

Consider the fully-constrained muffin problem family of order $n=2$. Here $t = v = 2$, $u = 3$, $\lambda=0$, and $\gamma=1/2$. Consider $P(x)$ in the family that is an $N$-problem with $N \ge 2$. Then $x_i < x < x_{i-1}$ for some $i > 0$ and $P(x)$ can be (at least) doubly reduced. The first reduced problem $P^{\prime}(x^{\prime})$, where $x_j^{\prime} < x^{\prime} < x_{j - 1}^{\prime}$ for some $j > 0$, will have $t^{\prime} = 3$ and $u^{\prime} = v^{\prime} = 2$, regardless of the value of $i$. Then the second reduced problem $P^{\prime \prime}$ will have $t^{\prime \prime} = 2$, $u^{\prime \prime} = j + 3$, and $v^{\prime \prime} = j + 2$. Then the standardized version of this second reduced problem will be a fully-constrained muffin problem of order $j + 2$. Because $j \ge 1$, this order will be at least 3. In other words, any muffin problem of order $n=2$ can be related to a muffin problem of order greater than 2. In this section we make this relationship concrete. 

In fact, we will show that for any muffin problem of order greater than 2, for every value of $b \ge 1$ there is a muffin problem $P(x)$ of order 2 with $x_b < x < x_{b - 1}$ that is related to it. This implies that there are relationships between muffin problems of order 2: we will show an explicit method for constructing a solution to a muffin problem $P(x^{b-1})$ of order 2 with $x_{b - 1} < x^{b - 1} < x_{b - 2}$ given a solution to the related muffin problem $P(x^b)$ of order 2 with $x_b < x^b < x_{b - 1}$, for all $b \ge 2$, and vice versa.

\subsubsection{Relating a muffin problem of order $2$ to a muffin problem of order $> 2$} \label{sec:relationship_2_to_>2}

The domain of $x$ for the fully-constrained muffin problem family of order $n=2$ is $(\lambda + \gamma v, \gamma u] = (1, 3/2]$. Further,
\[
  x_b = \frac{[\lambda + \gamma(v - 1)t]ub + \gamma tu}{[(v - 1)t + u - v]b + t} = \frac{3b + 3}{3b + 2},
\]
so $x_0 = 3/2$, $x_1 = 6/5$, and $x_{\infty} = 1$. 

Consider a problem $P(x)$ in the family with
\[
  x_b = \frac{3b + 3}{3b + 2} < x \le \frac{3b}{3b - 1} = x_{b - 1}.
\]
If we write $d = m - s$, then
\[
  x = \frac{m}{s} = \frac{s + d}{s} = \frac{s/d + 1}{s/d},
\]
so
\[
  3b - 1 \le \frac{s}{d} < 3b + 2 \quad \Leftrightarrow \quad 3db - d \le s < 3db + 2d.
\]
So write $s = 3db + a - d$, where $0 \le a < 3d$. Then $m = 3db + a$. Any problem in the family can be written in the form $(m, s) = (3db + a, 3db + a - d)$ for some $d \ge 1$, $b \ge 1$, and $0 \le a < 3d$. 

The problem is a 0-problem when $a = 0$, so that $s = 3db - d$, $m = 3db$, and $x = x_{b - 1}$, for $b = 2, 3, \ldots$ (recall that $x = x_0 = 3/2$ is a special case with $f(m,s) = 1/2$). Then
\[
  f(m, s) = \frac{x}{3} = \frac{b}{3b - 1} = \frac{db}{3db + a - d}.
\]

When $x_b < x < x_{b - 1}$ or equivalently, $a > 0$, the problem is not a 0-problem. The reduced problem has:
\[
\begin{split}
  x^{\prime} &= [(v - 1)b + 1]x_t - bx_v = 1 - b(x - 1)\\
  x_t^{\prime} &= x_u = x\\
  x_v^{\prime} &= [(v - 1)(b - 1) + 1]x_t - (b - 1)x_v = 1 - (b - 1)(x - 1)\\
  \lambda^{\prime} &= \lambda + \gamma(v - 1)t - x = 1 - x\\
  \gamma^{\prime} &= \frac{x}{u} = \frac{x}{3}\\
  u^{\prime} &= [(v - 1)b + 1]t - bv = 2\\
  t^{\prime} &= u = 3\\
  v^{\prime} &= [(v - 1)(b - 1) + 1]t - (b - 1)v = 2.
\end{split}
\]

Write $g(m, s) = g(P)$ for the optimal value of the fully-constrained muffin problem. Recall that the problem is a 1-problem of type 1 with $f(m, s) = 1/3$ and $g(m, s) = [1 - b(x - 1)]/2$ when $x^{\prime} = 1 - b(x - 1) \le x_{\infty}^{\prime} = 2/3 \Leftrightarrow x \ge 1 + 1/(3b) = (3b + 1)/(3b)$, i.e., when
\[
  \frac{3b + 1}{3b} \le x < \frac{3b}{3b - 1} = x_{b - 1},
\]
or, equivalently, when $0 < a \le d$. Then
\[
  g(m, s) = \frac{s - db}{2s} = \frac{db + (a - d)/2}{3db + a - d}.
\]

When
\[
  x_b = \frac{3b + 3}{3b + 2} < x < \frac{3b + 1}{3b},
\]
i.e., when $d < a < 3d$, the problem is an $N$-problem with $N \ge 1$. It is only a 1-problem of type 2 when $x^{\prime} = x_j^{\prime}$ for some $j \ge 1$, where:
\[
  x_j^{\prime} = \frac{[\lambda^{\prime} + \gamma^{\prime}(v^{\prime} - 1)t^{\prime}]u^{\prime}j + \gamma^{\prime} t^{\prime}u^{\prime}}{[(v^{\prime} - 1)t^{\prime} + u^{\prime} - v^{\prime}]j + t^{\prime}} = \frac{2(j + x)}{3(j + 1)}.
\]
Now
\[
\begin{split}
  x_j^{\prime} &< x^{\prime}\\
  \Leftrightarrow \frac{2(j + x)}{3(j + 1)} &< 1 - b(x - 1)\\
  \Leftrightarrow 2(j + x) &< 3(j + 1)[1 - b(x - 1)] = 3(j + 1)(b + 1) - 3(j + 1)bx\\
  \Leftrightarrow [3(j + 1)b + 2]x &< 3(j + 1)b + j + 3\\
  \Leftrightarrow x &< \frac{3(j + 1)b + j + 3}{3(j + 1)b + 2},
\end{split}
\]
and note that when $j = 0$, the right-hand side is $(3b + 3)/(3b + 2)$, and when $j \rightarrow \infty$, the right-hand side tends to $(3b + 1)/(3b)$. For what values of $a$ and $d$ is the problem a 1-problem of type 2:
\[
\begin{split}
  \frac{3db + a}{3db + a - d} &= \frac{3(j + 1)b + j + 3}{3(j + 1)b + 2}\\
  \Leftrightarrow [3(j + 1)b + 2](3db + a) &= [3(j + 1)b + j + 3](3db + a - d)\\
  \Leftrightarrow 3(j + 1)db + 2d &= 3(j + 1)db + (j + 1)(a - d)\\
  \Leftrightarrow 2d &= (j + 1)(a - d)\\
  \Leftrightarrow a &= \frac{j + 3}{j + 1} d.
\end{split}
\]
In this case
\[
\begin{split}
  f(m, s) = \frac{x_j^{\prime}}{u^{\prime}} = \frac{j + x}{3(j + 1)} &= \frac{1}{s} \frac{3db(j + 1) + (a - d)(j + 1) + d}{3(j + 1)}\\
  &= \frac{db + (a - d)/3 + d/[3(j + 1)]}{3db + a - d}\\
  &= \frac{db + (a - d)/2}{3db + a - d} > \frac{1}{3}.
\end{split}
\]
Otherwise, the problem is an $N$-problem with $N \ge 2$ and there is some $j \ge 1$ for which: 
\[
  \frac{3(j + 1)b + j + 3}{3(j + 1)b + 2} < x < \frac{3jb + j + 2}{3jb + 2}.
\]
We can further reduce the problem to obtain:
\[
\begin{split}
  x^{\prime \prime} &= [(v^{\prime} - 1)j + 1]x_t^{\prime} - jx_v^{\prime} = (j + 1)x - j + j(b - 1)(x - 1) = x + jb(x - 1)\\
  x_t^{\prime \prime} &= x_u^{\prime} = x^{\prime}\\
  x_v^{\prime \prime} &= [(v^{\prime} - 1)(j - 1) + 1]x_t^{\prime} - (j - 1)x_v^{\prime} = x + (j - 1)b(x - 1)\\
  \lambda^{\prime \prime} &= \lambda^{\prime} + \gamma^{\prime}(v^{\prime} - 1)t^{\prime} - x^{\prime} = b(x - 1)\\
  \gamma^{\prime \prime} &= \frac{x^{\prime}}{u^{\prime}} = \frac{1}{2}[1 - b(x - 1)]\\
  u^{\prime \prime} &= [(v^{\prime} - 1)j + 1]t^{\prime} - jv^{\prime} = 3(j + 1) - 2j = j + 3\\
  t^{\prime \prime} &= u^{\prime} = 2\\
  v^{\prime \prime} &= [(v^{\prime} - 1)(j - 1) + 1]t^{\prime} - (j - 1)v^{\prime} = 3j - 2(j - 1) = j + 2.
\end{split}
\]
This doubly-reduced problem has $t^{\prime \prime} = 2$ and $u^{\prime \prime} = v^{\prime \prime} + 1$, so transforming it to standard form results in a muffin problem. The equivalence function is
\[
  h(y) = \frac{1}{2} \frac{(u^{\prime \prime} - v^{\prime \prime})y - \lambda^{\prime \prime}}{(u^{\prime \prime} - v^{\prime \prime})\gamma^{\prime \prime} - \lambda^{\prime \prime}} = \frac{1}{2} \frac{y - b(x - 1)}{[1 - b(x - 1)]/2 - b(x - 1)} = \frac{y - b(x - 1)}{1 - 3b(x - 1)},
\]
and the inverse transformation is
\[
  h^{-1}(z) = 2\gamma^{\prime \prime} z + (1 - 2z)\frac{\lambda^{\prime \prime}}{u^{\prime \prime} - v^{\prime \prime}} = [1 - b(x - 1)]z + (1 - 2z)b(x - 1) = b(x - 1) + [1 - 3b(x - 1)]z.
\]
When $y_1 + \cdots + y_{j + 3} = x^{\prime \prime}$, we find
\[
  x_{\text{std}}^{\prime \prime} = h(y_1) + \cdots + h(y_{j+3}) = \frac{x^{\prime \prime} - (j + 3)b(x - 1)}{1 - 3b(x - 1)} = \frac{x - 3b(x - 1)}{1 - 3b(x - 1)} = \frac{1 - (3b - 1)(x - 1)}{1 - 3b(x - 1)}.
\]
Observe that
\[
  x_{\text{std}}^{\prime \prime}(x_j^{\prime}) = \frac{3(j + 1)b + 2 - (3b - 1)(j + 1)}{3(j + 1) + 2 - 3b(j + 1)}
  = \frac{j + 3}{2}.
\]
This is independent of $b$. Then for every $b \ge 1$, the set of muffin problems with $n = 2$ and $x_b < x < (3b + 1)/3b$ is mapped to the full set of muffin problems with $n > 2$.

Finally, the doubly-reduced problem has
\[
  x_{\text{std}}^{\prime \prime} = \frac{s - (3b - 1)(m - s)}{s - 3b(m - s)} = \frac{a}{a - d}.
\]
Then
\[
\begin{split}
  f(m, s) = f(3db + a, 3db + a - d) &= h^{-1}(f(a, a - d))\\
  &= b(x - 1) + [1 - 3b(x - 1)]f(a, a - d)\\
  &= \frac{1}{s}[db + (s - 3db)f(a, a - d)]\\
  &= \frac{1}{s}[db + (a - d)f(a, a - d)]\\
  &= \frac{db + (a - d)f(a, a - d)}{3db + a - d}.
\end{split}
\]
Note that this is greater than 1/3 because $f(a, a - d) > 1/3$ and $a > d$.

We summarize the results of this section with the following theorem.

\begin{theorem} \label{thm:n=2->n>2} For a fully-constrained muffin problem $(m, s) = (3db + a, 3db + a - d)$ with $d \ge 1$, $0 \le a < 3d$, and $b \ge 1$, we have
\[
  f(m, s) = \frac{db + X(a, d)}{3db + a - d},
\]
where
\[
  X(a, d) = 
  \begin{cases}
  0 & \text{if } a = 0\\
  (a - d)/2 & \text{if } 0 < a \le d \text{ or } 2d \text{ is divisible by } a - d\\
  (a - d)f(a, a - d) & \text{otherwise.}\\
  \end{cases}
\]
\end{theorem}

\subsubsection{Relationships between muffin problems of order $2$} \label{sec:relationship_2_to_2}

Consider a fully-constrained muffin problem $P(x)$ with $n = 2$ and $x_b < x < x_{b - 1}$ for some $b > 0$. For this problem, let $m \equiv m^b = 3db + a$ and $s \equiv s^b = 3db + a - d$. The previous section showed that
\[
  f(m^b, s^b) = \frac{db + X(a, d)}{s^b}.
\]
Then it is immediate that for $b \ge 2$:
\[
  s^bf(m^b, s^b) = s^{b - 1}f(m^{b - 1}, s^{b - 1}) + d.
\]
This implies a relationship between the muffin problems $P = (m^b, s^b)$ and $P^{\prime \prime} = (m^{b - 1}, s^{b - 1})$. In the following subsections we establish this relationship directly by showing that given a solution to $P$, we can construct a solution to $P^{\prime \prime}$, and vice versa.

Starting with a solution to $P$, we reduce the solution to a solution of a reduced problem $P^{\prime}$ (this is a different reduction than the one presented earlier in the paper). Standardizing the solution gives a solution to $P^{\prime \prime}$. In the other direction, we first apply the inverse of the equivalence function to the solution to $P^{\prime \prime}$ to obtain a solution to $P^{\prime}$. We then expand this to obtain a solution to $P$.

\subsubsection{The Alternative Reduction} \label{sec:alternative_reduction}

We begin with a lemma that allows us to restrict attention to certain problem solutions.

\begin{lemma} \label{lem:no_UU} For a fully-constrained muffin problem $P(x) = (T; U, V)$ with $x_b < x < x_{b - 1}$ for some $b \ge 2$ (so that $x < x_1 = 6/5$), there exists an optimal solution in which no row of $T$ contains two $U$-elements.
\end{lemma}
\textbf{Proof.} One optimal solution is given by Algorithm \ref{alg:recursive}. In this solution $(T, V)$ is divided into maximal, inseparable pairs, which are either $b$-pairs or $(b - 1)$-pairs. A $k$-pair contains $[(v - 1)k + 1]t - kv = 2$ elements from $U$. If $k \ge 1$, then these two elements must be in different rows of $T$, otherwise the pair would be separable. But by assumption, $b - 1 \ge 1$, so the conclusion follows.
\hspace{\stretch{1}}$\blacksquare$

%\begin{lemma} \label{lem:no_UU} For a fully-constrained muffin problem with $x < 6/5$, in an optimal solution every element of $U$ has size less than 1/2.
%\end{lemma}
%\textbf{Proof.} Because $x < 6/5$ we have $n = 2$ and $b \ge 2$. Let $Y = f(m, s, 2, 2) = g(m, s)$. By considering a row of $U$, the largest value that an element of $U$ can take is if the other two elements in the row take size $Y$: this largest value is then $x - 2Y$. If the problem is a 0-problem, then $Y = x/3$, so $x - 2Y = x/3 < 1/2$. Otherwise, by Theorem \ref{thm:f<1/3_n=2_1problem}, there are two cases to consider.
%
%\textbf{Case 1.} $x \le 7/6$ and $Y > 1/3$. Then this largest value is less than $7/6 - 2/3 = 1/2$.
%
%\textbf{Case 2.} The problem is a 1-problem with $Y = g(m, s) = [1 - b(x - 1)]/2$ and $x < x_{b - 1} = 3b/(3b - 1)$. Then
%\[
%  x - 2Y = x - 1 + b(x - 1) = (b + 1)(x - 1) < (b + 1)(x_{b - 1} - 1) = \frac{b + 1}{3b - 1} \le \frac{1}{2},
%\]
%because $b \ge 3$. What about $b = 2$??????
%\hspace{\stretch{1}}$\blacksquare$

The $P = (T; U, V)$ problem has $b \ge 2$ (so $2s_v/r_v > 1 \Leftrightarrow n_v > n_u$). By Lemma \ref{lem:no_UU} there is an optimal solution in which no row of $T$ contains two $U$-elements, so restrict attention to solutions in which each $U$-element is inserted into a separate row of $T$. Then $T$ consists of two types of rows: those with one element from $U$ and one from $V$ , and those with both elements from $V$. Write $T_1$ and $T_2$ for the submatrices of $T$ consisting of these two types of rows.

Given such a solution to $P$ in which no row of $T$ contains two $U$-elements, we construct a solution to an alternative reduced problem $P^{\prime} = (T^{\prime}, U^{\prime}, V^{\prime})$ as follows:
\begin{itemize}
\item $T^{\prime}$: for every row $[z, x - z]$ in $V$, there is a row $[1 - z, 1 + z - x]$ in $T^{\prime}$;
\item $U^{\prime} = U$;
\item $V^{\prime}$: for every row $[z, 1 - z]$ in $T_2$, there is a row $[1 - z, z]$ in $V^{\prime}$.
\end{itemize}
It is easy to see that $n_t^{\prime} = n_u^{\prime} + n_v^{\prime}$. Also, if $z \in V^{\prime}$, then $1 - z \in T_2$, so $1 - z \in V$, so $z \in T^{\prime}$; and if $z \in U^{\prime} = U$, then $1 - z \in T_1$ and $V$, so $z \in T^{\prime}$. Further, confirm that
\[
  \frac{x_t^{\prime}}{t^{\prime}} = 1 - \frac{x_v}{2} < \frac{1}{2} = \frac{x_v^{\prime}}{v^{\prime}},
\]
because $1 < x_v = x$; and
\[
  \frac{x_t^{\prime}}{t^{\prime}} = 1 - \frac{x_v}{2} = 1 - \frac{x}{2} > \frac{x}{3} = \frac{x_u^{\prime}}{u^{\prime}},
\]
because $x < 6/5$. Then $P^{\prime}$ is a 3M-DAP. 

If the smallest element in the solution to $P$ appears in $U$, then because every element of $P^{\prime}$ is an element of $P$ and $U^{\prime} = U$, the smallest element in the constructed solution to $P^{\prime}$ appears in $U^{\prime}$. If we start with an optimal solution to $P$, then the value $g(P)$ is an element of $U$, so $g(P^{\prime}) = g(P)$ and this value is an element of $U^{\prime}$.

Now $P$ has $s_t = m^b$, $s_u = 2d$, and $s_v = s^b - 2d = 3d(b - 1) + a = m^{b - 1}$. Then the reduced problem $P^{\prime}$ has $s_t^{\prime} = s_v = m^{b - 1}$, $s_u^{\prime} = s_u = 2d$, and $s_v^{\prime} = s_t - n_u = m^b - 6d$, so $s_u^{\prime} + s_v^{\prime} = 2d + m^b - 6d = m^b - 4d = 3db + a - 4d = 3d(b - 1) + a - d = s^{b - 1}$, as required.

The reduced problem is not in standard form. The function that transforms it to standard form is:
\[
  h(y) = \frac{1}{2} \frac{(u^{\prime} - v^{\prime})y - (x_u^{\prime} - x_v^{\prime})}{(u^{\prime} - v^{\prime})x_t^{\prime}/t^{\prime} - (x_u^{\prime} - x_v^{\prime})} = \frac{1}{2} \frac{y - (x - 1)}{1 - x/2 - (x - 1)} = \frac{y - (x - 1)}{4 - 3x} = \frac{s^by - d}{4s^b - 3m^b} = \frac{s^by - d}{s^{b - 1}}.
\]
Confirm that if $y_1 + y_2 + y_3 = x_u^{\prime} = x = m^b/s^b$, then
\[
  h(y_1) + h(y_2) + h(y_3) = \frac{m^b - 3d}{s^{b - 1}} = \frac{m^{b - 1}}{s^{b - 1}}.
\]
Then the reduced standardized problem is indeed the fully-constrained muffin problem $P^{\prime \prime}$ and the constructed solution is a feasible solution therefor.

Write $Y^b$ for the smallest value in the solution to the $P = (m^b, s^b)$ problem and $Y^{b - 1}$ as the smallest value in the constructed solution to the reduced, standardized problem, i.e., the $P^{\prime \prime} = (m^{b - 1}, s^{b - 1})$ problem. Then
\[
  Y^{b - 1} = h(Y^b) = \frac{s^bY^b - d}{s^{b - 1}} \quad \Leftrightarrow \quad Y^b = \frac{s^{b - 1}Y^{b - 1} + d}{s^b}.
\]
If the original solution to $P$ is optimal, then we obtain
\[
  g(m^b, s^b) = Y^b = \frac{s^{b - 1}Y^{b - 1} + d}{s^b} \le \frac{s^{b - 1}g(m^{b - 1}, s^{b - 1}) + d}{s^b}.
\]

\subsubsection{An Equivalent Expansion} \label{sec:equivalent_expansion}

Here we go in the other direction: again with $b \ge 2$, given a solution to the $P^{\prime \prime} = (m^{b - 1}, s^{b - 1})$ problem, construct a solution to the $P = (m^b, s^b)$ problem.

Start by transforming the solution to $P^{\prime \prime}$ using the inverse transformation from above:
\[
\begin{split}
  h^{-1}(z) = 2z\frac{x_t^{\prime}}{t^{\prime}} + (1 - 2z)\frac{x_u^{\prime} - x_v^{\prime}}{u^{\prime} - v^{\prime}} &= (2 - x)z + (x - 1)(1 - 2z)\\
  &= (4 - 3x)z + (x - 1)\\
  &= \frac{1}{s^b}[(3d(b - 1) + a - d)z + d]\\
  &= \frac{1}{s^b}(s^{b - 1}z + d).
\end{split}
\]
Applying this inverse transformation to the solution of $P^{\prime \prime}$ results in a solution to $P^{\prime}$. The rowsums of $T^{\prime}$ will be
\[
  \frac{s^{b - 1} + 2d}{s^b} = \frac{s^b - 3d + 2d}{s^b} = \frac{s^b - d}{s^b} = \frac{2s^b - m^b}{s^b} = 2 - x,
\]
while the rowsums for $V^{\prime}$ will be
\[
  \frac{m^{b - 1} + 2d}{s^b} = \frac{m^b - 3d + 2d}{s^b} = \frac{m^b - d}{s^b} = 1,
\]
and the rowsums for $U^{\prime}$ will be
\[
  \frac{m^{b - 1} + 3d}{s^b} = \frac{m^b - 3d + 3d}{s^b} = x.
\]
We now construct a solution to $P$ from the solution to $P^{\prime}$, as follows:
\begin{itemize}
\item $T^1$:  for each element of $U^{\prime}$ (with value $y$), there is a row $[y, 1 - y]$ in $T^1$;
\item $T^2$: for every row $[y, 1 - y]$ in $V^{\prime}$, there is a row $[1 - y, y]$ in $T^2$;
\item $U = U^{\prime}$;
\item $V$: for every row $[y, 2 - x - y]$ in $T^{\prime}$, there is a row $[1 - y, x + y - 1]$ in $V$.
\end{itemize}
Then $n_t = 2n_u^{\prime} + n_v^{\prime} = n_u + n_u^{\prime} + n_v^{\prime} = n_u + n_t^{\prime} = n_u + n_v$. Also, if $y \in V$, then $1 - y \in T^{\prime}$, so either $1 - y \in U^{\prime}$, in which case $y \in T$, or $1 - y \in V^{\prime}$, so $y \in T$; and if $y \in U$, then $y \in U^{\prime}$, so $y \in T$. Further, confirm that
\[
  \frac{x_u}{u} = \frac{x}{3} < \frac{x_t}{t} = \frac{1}{2} < \frac{x}{2} = \frac{x_v}{v}.
\]
Then $P$ is indeed the fully-constrained muffin problem $(m^b, s^b)$, and the constructed solution is a feasible solution for $P$.

We cannot immediately claim that the smallest element in the solution to $P^{\prime \prime}$ will match to the smallest element in the constructed solution to $P$. This is because the expansion step might introduce a smaller element. The following lemma shows that, provided all elements of $U^{\prime \prime}$ are nonnegative, this cannot happen.

\begin{lemma} \label{lem:U'<1-Y} In the transformation and expansion of a fully-constrained muffin problem with $n=2$, as described above, suppose that all elements of $U^{\prime \prime}$ are nonnegative. Then any element added to the solution is at least as large as the smallest element of $U^{\prime}$.
\end{lemma}
\textbf{Proof.} Let $Y \ge 0$ be the smallest element of $U^{\prime \prime}$. Then the smallest element of $U^{\prime}$ is:
\[
  Y^{\prime} = h^{-1}(Y) = \frac{s^{b - 1}Y + d}{s^b} \ge \frac{d}{s^b} = x - 1.
\]
The additional elements added to the solution in the expansion step are $\{1 - y : y \in U^{\prime}\}$. For any $y \in U^{\prime}$, consider the row of $U^{\prime}$ containing $y$. Then:
\[
  y \le x - 2Y^{\prime} \le 1 - Y^{\prime}
\]
so $1 - y \ge Y^{\prime}$, as claimed.
\hspace{\stretch{1}}$\blacksquare$

Assume that all elements of $U^{\prime \prime}$ are nonnegative. Write $Z^{b - 1}$ for the smallest value in the solution to the $P^{\prime \prime} = (m^{b - 1}, s^{b - 1})$ problem and $Z^b$ as the smallest value in the constructed solution to the expanded problem $P = (m^b, s^b)$ problem. Then
\[
  Z^b = h^{-1}(Z^{b - 1}) = \frac{s^{b - 1}Z^{b - 1} + d}{s^b} \quad \Leftrightarrow \quad Z^{b - 1} = \frac{s^bZ^b - d}{s^{b - 1}}.
\]
If the original solution to $P^{\prime \prime}$ is optimal, then $Z^{b - 1} = g(m^{b - 1}, s^{b - 1}) > 0$, so Lemma \ref{lem:U'<1-Y} applies, and we obtain
\[
  g(m^b, s^b) \ge Z^b = \frac{s^{b - 1}Z^{b - 1} + d}{s^b} = \frac{s^{b - 1}g(m^{b - 1}, s^{b - 1}) + d}{s^b}.
\]
Combining this with the conclusion from the alternative reduction leads us to conclude that
\[
  g(m^b, s^b) = \frac{s^{b - 1}g(m^{b - 1}, s^{b - 1}) + d}{s^b},
\]
as we have previously established.

\subsection{Proof of Conjecture \ref{conj:f(km,ks)=f(m,s)}} \label{conj1_proof}

We begin by showing the equivalent result for 3M-DAPs. The proof of the conjecture will immediately follow. For a 3M-DAP $P$ and $k > 1$, let $kP$ be the 3M-DAP in which the three matrices have the same number of columns and the same rowsums as the given problem but each matrix has $k$ times as many rows.

\begin{theorem} \label{thm:g(kP)=g(P)} Let $P$ be a 3M-DAP and $k > 1$. Then $g(kP)=g(P)$.
\end{theorem}
\textbf{Proof.} The proof is by induction on $N$, the problem type of the $P$ problem.

Suppose $P$ is a 0-problem. Then either $n_u \le (t - 2)s_t + (v - 2)s_v \Leftrightarrow kn_u \le (t - 2)ks_t + (v - 2)ks_v$, or $(t - 2)s_t < n_u \le (t - 1)s_t$ and $(v - 2)s_v < r_t \le (v - 1)s_v$ and $2s_v$ is divisible by $r_v$ which implies that $(t - 2)ks_t < kn_u \le (t - 1)ks_t$ and $(v - 2)ks_v < kr_t \le (v - 1)ks_v$ and $2ks_v$ is divisible by $kr_v$. In either case, $kP$ is also a 0-problem, and we have $g(kP) = x_u/u = g(P)$.

Suppose $P$ is an $N$-problem with $N \ge 1$. Then Algorithm \ref{alg:recursive} reduces this to a $P^{\prime}$ problem. Now $kP$ cannot be a 0-problem because otherwise $P$ would also be a 0-problem. Then Algorithm \ref{alg:recursive} reduces this to the $kP^{\prime}$ problem. Applying the inductive hypothesis gives
\[
  g(kP) = g(kP^{\prime}) = g(P^{\prime}) = g(P).
\]
\hspace{\stretch{1}}$\blacksquare$

\begin{theorem} \label{thm:f(km,ks)=f(m,s)} For all $k \ge 2, m \ge s \ge 1$, $f(km, ks) = f(m, s)$.
\end{theorem}
\textbf{Proof.} It is an immediate corollary of the previous theorem that $g(km, ks) = g(m, s)$. By Theorems \ref{thm:f>1/3_n>=3} and \ref{thm:f<1/3_n=2_1problem} whether $f(m, s)$ is greater than or is equal to 1/3 depends solely on the value of $m/s$. The conclusion follows.
\hspace{\stretch{1}}$\blacksquare$

%\newpage

\section*{References}

\begin{description}[noitemsep]

\item
Antonick, G. (editor), 2013. The New York Times Numberplay Online Blog, August 19, 2013. \url{wordplay.blogs.nytimes.com/2013/08/19/cake}

\item
Blachman, N. (editor), 2016. Julia Robinson Mathematics Festival: A Sample of Mathematical Puzzles.

\item
Cormen, T. H.,  C. E. Leiserson, R. L. Rivest, and C. Stein, 2001. \textit{Introduction to Algorithms}. Second Edition. MIT Press and McGraw-Hill. ISBN 0-262-03293-7. 

\item
Cui, G.,  J. P. Dickerson, N. Durvasula, W. Gasarch, E. Metz, J. Prinz, N. Raman, D. Smolyak, and S. H. Yoo, 2018. A Muffin-Theorem Generator. \textit{Proceedings of the 9th International Conference on Fun with Algorithms (FUN 2018)}, \textbf{100}, 15:1--15:19. 
	\url{doi.org/10.4230/LIPIcs.FUN.2018.15}

\item
Dai, J., D. Li, Yi Liao, and F. Zhu, 1996. An extension to the single bottleneck transportation problem. \textit{International Journal of Systems Science}, \textbf{27(6)}, 577--581. \url{doi.org/10.1080/00207729608929252}

\item
Garfinkel, R. S. and M. R. Rao, 1971. The bottleneck transportation problem. \textit{Naval Research Logistics Quarterly}, John Wiley \& Sons, \textbf{18(4)}, 465--472, December. \url{doi.org/10.1002/nav.3800180404}

\item
Gasarch, W. 2019. The Muffin Website.
\url{www.cs.umd.edu/users/gasarch/MUFFINS/muffins.html}

\item
Gasarch, W., E. Metz, J. Prinz, and D. Smolyak, 2020. \textit{Mathematical Muffin Morsels: Nobody Wants a Small Piece}. Problem Solving in Mathematics and Beyond: Volume 16. World Scientific.

%\item
%Andelman, N., and Y. Mansour, 2004. Auctions with budget constraints. Proceedings of 9th Scandinavian Workshop on Algorithm Theory SWAT, pages 26--38.
%
%\item
%Buchbinder, N., and J. Naor, 2007. The Design of Competitive Online Algorithms via a Primal-Dual Approach. \textit{Foundations and Trends in Theoretical Computer Science} \textbf{3}, Nos. 2-3, 93-263. 
%
%\item
%Chakrabarty, D., and G. Goel, 2008. On the approximability of budgeted allocations and improved lower bounds for submodular welfare maximization and GAP. Proc. of IEEE FOCS, 687--696.
%

\end{description}

\newpage

\appendix

\section*{Appendix A: Further results}
\stepcounter{section} \label{AppA}

\subsection{The Greedy Algorithm for 3M-DAP with $t=u=v=2$ is Optimal} \label{App_Greedy_Opt}

Assume $s_v > 0$, otherwise the problem is trivial. Let $g = \text{gcd}(s_u, s_v)$ and write $s_u = ga$ and $s_v = gb$. Then $a$ and $b$ are coprime. Further, $s \equiv s_t = s_u + s_v = g(a + b)$ is also divisible by $g$. We have
\begin{equation} \label{abc}
  (a + b)x_t = ax_u + bx_v \Leftrightarrow a(x_t - x_u) = b(x_v - x_t) \Leftrightarrow x_v - x_u = \frac{a + b}{b}(x_ t - x_u),
\end{equation}
We can map the problem to an equivalent problem $P^{\prime}$ using the equivalence function
\[
  h(y) = \frac{b}{x_t - x_u}\left(y - \frac{x_t}{2}\right) = \frac{a}{x_v - x_t}\left(y - \frac{x_t}{2}\right).
\]
Then $P^{\prime}$ has $x_t^{\prime} = 0$, $x_u^{\prime} = -b$, and $x_v^{\prime} = a$.

\begin{lemma} \label{P'_opt_LB} An upper bound on the solution to problem $P^{\prime}$ is given by $-(a + b - 1)/2$.
\end{lemma}
\textbf{Proof.} Given a solution to $P^{\prime}$, select the piece whose size is closest to 0. Call this $y_1^{\prime}$. This piece must show up in either $U^{\prime}$ or $V^{\prime}$. If it appears in $U^{\prime}$, the other piece in that row of $U^{\prime}$ is $z_1^{\prime} = x_u^{\prime} - y_1^{\prime} = - b - y_1^{\prime}$. This piece belongs to some muffin whose other piece has size $y_2^{\prime} = x_t^{\prime} - z_1^{\prime} = y_1^{\prime} + b$. Alternatively, if $y_1^{\prime}$ appears in $V^{\prime}$, then we find $z_1^{\prime} = x_v^{\prime} - y_1^{\prime} = a - y_1^{\prime}$ and $y_2^{\prime} = x_t^{\prime} - z_1^{\prime} = y_1^{\prime} - a$. Then $y_2^{\prime}$ differs from $y_1^{\prime}$ by either $b$ or $-a$.

Proceeding in this manner we obtain $y_{i + 1}^{\prime} = - z_i^{\prime} = y_1^{\prime} + qb - pa$ for integer $p$ and $q$ with $p + q = i$. The process terminates when $z_i^{\prime} = -y_1^{\prime}$, i.e., when $qb - pa = 0$. Because $a$ and $b$ are coprime, this requires $p = db, q = da$ for some integer $d \ge 1$. We can then repeat this process until all muffins are assigned.

Suppose there are $1 \le i < j \le d(a + b)$ with $y_i^{\prime} = y_j^{\prime}$. Then there must exist $m, n > 0$ with $m + n = j - i$ and $nb - ma = 0$, or $ma = nb$. Because $a$ and $b$ are coprime, this requires $m = kb, n = ka$ for some integer $k \ge 1$. Then $m + n = k(a + b) \ge a + b$. Then the first $a + b$ values of the $\{y_i^{\prime}\}$ must be distinct. Now
\begin{align*}
  \min \{y_i^{\prime}, z_i^{\prime} : 1 \le i \le s_t\} &\le \min \{y_i^{\prime}, z_i^{\prime} : 1 \le i \le a + b\}\\
  &= \min \{y_i^{\prime}, -y_i^{\prime} : 1 \le i \le a + b\}\\
  &= \min \{-|y_i^{\prime}| : 1 \le i \le a + b\}\\
  &= \min \{-y_{\max}^{\prime}, y_{\min}^{\prime}\},
\end{align*}
where $y_{\max}^{\prime} = \max \{y_i^{\prime}: 1 \le i \le a + b\}$ and $y_{\min}^{\prime} = \min \{y_i^{\prime}: 1 \le i \le a + b\}$. Because all the $\{y_i^{\prime}\}$ differ from one another by integer values, we must have
\[
  y_{\max}^{\prime} - y_{\min}^{\prime} \ge a + b - 1.
\]
Then
\[
  \min\{-y_{max}^{\prime}, y_{\min}^{\prime}\} \le -\frac{a + b - 1}{2}.
\]
\hspace{\stretch{1}}$\blacksquare$

\begin{theorem} \label{greedy_opt_P'} The greedy algorithm applied to problem $P^{\prime}$ produces an optimal solution. The size of the smallest piece is $-(a + b - 1)/2$.
\end{theorem}
\textbf{Proof.} Part 1 of the greedy algorithm starts by dividing the first muffin $[y_1^{\prime}, z_s^{\prime}] = [0, 0]$. The state of the algorithm at any point can be represented by $(p, q)$, where $p$ is the number of students requiring total $-b$ ($U$-students) remaining to be assigned and $q$ is the number of students requiring total $a$ ($V$-students) remaining to be assigned. The value in this state is
\[
  w_{p + q} \equiv w(p, q) \doteq y_{s+1-p-q}^{\prime} = qa - pb.
\]
The initial state is $(ga, gb)$ and the initial value is $w_s = w_{ga + gb} = w(ga, gb) = y_1^{\prime} = 0$.

In state $(p, q)$ the greedy algorithm assigns whichever type of student results in the smallest value of $|w_{p + q - 1}|$, i.e., assigns a $U$-student if $|w(p - 1, q)| < |w(p, q - 1)|$, else assigns a $V$-student. 

When $q$ is sufficiently large that $w(p, q) \ge 0$, then $|w(p, q - 1)| \le \max\{w(p, q), a\}$, while $w(p - 1, q) = w(p, q) + b > \max\{w(p, q), b\} \ge \max\{w(p, q), a\} \ge |w(p, q - 1)|$, so the greedy algorithm assigns a $V$-student.

When $q$ is sufficiently small that $w(p - 1, q) \le 0$, then $w(p, q - 1) = w(p, q) - a = w(p - 1, q) - b - a < w(p - 1, q) \le 0$, so the greedy algorithm assigns a $U$-student.

When $w(p, q) < 0 < w(p - 1, q)$, then $w(p, q - 1) = w(p, q) - a < 0 < w(p - 1, q)$, so the greedy algorithm will assign a $U$-student when
\[
\begin{split}
  w(p - 1, q) < -w(p, q - 1)  &\Leftrightarrow aq - b(p - 1) < bp - a(q - 1)\\
  &\Leftrightarrow a(2q - 1) < b(2p - 1)\\
  &\Leftrightarrow q < \frac{b}{2a}(2p - 1) + \frac{1}{2}.
\end{split}
\]
In summary, in state $(p, q)$ the greedy algorithm assigns a $U$-student if
\[
  q \le \left \lceil \frac{b}{2a}(2p - 1) + \frac{1}{2} \right \rceil - 1 \equiv q^*(p),
\]
and a $V$-student otherwise. Note that $q^*(p + 1) \ge q^*(p)$.

The greedy path is the set of states $G = \{(p, q) : q^*(p) \le q \le q^*(p + 1)\}$. If $(p, q) \in G$ and $q > q^*(p)$, then a $V$-student will be assigned so the next state is $(p, q - 1)$ and this is on the greedy path. Alternatively, if $(p, q) \in G$ and $q = q^*(p)$, then a $U$-student will be assigned so the next state is $(p - 1, q)$ and this is again on the greedy path. Thus, once we reach a state on the greedy path, the greedy path will be followed for the remainder of the algorithm. Now
\[
\begin{split}
  q^*(a) &= \left \lceil \frac{b}{2a}(2a - 1) + \frac{1}{2} \right \rceil - 1 = \left \lceil b - \frac{b - a}{2a} \right \rceil - 1 \le b, \quad \quad \quad \text{and}\\
  q^*(a + 1) &= \left \lceil \frac{b}{2a}(2a + 1) + \frac{1}{2} \right \rceil = \left \lceil b + \frac{a + b}{2a} \right \rceil - 1 \ge b,
\end{split}
\]
so $q^*(a) \le b \le q^*(a + 1)$, so the initial state $(a, b) \in G$ and all future states will be on the greedy path.

On the greedy path, for a given value of $p$, the largest value state is $(p, q^*(p + 1))$. Here
\[
\begin{split}
  w(p, q^*(p + 1)) &= aq^*(p + 1) - bp\\
  &= a\left \lceil \frac{b}{2a}(2p + 1) + \frac{1}{2} \right \rceil - a - bp\\
  &< a\left( \frac{b}{2a}(2p + 1) + \frac{1}{2} \right) - bp\\
  &= b\left(p + \frac{1}{2}\right) + \frac{a}{2} - bp\\
  &= \frac{a + b}{2}.
\end{split}
\]
Similarly, for a given value of $p$, the smallest value state on the greedy path is $(p, q^*(p))$. Here
\[
\begin{split}
  w(p, q^*(p)) &= aq^*(p) - bp\\
  &= a\left \lceil \frac{b}{2a}(2p - 1) + \frac{1}{2} \right \rceil - a - bp\\
  &\ge a\left( \frac{b}{2a}(2p - 1) + \frac{1}{2} \right) - a - bp\\
  &= b\left(p - \frac{1}{2}\right) + \frac{a}{2} - a - bp\\
  &= -\frac{a + b}{2}.
\end{split}
\]
Conclude that for any state $(p, q) \in G$:
\[
  -\frac{a + b}{2} \le w(p, q) < \frac{a + b}{2}.
\]
There are $a + b$ unique integers in this interval and each such integer must be taken on by some $w_i$ because (just as we argued in the proof of the previous lemma) the first $a + b$ values of $w_i$ are distinct integers.

Because $a$ and $b$ are coprime, there are two cases. First, when one of $a$ and $b$ is even and the other is odd, we have
\[
  y_{\min}^{\prime} = - \frac{a + b - 1}{2} \quad \text{and} \quad z_{\min}^{\prime} = - y_{\max}^{\prime} = - \frac{a + b - 1}{2}.
\]
Then part 2 of the greedy algorithm make no changes so the algorithm produces a smallest piece with size $-(a + b - 1)/2$.

Second, when both $a$ and $b$ are odd, 
\[
  y_{\min}^{\prime} = - \frac{a + b}{2} \quad \text{and} \quad z_{\min}^{\prime} = - y_{\max}^{\prime} = - \frac{a + b - 2}{2}.
\]
Then in part 2 of the greedy algorithm $\varepsilon = 1/2$, so the pieces are adjusted so that once again the smallest piece has size $-(a + b - 1)/2$.

Now appeal to Lemma \ref{P'_opt_LB} to conclude that the greedy algorithm produces an optimal solution.\hspace{\stretch{1}}$\blacksquare$

\begin{corollary} \label{greedy_opt} For any 3M-DAP with $t=u=v=2$, the greedy algorithm produces an optimal solution. The optimal value is
\[
  \frac{x_t}{2} - \frac{x_t - x_u}{b} \frac{a + b - 1}{2} = \frac{x_t}{2} - \frac{x_v - x_u}{2} \frac{a + b - 1}{a + b} = \frac{x_t}{2} + \frac{\lambda}{2} \frac{a + b - 1}{a + b} = \frac{x_u}{2} - \frac{x_t - x_u}{2} \frac{a - 1}{b}.
\]
\end{corollary}

It is immediate from Corollary \ref{greedy_opt} that such a problem is a 0-problem if and only if $a = 1$. This of course matches what we have learned earlier. First, there are no 0-problems of type 1 because $x_{\infty} = x_u + x_t - x_v = \lambda + \gamma v$. The 0-problems of type 2 occur when
\begin{align*}
  x = x_u = x_b &= \frac{(x_u + x_t - x_v)b + x_t}{b + 1}\\
  \Leftrightarrow x_u &= (b + 1)x_t - bx_v\\
  \Leftrightarrow x_u + bx_v &= (b + 1)x_t,
\end{align*}
which on comparison with \eqref{abc} confirms the conclusion.

We also note that the optimal value
\[
  \frac{x_t}{2} + \frac{\lambda}{2} \frac{a + b - 1}{a + b} > \frac{x_t}{2} + \frac{\lambda}{2}
\]
is indeed greater than the strict lower bound established in Lemma \ref{greedy_t=u=v=2}.

\subsection{Proof of Lemma \ref{lem:lambdaLB}} \label{lem:lambdaLB_proof}

In this section we prove Lemma \ref{lem:lambdaLB}, which requires us to consider more general DAPs than the 3M-DAPs to which we can restrict attention otherwise. We begin by generalizing Lemma \ref{greedy_t=u=v=2}.

Suppose that we have $s$ sources, each with supply $x_t$ and each of which must be divided into two pieces. Each of $s$ sinks will be assigned two of these pieces. The demands at the sinks may differ: suppose that sink $j$ has demand $x_j$. We must have $\sum_{j=1}^s x_j = sx_t$.

Consider the following greedy algorithm for this type of problem. Construct two vectors $\boldsymbol{y}$ and $\boldsymbol{z}$, each of length $s$. Divide the first row of $T$, $[y_1, z_s] = [x_t/2, x_t/2]$. Divide the $i^{th}$ row of $T$, for $2 \le i \le s$, $[y_i, z_{i-1}]$. Each sink receives pieces $[y_j, z_j]$ for some $j$.

\fontsize{8}{8}
\begin{algorithm}
\SetKwInOut{Input}{input}
\SetKwInOut{Output}{output}
\Input{$T$, an empty matrix of dimensions $s \times 2$\\
$x_t$, the required row sum for each row of $T$\\
$J = \{1, \ldots, s\}$, the indices of the demand vectors\\
$\{\boldsymbol{u}_j : j \in J\}$, empty demand vectors, each containing 2 elements\\
$\{x_j : j \in J\}$, the required row sum for each demand vector\\
$sx_t = \sum_{j \in J} x_j$\\
}
\Output{the filled supply matrix $T$ and demand vectors $\{\boldsymbol{u}_j : j \in J\}$\\
the multiset of entries in $T$ is the union of the multiset of entries in\\
the $\{\boldsymbol{u}_j : j \in J\}$\\
}
\BlankLine

\Begin(\textbf{part 1}){
divide the first row of $T$, $[x_t/2, x_t/2]$, so that $y_1 = z_s = x_t/2$\\
\For{$1 \le i \le s$}{
	let $\{x^{\prime}_k\}$ represent the unique values among the remaining $\{x_j : j \in J\}$ arranged in ascending order: $x^{\prime}_1 < x^{\prime}_2 < \cdots < x^{\prime}_r$, where $r \le s + 1 - i$\\
	set $x^{\prime}_{r+1} = +\infty$\\
	set $k^* = \min\left\{k: 1 \le k \le r, y_i < \frac{1}{2}(x^{\prime}_k + x^{\prime}_{k + 1} - x_t) \right\}$\\
	assign $y_i$ to a sink $j^*$ requiring allocation $x^{\prime}_{k^*}$\\
	assign $z_i = x^{\prime}_{k^*} - y_i$ as the second piece to sink $j^*$\\
        	\uIf{$i < s$}{
		set $J \leftarrow J \backslash \{j^*\}$\\
        		compute $y_{i+1} = x_t - z_i$\\
        		divide the $(i + 1)^{st}$ row of $T$, $[y_{i+1}, z_i]$\\
        	}
}
}

\Begin(\textbf{part 2}){
compute $y_{\min} = \min_i\{y_i\}$ and $z_{\min} = \min_i\{z_i\}$\\
compute $\varepsilon = (z_{\min} - y_{\min})/2$\\
$y_i \leftarrow y_i + \varepsilon, \forall i$\\
$z_i \leftarrow z_i - \varepsilon, \forall i$\\
}

\caption{Greedy algorithm for DAP $(T; \{\boldsymbol{u}_j : j \in J\})$ where $T$ has dimensions $s \times 2$ and for each of the $s$ sinks $\boldsymbol{u}_j$ has 2 elements}
\label{alg:greedy_v2}
\end{algorithm}

\normalsize

At the $i^{th}$ step of the algorithm, we will have a piece $y_i$ that we must assign to one of the remaining sinks. If we assign $y_i$ to a sink $j$ with demand $x_j$, then the other piece to be assigned to that sink has size $z_i = z_i^j = x_j - y_i$. The greedy algorithm assigns $y_i$ so as to make $z_i$ as close to $x_t/2$ as possible.

The following lemma provides bounds on the sizes of the pieces in the solution given by the greedy algorithm. Write $x_{\max} = \max_j \{x_j\}$, $x_{\min} = \min_j \{x_j\}$, and $\lambda = x_{\min} - x_{\max}$.

\begin{lemma} \label{greedy_t=u=v=2_gen} If $\lambda < 0$, then Algorithm \ref{alg:greedy_v2} produces $\boldsymbol{y}$ and $\boldsymbol{z}$ with $(x_t + \lambda)/2 < y_i, z_i < (x_t - \lambda)/2$, for all $i$. 
\end{lemma}
\textbf{Proof.} We begin by showing that part 1 produces $\boldsymbol{y}$ and $\boldsymbol{z}$ with $(x_t + \lambda)/2 \le y_i < (x_t - \lambda)/2$ and $(x_t + \lambda)/2 < z_i \le (x_t - \lambda)/2$, for all $i$. 

At step $i$ of the algorithm there are two unassigned supply pieces ($y_i$ and $z_s = x_t/2$), $s - i$ undivided rows of $T$, and $s - i + 1$ unassigned sinks. Then
\[
  y_i + (s - i + 1/2)x_t = \sum_{j \in J} x_j.
\]

The proof is by induction. We certainly have $(x_t + \lambda)/2 \le y_1 < (x_t - \lambda)/2$. Suppose it is true of $y_i$. If $x^{\prime}_{k^*} \ge x_t$, then $z_i = x^{\prime}_{k^*} - y_i \ge x_t - y_i > (x_t + \lambda)/2$.

If $k^* = 1$, then $x_j \ge x_t$ for all the remaining sinks so it must be that $y_i \ge x_t$. Then
\[
  (s - i + 1/2)x_t - z_i = y_i + (s - i + 1/2)x_t - x^{\prime}_1 = \sum_{j \in J} x_j - x^{\prime}_1 \ge (s - i)x^{\prime}_1 \ge (s - i)x_t,
\]
so $z_i \le x_t/2 < (x_t - \lambda)/2$. Otherwise $k^* > 1$, so
\[
\begin{split}
  z_i &= x^{\prime}_{k^*} - y_i\\
  &\le x^{\prime}_{k^*} - \frac{1}{2}(x^{\prime}_{k^*-1} + x^{\prime}_{k^*} - x_t)\\
  &= \frac{1}{2}(x_t + x^{\prime}_{k^*} - x^{\prime}_{k^*-1})\\
  &\le \frac{1}{2}(x_t - \lambda).
\end{split}
\]
Alternatively, suppose $x^{\prime}_{k^*} < x_t$. Then $z_i = x^{\prime}_{k^*} - y_i < x_t - y_i \le (x_t - \lambda)/2$.

If $k^* = r$, then $x_j < x_t$ for all the remaining sinks so it must be that $y_i < x_t$. Then
\[
  (s - i + 1/2)x_t - z_i = y_i + (s - i + 1/2)x_t - x^{\prime}_r = \sum_{j \in J} x_j - x^{\prime}_r \le (s - i)x^{\prime}_r < (s - i)x_t,
\]
so $z_i > x_t/2 > (x_t + \lambda)/2$. Otherwise $k^* < r$, so
\[
\begin{split}
  z_i &= x^{\prime}_{k^*} - y_i\\
  &> x^{\prime}_{k^*} - \frac{1}{2}(x^{\prime}_{k^*} + x^{\prime}_{k^*+1} - x_t)\\
  &= \frac{1}{2}(x_t + x^{\prime}_{k^*} - x^{\prime}_{k^*+1})\\
  &\ge \frac{1}{2}(x_t + \lambda).
\end{split}
\]
Conclude that $(x_t + \lambda)/2 < z_i \le (x_t - \lambda)/2$. Then $(x_t + \lambda)/2 \le y_{i + 1} < (x_t - \lambda)/2$.

Finally, in part 2, the size of the smallest piece is only unchanged if $z_{\min} = y_{\min}$, but then both are larger than $(x_t + \lambda)/2$. Otherwise, the size of the smallest piece strictly increases and so must be larger than $(x_t + \lambda)/2$.\hspace{\stretch{1}}$\blacksquare$

\begin{lemma} \label{lem:alg_greedy_v2_complexity} When the number of distinct demands (required row sums) is a fixed constant, i.e., independent of $n_t$, Algorithm \ref{alg:greedy_v2} has complexity $\Theta(n_t)$.
\end{lemma}
\textbf{Proof.} The analysis is identical to the analysis of the complexity of Algorithm \ref{alg:greedy} except that the identification of $k^*$ requires more than just one comparison. But because the number of distinct demands is assumed fixed, this step still has complexity $O(1)$. This step is repeated $n_t/2$ times, so the conclusion that the complexity is linear in $n_t$ remains.
\hspace{\stretch{1}}$\blacksquare$

\fontsize{8}{8}
\begin{algorithm}
\SetKwInOut{Input}{input}
\SetKwInOut{Output}{output}
\Input{$T$, an empty matrix of dimensions $s_t \times t$, $t$ even\\
$x_t$, the required row sum for each row of $T$\\
$J = \{1, \ldots, s\}$, the indices of the demand vectors\\
$\{\boldsymbol{u}_j : j \in J\}$, empty demand vectors with lengths $\{u_j : j \in J\}$\\
$\{x_j : j \in J\}$, the required row sum for each demand vector\\
$ts_t = \sum_{j \in J} u_j$, $s_tx_t = \sum_{j \in J} x_j$\\
}
\Output{the filled matrix $T$ and demand vectors $\{\boldsymbol{u}_j : j \in J\}$\\
the multiset of entries in $T$ is the union of the multiset of entries in\\
the $\{\boldsymbol{u}_j : j \in J\}$\\
}
\BlankLine

divide each row of $T$ into $t/2$ pairs of elements\\
each pair is a row in matrix $T^{\prime}$: dimensions $(ts_t/2) \times 2$, required row sums $x_t^{\prime} = 2x_t/t$\\
divide $\frac{1}{2}\sum_{j \in J} (u_j - 2) = ts_t/2 - |J|$ rows of $T^{\prime}$ $[x_t/t, x_t/t]$\\
let $T^{\prime \prime}$ be the matrix consisting of the remaining rows of $T^{\prime}$\\
\For{$j \in J$}{
	set the first $u_j - 2$ elements of $\boldsymbol{u}_j$ to value $x_t/t$\\
	let $\boldsymbol{u}_j^{\prime}$ be the vector consisting of the last two entries of $\boldsymbol{u}_j$\\
	set $x_j^{\prime} = x_j - (u_j - 2)x_t/t$\\
}
use Algorithm \ref{alg:greedy_v2} to solve the $(T^{\prime \prime}; \{\boldsymbol{u}_j^{\prime} : j \in J\})$ DAP\\
\caption{Algorithm to achieve lower bound for DAP $(T; \{\boldsymbol{u}_j : j \in J\})$ with $t$ even}
\label{alg:t_even_gen}
\end{algorithm}

\normalsize

As before, we can leverage the lower bound on the optimal solution for any DAP $(T; \{\boldsymbol{u}_j : j \in J\})$ in which all supplies and demands are divided into two to establish useful lower bounds on more general problems. Algorithm \ref{alg:t_even_gen} shows how to do this for any DAP where the supplies can be represented by one matrix with an even number of columns.  In the algorithm, for each $j \in J$, the two elements of $\boldsymbol{u}_j^{\prime}$ must sum to $x_j - (u_j - 2)x_t/t = x_j - (u_j - 2)x_t^{\prime}/2$. Set 
\[
\begin{split}
  x_{\max}^{\prime} &= \max_{j \in J}\{x_j - (u_j - 2)x_t/t\}\\
  x_{\min}^{\prime} &= \min_{j \in J}\{x_j - (u_j - 2)x_t/t\}\\
  \lambda^{\prime} &= x_{\min}^{\prime} - x_{\max}^{\prime}.
\end{split}
\]
The following lemma provides bounds on the sizes of the pieces in the solution given by Algorithm \ref{alg:t_even_gen}. 

\begin{lemma} \label{lem:LB_on_f;t_even_gen} Given a DAP $P = (T; \{\boldsymbol{u}_j : j \in J\})$ with $t$ even, $f(P) > \lambda^{\prime}/2 + x_t/t$. When $t = 2$, in an optimal solution the largest element has size less than $(x_t - \lambda^{\prime})/2$.
\end{lemma}
\textbf{Proof.} Apply Algorithm \ref{alg:t_even_gen} to solve the $(T; \{\boldsymbol{u}_j : j \in J\})$ problem. In the final step, on using Algorithm \ref{alg:greedy_v2} to solve the $(T^{\prime \prime}; \{\boldsymbol{u}_j^{\prime} : j \in J\})$ problem, Lemma \ref{greedy_t=u=v=2_gen} applies, whence the smallest piece in an optimal solution to this reduced problem is greater than
\[
  \frac{1}{2}(x_t^{\prime \prime} + \lambda^{\prime}) = \frac{1}{2}\lambda^{\prime} + \frac{x_t}{t}.
\] 
\hspace{\stretch{1}}$\blacksquare$

\begin{lemma} \label{lem:alg_t_even_gen_complexity} When the number of distinct demands (required row sums) is a fixed constant, i.e., independent of $n_t$, Algorithm \ref{alg:t_even_gen} has complexity $\Theta(n_t)$.
\end{lemma}
\textbf{Proof.} The analysis is identical to the analysis of the complexity of Algorithm \ref{alg:t_even}. The reduced DAP constructed has the same number of distinct demands as the input to the algorithm, which by assumption is a fixed constant. Then Lemma \ref{lem:alg_greedy_v2_complexity} applies and the conclusion follows.
\hspace{\stretch{1}}$\blacksquare$

\fontsize{8}{8}
\begin{algorithm}
\SetKwInOut{Input}{input}
\SetKwInOut{Output}{output}
\Input{$T$, $U$, $V$, empty matrices of dimensions $s_t \times t$, $s_u \times u$, and $s_v \times v$, respectively\\
$x_t$, $x_u$, $x_v$, the required row sums for each row of $T$, $U$, and $V$, respectively\\
$u = v + 1$\\
}
\Output{the filled matrices $T$, $U$, $V$\\
the multiset of entries in $T$ is the union of the multiset of entries in $U$ and $V$\\
}
\BlankLine

\uIf{$t$ even}{
use Algorithm \ref{alg:t_even} to solve the $(T; U, V)$ 3M-DAP\\
}
\uElseIf{$s_t \ge (u - 2)s_u + (v - 2)s_v$}{
// $u = 3, v = 2, s_t \ge s_u$\\
\uIf{$n_u < s_t$}{
fill all elements of $U$ with value $x_u/u$ and insert these into the first column of $T$\\
form a $(T^{\prime}; U^{\prime}, V^{\prime})$ 3M-DAP where:\\
$\quad T^{\prime} = V$; $x_t^{\prime} = x_v$\\
$\quad U^{\prime}$ is the $(s_t - n_u) \times t$ unfilled submatrix of $T$; $x_u^{\prime} = x_t$\\
$\quad V^{\prime}$ is the $n_u \times (t - 1)$ unfilled submatrix of $T$; $x_v^{\prime} = x_t - x_u/u$\\
$t^{\prime} = 2$ and $u^{\prime} = v^{\prime} + 1$, so use Algorithm \ref{alg:t_even} to solve the $(T^{\prime}; U^{\prime}, V^{\prime})$ 3M-DAP\\
}
\uElseIf{$n_u = s_t$}{
fill all elements of $U$ with value $x_u/u$ and insert these into the first column of $T$\\
fill all elements of $V$ and all remaining elements of $T$ with value $x_v/v$\\
}
\uElse{
// $n_u > s_t$, $s_t - s_u$ even\\
set $r_u = (s_t - s_u)/2$, so $0 \le r_u < s_u$\\
fill all elements of the first column of $T$ with value $x_u/u$ and insert these into the the first $r_u$ rows of $U$ and the remainder of the first column of $U$\\
form a $(T^{\prime}; U^{\prime}, V^{\prime})$ 3M-DAP where:\\
$\quad T^{\prime}$ is the $s_t \times (t - 1)$ unfilled submatrix of $T$; $x_t^{\prime} = x_t - x_u/u$\\
$\quad U^{\prime}$ is the $(s_u - r_u) \times 2$ unfilled submatrix of $U$; $x_u^{\prime} = 2x_u/3$\\
$\quad V^{\prime} = V$; $x_v^{\prime} = x_v$\\
$t^{\prime}$ is even and $u^{\prime} = v^{\prime} + 1$, so use Algorithm \ref{alg:t_even} to solve the $(T^{\prime}; U^{\prime}, V^{\prime})$ 3M-DAP\\
}
}
\nonl continued below...\\
\caption{Algorithm to achieve strict lower bound for 3M-DAP $(T; U, V)$ with $u=v+1$}
\label{alg:u=v+1}
\end{algorithm}

\normalsize

\addtocounter{algocf}{-1}

\fontsize{8}{8}
\begin{algorithm*}
\setcounter{AlgoLine}{23}
\BlankLine

\nonl ...continued from above\\
\uElse{
Fill the first column of $T$ with value $\mu$ where $\lambda < \mu \le x_t/t$\\
Insert the $s_t$ values of $\mu$ into the columns of $U$ and $V$ in the following order: first column of $U$, first column of $V$, second column of $U$, second column of $V$, etc.\\
\uCase{\emph{\textbf{1}} All rows of $U$ and $V$ have the same number $w \ge 2$ of unfilled elements}{
form a $(T^{\prime}; U^{\prime}, V^{\prime})$ 3M-DAP where:\\
$\quad T^{\prime}$ is the $s_t \times (t - 1)$ unfilled submatrix of $T$; $x_t^{\prime} = x_t - \mu$\\
$\quad U^{\prime}$ is the $s_u \times w$ unfilled submatrix of $U$; $x_u^{\prime} = x_u - (u - w)\mu$\\
$\quad V^{\prime}$ is the $s_v \times w$ unfilled submatrix of $V$; $x_v^{\prime} = x_v - (v - w)\mu$\\
$t^{\prime}$ is even and $u^{\prime} = v^{\prime} + 1$, so use Algorithm \ref{alg:t_even} to solve the $(T^{\prime}; U^{\prime}, V^{\prime})$ 3M-DAP\\
}
\uCase{\emph{\textbf{2}} All rows of $V$ have $w \ge 2$ unfilled elements; some number $0 < r < s_u$ rows of $U$ have $w$ and the remainder have $w + 1$ unfilled elements}{
form a $(T^{\prime}; U_1^{\prime}, U_2^{\prime}, V^{\prime})$ DAP where:\\
$\quad T^{\prime}$ is the $s_t \times (t - 1)$ unfilled submatrix of $T$; $x_t^{\prime} = x_t - \mu$\\
$\quad U_1^{\prime}$ is the $r \times w$ unfilled submatrix of $U$; $x_{u_1}^{\prime} = x_u - (u - w)\mu$\\
$\quad U_2^{\prime}$ is the $(s_u - r) \times (w + 1)$ unfilled submatrix of $U$; $x_{u_2}^{\prime} = x_u - (u - w - 1)\mu$\\
$\quad V^{\prime}$ is the $s_v \times w$ unfilled submatrix of $V$; $x_v^{\prime} = x_v - (v - w)\mu$\\
$t^{\prime}$ is even so use Algorithm \ref{alg:t_even_gen} to solve the $(T^{\prime}; U_1^{\prime}, U_2^{\prime}, V^{\prime})$ DAP\\
}
\uCase{\emph{\textbf{3}} All rows of $U$ have $w + 1 > 2$ unfilled elements; some number $0 < r < s_v$ rows of $V$ have $w$ and the remainder have $w + 1$ unfilled elements}{
form a $(T^{\prime}; U^{\prime}, V_1^{\prime}, V_2^{\prime})$ DAP where:\\
$\quad T^{\prime}$ is the $s_t \times (t - 1)$ unfilled submatrix of $T$; $x_t^{\prime} = x_t - \mu$\\
$\quad U^{\prime}$ is the $s_u \times (w + 1)$ unfilled submatrix of $U$; $x_u^{\prime} = x_u - (v - w - 1)\mu$\\
$\quad V_1^{\prime}$ is the $r \times w$ unfilled submatrix of $V$; $x_{v_1}^{\prime} = x_v - (v - w)\mu$\\
$\quad V_2^{\prime}$ is the $(s_v - r) \times (w + 1)$ unfilled submatrix of $V$; $x_{v_2}^{\prime} = x_v - (v - w - 1)\mu$\\
$t^{\prime}$ is even so use Algorithm \ref{alg:t_even_gen} to solve the $(T^{\prime}; U^{\prime}, V_1^{\prime}, V_2^{\prime})$ DAP\\
}
}
\caption{Algorithm to achieve strict lower bound for 3M-DAP $(T; U, V)$ with $u=v+1$, CONT.}
\label{alg:u=v+1_part2}
\end{algorithm*}

\normalsize

Now we are in a position to prove Lemma \ref{lem:lambdaLB}. Algorithm \ref{alg:u=v+1} produces a solution to a 3M-DAP $(T; U, V)$ with $u = v + 1$ with smallest piece strictly larger than $\lambda$. We restate the lemma with proof here.

\begin{lemma} \label{lem:lambdaLB_repeat} Given a 3M-DAP $P = (T; U, V)$ with $u = v + 1$, $f(P) > \lambda = x_u - x_v$.
\end{lemma}
\textbf{Proof.} Apply Algorithm \ref{alg:u=v+1} to solve $P$. When $t$ is even, apply Lemma \ref{lem:LB_on_f;t_even} to conclude that
\[
  f(P) > \frac{1}{2}\lambda + \frac{1}{2}\frac{x_t}{t} > \lambda,
\]
because $x_t/t > x_u/u > (x_u - x_v)/(u - v) = \lambda$.

Now we investigate when $t$ is odd. First, suppose that $s_t \ge (u - 2)s_u + (v - 2)s_v$. If $v \ge 3$, then $u \ge 4$, so $(v - 2)s_v \ge n_v/3$ and $(u - 2)s_u \ge n_u/2$, so $(u - 2)s_u + (v - 2)s_v > (n_u + n_v)/3 = n_t/3$. But $t \ge 3$, so $s_t \le n_t/3$. But this contradicts our supposition so we must have $v = 2$ and $u = 3$. Then the supposition is that $s_t \ge s_u$.

If $n_u = s_t$, Algorithm \ref{alg:u=v+1} produces a solution with only two values $x_u/u$ and $x_v/v$. 

If $n_u < s_t$, we obtain $f(P) \ge \min\{x_u/u, f(P^{\prime})\}$, where the subproblem $P^{\prime} = (T^{\prime}; U^{\prime}, V^{\prime})$ has $t^{\prime} = 2$, $u^{\prime} = v^{\prime} + 1$, and $\lambda^{\prime} = x_u/u$. By the first part of the proof, $f(P^{\prime}) > x_u/u$, so the conclusion follows.

If $n_u > s_t$, we need that $s_t - s_u = n_t - s_u - (t - 1)s_t$ is even, which follows because both $n_t - s_u = n_u - s_u + n_v = 2(s_u + s_v)$ and $t - 1$ are even. Then we obtain $f(P) \ge \min\{x_u/u, f(P^{\prime})\}$, where the subproblem $P^{\prime} = (T^{\prime}; U^{\prime}, V^{\prime})$ has $t^{\prime}$ even, $u^{\prime} = v^{\prime} + 1$, and $\lambda^{\prime} = 2x_u/3 - x_v$. By Lemma \ref{lem:LB_on_f;t_even}:
\[
\begin{split}
  f(P^{\prime}) &> \frac{1}{2}\lambda^{\prime} + \frac{x_t^{\prime}}{t^{\prime}}\\
  &= \frac{x_u}{3} - \frac{x_v}{2} + \frac{x_t - x_u/3}{t - 1}\\
  &= \lambda - \frac{2x_u}{3} + \frac{x_v}{2} + \frac{x_t - x_u/3}{t - 1}\\
  &> \lambda,
\end{split}
\]
because $x_v/2 > x_u/3$ and $(x_t - x_u/3)/(t - 1) > x_t/t > x_u/3$. Then $f(P) > \lambda$, as claimed.

Second, suppose that $s_t < (u - 2)s_u + (v - 2)s_v$. Because $\mu$ is such that $\lambda < \mu \le x_t/t$, we have
\[
  \frac{x_t}{t} > \lambda \Leftrightarrow \frac{x_t - \lambda}{t - 1} > \frac{x_t}{t} \Rightarrow \frac{x_t - \lambda}{t - 1} > \mu.
\]

\textbf{Case 1.} We obtain $f(P) \ge \min\{\mu, f(P^{\prime})\}$, where the subproblem $P^{\prime} = (T^{\prime}; U^{\prime}, V^{\prime})$ has $t^{\prime}$ even and $\lambda^{\prime} = x_u - (u - w)\mu - x_v + (v - w)\mu = \lambda - (u - v)\mu = \lambda - \mu$. By Lemma \ref{lem:LB_on_f;t_even}:
\[
\begin{split}
  f(P^{\prime}) &> \frac{1}{2}\lambda^{\prime} + \frac{x_t^{\prime}}{t^{\prime}}\\
  &= \frac{1}{2}(\lambda - \mu) + \frac{x_t - \mu}{t - 1}\\
  &= \lambda + \frac{x_t - \mu}{t - 1} - \frac{1}{2}(\lambda + \mu)\\
  &> \lambda + \frac{x_t}{t} - \mu\\
  &\ge \lambda.
\end{split}
\]
Then $f(P) > \lambda$, as claimed.

\textbf{Case 2.} We obtain $f(P) \ge \min\{\mu, f(P^{\prime})\}$, where the subproblem $P^{\prime} = (T^{\prime}; U^{\prime}_1, U^{\prime}_2, V^{\prime})$ has $t^{\prime}$ even and three possible values for the $x_j^{\prime \prime} = x_j^{\prime} - (u_j^{\prime} - 2)x_t^{\prime}/t^{\prime}$:
\begin{align*}
  x_v^{\prime \prime} &= x_v - (v - w)\mu - (w - 2)(x_t - \mu)/(t - 1)\\
  x_{u_1}^{\prime \prime} &= x_u - (u - w)\mu - (w - 2)(x_t - \mu)/(t - 1)\\
  x_{u_2}^{\prime \prime} &= x_u - (u - w - 1)\mu - (w - 1)(x_t - \mu)/(t - 1).
\end{align*}
Then
\[
\begin{split}
  x_{u_1}^{\prime \prime} - x_v^{\prime \prime} &= x_u - (u - w)\mu - (w - 2)(x_t - \mu)/(t - 1) - x_v + (v - w)\mu + (w - 2)(x_t - \mu)/(t - 1)\\
  &= \lambda - \mu < 0, \quad \quad \text{and}\\
  x_{u_2}^{\prime \prime} - x_v^{\prime \prime} &= x_u - (u - w - 1)\mu - (w - 1)(x_t - \mu)/(t - 1)\\
  &\quad \quad \quad \quad \quad \quad \quad \quad \quad - x_v + (v - w)\mu + (w - 2)(x_t - \mu)/(t - 1)\\
  &= \lambda - (x_t - \mu)/(t - 1) \le \lambda - \mu = x_{u_1}^{\prime \prime} - x_v^{\prime \prime} < 0.
\end{split}
\]
Then $x_{\min}^{\prime \prime} = x_{u_2}^{\prime \prime}$, $x_{\max}^{\prime \prime} = x_v^{\prime \prime}$, and $\lambda^{\prime \prime} = x_{u_2}^{\prime \prime} - x_v^{\prime \prime} = \lambda - (x_t - \mu)/(t - 1)$. Then by Lemma \ref{lem:LB_on_f;t_even_gen}:
\[
  f(P^{\prime}) > \frac{1}{2}\lambda^{\prime \prime} + \frac{x_t^{\prime}}{t^{\prime}} = \frac{1}{2}\lambda + \frac{1}{2}\frac{x_t - \mu}{t - 1} > \lambda.
\]
Then $f(P) > \lambda$, as claimed.

\textbf{Case 3.} We obtain $f(P) \ge \min\{\mu, f(P^{\prime})\}$, where the subproblem $P^{\prime} = (T^{\prime}; U^{\prime}, V^{\prime}_1, V^{\prime}_2)$ has $t^{\prime}$ even and three possible values for the $x_j^{\prime \prime} = x_j^{\prime} - (u_j^{\prime} - 2)x_t^{\prime}/t^{\prime}$:
\begin{align*}
  x_u^{\prime \prime} &= x_u - (u - w - 1)\mu - (w - 1)(x_t - \mu)/(t - 1)\\
  x_{v_1}^{\prime \prime} &= x_v - (v - w)\mu - (w - 2)(x_t - \mu)/(t - 1)\\
  x_{v_2}^{\prime \prime} &= x_v - (v - w - 1)\mu - (w - 1)(x_t - \mu)/(t - 1).
\end{align*}
Then
\[
\begin{split}
  x_u^{\prime \prime} - x_{v_2}^{\prime \prime} &= x_u - (u - w - 1)\mu - (w - 1)(x_t - \mu)/(t - 1)\\
  &\quad \quad \quad \quad \quad \quad \quad \quad \quad - x_v + (v - w - 1)\mu + (w - 1)(x_t - \mu)/(t - 1)\\
  &= \lambda - \mu < 0, \quad \quad \text{and}\\
  x_u^{\prime \prime} - x_{v_1}^{\prime \prime} &= x_u - (u - w - 1)\mu - (w - 1)(x_t - \mu)/(t - 1)\\
  &\quad \quad \quad \quad \quad \quad \quad \quad \quad - x_v + (v - w)\mu + (w - 2)(x_t - \mu)/(t - 1)\\
  &= \lambda - (x_t - \mu)/(t - 1) \le \lambda - \mu = x_u^{\prime \prime} - x_{v_2}^{\prime \prime} < 0.
\end{split}
\]
Then $x_{\min}^{\prime \prime} = x_u^{\prime \prime}$, $x_{\max}^{\prime \prime} = x_{v_1}^{\prime \prime}$, and $\lambda^{\prime \prime} = x_u^{\prime \prime} - x_{v_1}^{\prime \prime} = \lambda - (x_t - \mu)/(t - 1)$. Then by Lemma \ref{lem:LB_on_f;t_even_gen}:
\[
  f(P^{\prime}) > \frac{1}{2}\lambda^{\prime \prime} + \frac{x_t^{\prime}}{t^{\prime}} = \frac{1}{2}\lambda + \frac{1}{2}\frac{x_t - \mu}{t - 1} > \lambda.
\]
Then $f(P) > \lambda$, as claimed.
\hspace{\stretch{1}}$\blacksquare$

\begin{lemma} \label{lem:alg_u=v+1_complexity} Algorithm \ref{alg:u=v+1} has complexity $\Theta(n_t)$.
\end{lemma}
\textbf{Proof.} Algorithm \ref{alg:u=v+1} operates in one of three ways, depending on the problem characteristics. In the first way, when $s_t = n_u \ge (u - 2)s_u + (v - 2)s_v$, the problem is solved directly by assigning the value $x_u/u$ to all elements of the first column of $T$ and to all elements of $U$, and assigning the value $x_v/v$ to all remaining elements of $T$ and to all elements of $V$. Clearly, this direct solution has complexity $\Theta(n_t)$.

The second method involves assigning a given value to all elements of the first column of either $U$ or $T$, and similarly assigning this value to elements in either $T$ or in $U$ and $V$, so that the remaining elements form a reduced 3M-DAP that can be solved in linear time using Algorithm \ref{alg:t_even}.

The third method involves assigning a given value to all elements of the first column of $T$ and similarly assigning this value to elements in $U$ and $V$, so that the remaining elements form a reduced DAP that can be solved using Algorithm \ref{alg:t_even_gen}. This reduced DAP has only three distinct demands so Lemma \ref{lem:alg_t_even_gen_complexity} applies: Algorithm \ref{alg:t_even_gen} runs in linear time.
\hspace{\stretch{1}}$\blacksquare$

\end{document}